# Contact Curves In Projective Space

Amorim, É[*] Vainsencher, I.[†]

2019

**Abstract**

The projective space $\mathbb{P}^3$ admits a contact structure arising from a non integrable distribution of planes determined by a symplectic form in $\mathbb{C}^4$. We are interested in the rational curves of degree d which are tangent to that contact distribution in $\mathbb{P}^3$: the contact curves. To explore its geometry, we construct the parameter space $\mathcal{L}_d$ using Kontsevich's stable maps stack, $\overline{\mathcal{M}}_{0,0}(\mathbb{P}^3, d)$. The intersection theory on this space allows us to define the virtual invariant $N_d$, associated with the number of degree d contact curves incident to $2d+1$ lines. Using Bott's localization formula for stable maps, we determine a general formula for $N_d$. We explicitly calculate it for contact curves up to degree 4 - confirming the known cases of contact lines and conics and introducing the new numbers for cubics and quartics.

## 1 Introduction

Consider the symplectic structure $(\mathbb{C}^4, \omega)$, where $\omega$ is a symplectic bilinear form over $\mathbb{C}^4$. The isomorphism $v \mapsto \omega(v, \cdot)$ defines a differential form over $\mathbb{P}^3 = \mathbb{P}(\mathbb{C}^4)$, also denoted by $\omega$, such that $\omega_p := \omega(p, \cdot)$ for all $p \in \mathbb{P}^3$. As $\omega_p(p) = 0$, the morphism $p \mapsto \mathrm{Nuc}\,\omega_p$ allow us to realize $\omega$ as a distribution of planes in $\mathbb{P}^3$. Also, $\omega \wedge (d\omega)^{\wedge 3}$ is non-degenerate. In these conditions, we have a **contact structure** in $\mathbb{P}^3$ (for general definitions, see [ADN01], §4). We are interested in curves that are tangent to this distribution:

**Definition 1.1.** *Let $\mathbb{P}^3$ be with contact structure $\omega$. We say $C$ is a* ***contact curve*** *if it is tangent to the planes of contact distribution at each its non-singular points:*

$$T_p C \subset \mathrm{Nuc}\,\omega_p,\ \forall p \in C \setminus \mathrm{Sing}\,C$$

In the symplectic context, a contact curve C as above defined is also a **involutive** and **legendrian** curve. These definitions are distinct on a general odd dimensional projective space with a contact structure, but for one dimensional varieties in $\mathbb{P}^3$ they all coincides.

[*]Departamento de Matemática, Federal Center for Technological Education of Minas Gerais - CEFET-MG. Av. Amazonas 7675, 30510-000, Belo Horizonte-MG, Brazil. *email*: eden.amorim@cefetmg.br

[†]Departamento de Matemática - Universidade Federal de Minas Gerais - UFMG, Av. Antônio Carlos, 6627, 31270-901, Belo Horizonte, Minas Gerais, Brazil. Caixa Postal 702. email: israel@mat.ufmg.br

The involutiveness on a contact structure is defined in [ADN01]. For to explore a legendrian description for these curves, see [Buc06]. In [LV11], the enumerative aspects of involutive curves are explored through the contact condition and it is the origin of this work.

In this paper, we construct a parameter space for space contact curves and calculate some of its invariants using Kontsevich's stable maps.

A rational degree $d$ $n$-marked stable map $f : (C; p_1, \cdots, p_n) \to \mathbb{P}^3$ is a morphism from a connected nodal curve $C$ with $n$ non-singular marked points to space $\mathbb{P}^3$, such that it has a finite number of automorphisms fixing marks. The moduli problem for stable maps give us the stack of rational degree $d$ $n$-marked stable maps over $\mathbb{P}^3$, $\overline{\mathcal{M}}_{0,n}(\mathbb{P}^3, d)$. It is a smooth algebraic stack and has a coarse moduli space (as scheme)

$$m : \overline{\mathcal{M}}_{0,n}(\mathbb{P}^3, d) \to \overline{M}_{0,n}(\mathbb{P}^3, d). \tag{1}$$

Roughly, the stack structure allows consider families of stable maps in $\mathbb{P}^3$ with its automorphisms. Although there is no universal family for this moduli, there are structural morphisms alike:

$$\begin{array}{ccc} \overline{M}_{0,1}(\mathbb{P}^3, d) & \xrightarrow{\nu} & \mathbb{P}^3 \\ \downarrow{\scriptstyle \pi} & & \\ \overline{M}_{0,0}(\mathbb{P}^3, d) & & \end{array} \tag{2}$$

with $\nu$ the evaluation on the mark and $\pi$ the forgetful map, that exclude the mark and stabilize the curve.

The moduli and stack structure for stable maps are well known and described, see [FP97].

As a smooth algebraic stack, there is a intersection theory for $\overline{\mathcal{M}}_{0,n}(\mathbb{P}^3, d)$ over rational numbers (see [Vis89], [AGV08] §2). We denote by $A_*(\overline{\mathcal{M}}_{0,n}(\mathbb{P}^3, d))_{\mathbb{Q}}$ the Chow group of stack over rational numbers. Although the moduli $\overline{M}_{0,n}(\mathbb{P}^3, d)$ isn't smooth in general, it has a well behaved Chow group $A_*(\overline{M}_{0,n}(\mathbb{P}^3, d))_{\mathbb{Q}}$ due the proper pushforward of moduli morphism $m$:

$$m_*[\mathcal{V}] = \frac{1}{a(\mathcal{V})}[V],$$

with $a(\mathcal{V})$ the order of automorphism group over the generic point of integrally closed substack $\mathcal{V}$ with moduli $V$.

The proper pushforward $m_*$ give isomorphisms between Chow groups of the stack and its moduli. Also, the bivariant intersection structure is defined for stack and moduli spaces, as the duality. The pullback $m^*$ give isomorphism between the bivariant Chow rings of the stack and its moduli, over rational numbers.

Also, there is well known equivariant intersection theory for stable maps (for construction, see [GP99], appendix C). In our case, let $T = (\mathbb{C}^*)^4$ be a torus acting over $\mathbb{P}^3$. There is an induced action over $\overline{\mathcal{M}}_{0,n}(\mathbb{P}^3, d)$ and a localization formula. Let $\mathcal{M}_\Gamma$ be the components of fixed point locus for the T-action. These componentes are identified by combinatorial data given by certain graphs $\Gamma$ (see [Kon95], §3.2 or [CK99] §9.2).

For $\alpha$ in localization of equivariant Chow ring, we have:

$$\int_{\overline{\mathcal{M}}_{0,n}(\mathbb{P}^3,d)} \alpha = \sum_{\Gamma} \frac{1}{a(\Gamma)} \int_{\mathcal{M}_\Gamma} \frac{\iota_\Gamma^* \alpha}{\mathrm{Euler}^T(\mathcal{N}_\Gamma^{virt})}, \tag{3}$$

with $a(\Gamma)$ is the order of the automorphism group of $\mathcal{M}_\Gamma$ and $\mathcal{N}_\Gamma^{virt}$ its virtual normal bundle.

The stable maps fixed by the T-action are those whose image covers the coordinate lines of $\mathbb{P}^3$. The graphs $\Gamma$ are colored weighted trees that describe the combinatorial structure for fixed point locus. Each vertex $v$ has a color $i_v \in \{0,1,2,3\}$ corresponding to one of the four base points of $\mathbb{P}^3$. Each edge $e$ has a weight $d_e$, such that $\sum d_e = d$, corresponding to the integer partitions of $d$ over the coordinate lines of $\mathbb{P}^3$. Each tree $\Gamma$ determined by these conditions give a component. For example, the figure 1 show the image of a fixed stable map and its tree $\Gamma$, representing one component.

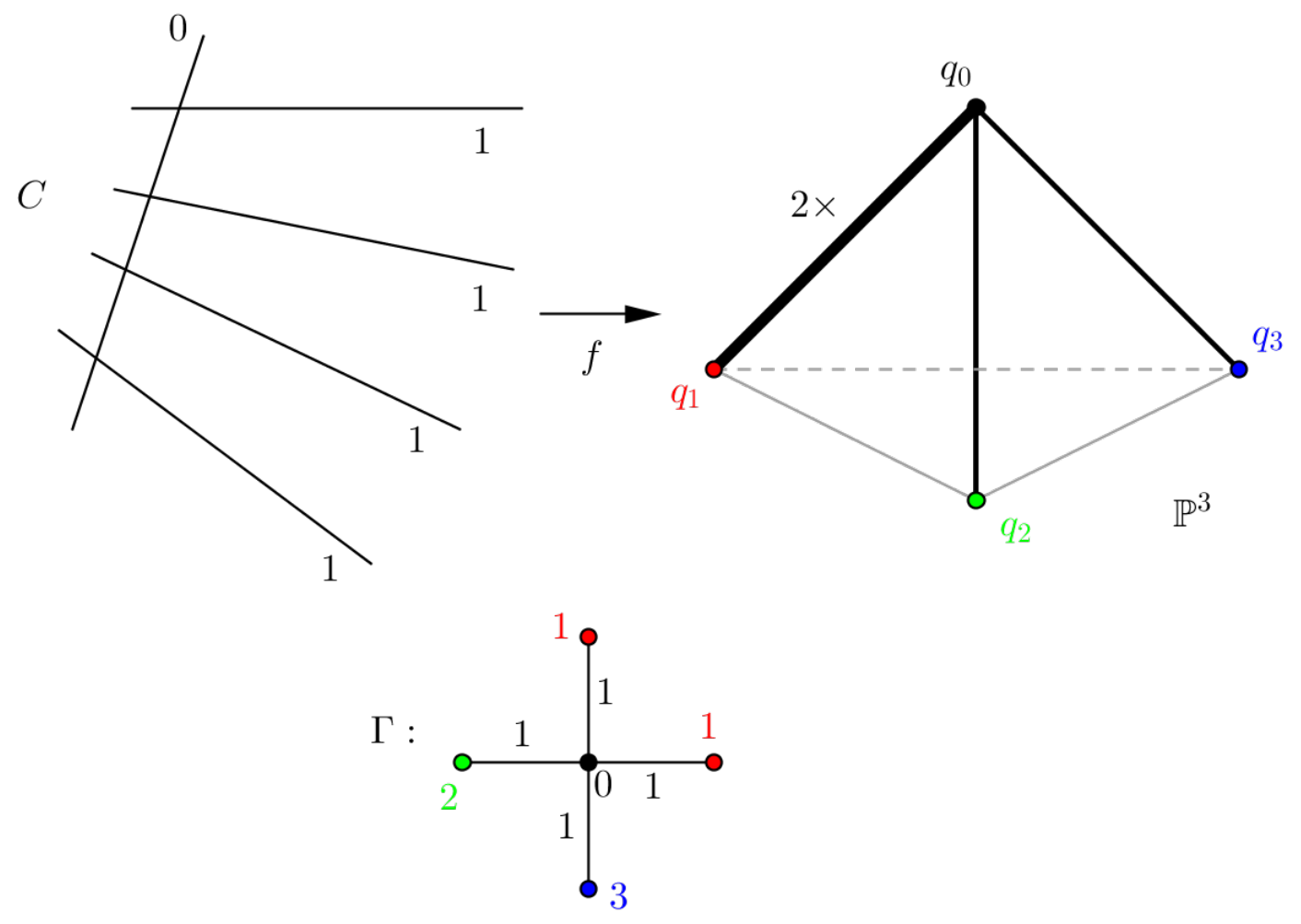

Fig. 1: Degree 4 fixed stable map on a positive dimensional component

If $\lambda_i$, $i = 1, \cdots, 4$ are the weights for the action $T = (\mathbb{C}^*)^4$ over $\mathbb{P}^3$, we express $1/\mathrm{Euler}^T(\mathcal{N}_\Gamma^{\mathrm{virt}}) = V(\Gamma) \cdot E(\Gamma)$ (see [Kon95], §3.3.3), where

$$V(\Gamma) = \prod_v \left( \left( \prod_{j \neq i_v} (\lambda_{i_v} - \lambda_j) \right)^{\mathrm{val}\, v - 1} \left( \sum_{e=(v,w)} \frac{d_e}{\lambda_{i_v} - \lambda_{i_w}} \right)^{\mathrm{val}\, v - 3} \prod_{e=(v,w)} \frac{d_e}{\lambda_{i_v} - \lambda_{i_w}} \right) \tag{4}$$

and

$$E(\Gamma) = \prod_{e=(v,w)} \left( \frac{(-1)^{d_e} \left(\frac{d_e}{\lambda_{i_v} - \lambda_{i_w}}\right)^{2d_e}}{d_e!^2} \prod_{k \neq i_v, i_w} \prod_{\alpha=0}^{d_e} \frac{1}{\frac{\alpha\lambda_{i_v} + (d_e - \alpha)\lambda_{i_w}}{d_e} - \lambda_k} \right). \tag{5}$$

## 2 A Space For Contact Curves

The contact strucutre in $\mathbb{P}^3$ is described, as in [Bar77], §7, by the exact sequence

$$\mathcal{D} \rightarrowtail T\mathbb{P}^3 \overset{\omega}{\twoheadrightarrow} \mathcal{O}_{\mathbb{P}^3}(2)\ . \tag{6}$$

Consider a family of stable maps, without marks,

$$\begin{array}{ccc} \mathcal{C} & \xrightarrow{\varphi} & \mathbb{P}^3 \\ \downarrow{\scriptstyle \pi} & & \\ B & & \end{array} \tag{7}$$

where the fiber by $\pi$ over $b \in B$ is a stable map $\varphi_b$. Let $\omega_\pi$ be the relative dualizing sheaf over $\mathcal{C}$ and consider the morphism $\omega_\pi^\vee \to \varphi^* T\mathbb{P}^3$. Over $b \in B$, this morphism corresponds to derivative map of $\varphi_b$, in the non-singular locus. Thus, by the diagram

$$\begin{array}{ccccc} \varphi^*\mathcal{D} & \rightarrowtail & \varphi^* T\mathbb{P}^3 & \twoheadrightarrow & \varphi^*\mathcal{O}_{\mathbb{P}^3}(2), \\ & \nwarrow & \uparrow & \nearrow_{s} & \\ & & \omega_\pi^\vee & & \end{array} \tag{8}$$

the contact condition for the stable map $\varphi_b$ corresponds to the vanish of section $s$. Considering the section

$$\mathcal{O} \to \omega_\pi \otimes \varphi^*\mathcal{O}_{\mathbb{P}^3}(2),$$

the contact curves locus at that family is the zero locus of the adjoint section, still denoted by $s$ (see [AK77] §2.3)

$$\mathcal{O} \to \pi_*(\omega_\pi \otimes \varphi^*\mathcal{O}_{\mathbb{P}^3}(2)).$$

The sheaf $\pi_*(\omega_\pi \otimes \varphi^*\mathcal{O}_{\mathbb{P}^3}(2))$ is in fact vector bundle, because $H^1(C, \omega_C \otimes f^*\mathcal{O}_{\mathbb{P}^3}(2)) = 0$. Its rank is $2d - 1$ (see [Amo14] §2.4).

Now, consider the stack family

$$\begin{array}{ccc} \overline{\mathcal{M}}_{0,1}(\mathbb{P}^3, d) & \xrightarrow{\nu} & \mathbb{P}^3 \\ \downarrow{\scriptstyle \pi} & & \\ \overline{\mathcal{M}}_{0,0}(\mathbb{P}^3, d) & & \end{array}$$

Define a vector bundle $\mathcal{E}_d$ over $\overline{\mathcal{M}}_{0,0}(\mathbb{P}^3, d)$ such that, for each family as (7), there is the bundle $\pi_*(\omega_\pi \otimes \varphi^*\mathcal{O}_{\mathbb{P}^3}(2))$, compatible with base change (for general definition, see [Vis89], §7). So, the fiber over a stable map $(C; f)$ is $H^0(C, \omega_C \otimes f^*\mathcal{O}_{\mathbb{P}^3}(2))$. We write $\mathcal{E}_d = \pi_*(\omega_\pi \otimes \nu^*\mathcal{O}_{\mathbb{P}^3}(2))$.

Let $\overline{s}$ be the section of $\mathcal{E}_d$ associated to section $s$. Define $\mathcal{L}_d = \mathcal{Z}(\overline{s})$, the zero locus of $\overline{s}$. Let $\mathcal{L}_d \to L_d$ be its associated moduli space. As substack, $\mathcal{L}_d$ is the functor such that, for each scheme $B$, associates the isomorphism class os stable map families over $B$ which the section $s$ in (8) vanishes along the fibers.

Now, we define the virtual fundamental class of $\mathcal{L}_d$. If $\iota : \mathcal{L}_d \hookrightarrow \overline{\mathcal{M}}_{0,0}(\mathbb{P}^3, d)$, so

$$\iota_*([\mathcal{L}_d]^{\text{virt}}) = c_{2d-1}(\mathcal{E}_d) \cap [\overline{\mathcal{M}}_{0,0}(\mathbb{P}^3, d)].$$

Since $\dim \overline{\mathcal{M}}_d = 4d$ and $\mathcal{E}_d$ has rank $2d - 1$, the expected dimension of $\mathcal{L}_d$ is $2d + 1$.

Thus, we have

**Theorem 2.1.** *The space of rational degree* $d$ *contact curves in* $\mathbb{P}^3$ *is parametrized by the* stack $\iota : \mathcal{L}_d \hookrightarrow \overline{\mathcal{M}}_{0,0}(\mathbb{P}^3, d)$, *the zero locus of a section of vector bundle* $\mathcal{E}_d = \pi_*(\omega_\pi \otimes \nu^*\mathcal{O}_{\mathbb{P}^3}(2))$. *In particular,* $\mathcal{L}_d$ *has expected dimension* $2d + 1$ *and its virtual fundamental class is such that*

$$\iota_*([\mathcal{L}_d]^{\text{virt}}) = c_{2d-1}(\mathcal{E}_d) \cap [\overline{\mathcal{M}}_{0,0}(\mathbb{P}^3, d)] \in A_{2d+1}(\overline{\mathcal{M}}_{0,0}(\mathbb{P}^3, d))_{\mathbb{Q}}.$$

Finally, we define a virtual invariant related to contact curves:

**Definition 2.1.** *For each positive integer* $d$, *we define the line incidence contact virtual invariant by*

$$N_d := \int_{\mathcal{L}_d} H^{2d+1} = \int_{\overline{\mathcal{M}}_d} c_{2d-1}(\mathcal{E}_d) \cup H^{2d+1}, \tag{9}$$

*where* $H$ *is the divisor for incidence of a stable map to a line in general position in* $\mathbb{P}^3$.

Well-known types of invariants over stable maps are in [Kon95], §2. The first type are Gromov-Witten invariants, such as the number of rational plane curves of a fixed degree passing through enough generic points. These numbers are defined and calculated using only incidence conditions. The second type comes from Mirror Symmetry, such as the number of rational curves on a quintic 3-fold in $\mathbb{P}^4$, or invariants over multiple coverings of rational curves on Calabi-Yau 3-folds (see also [CK99] §7.1, 9.2). In this case, these invariantes are defined by Euler classes of related sheaves.

Therefore, the invariant $N_d$ (9) is a third case, defined by cohomologic union of incidence conditions and Euler classes of sheaves: while the Euler class $c_{2d-1}(\mathcal{E}_d)$ is related to contact condition, the line incidence class $H$ is used to complete to correct codimension.

In this way, virtually, $N_d$ is an invariant related to contact space curves of a fixed degree passing through enough generic lines. In the cases which we confirm the enumerative significance for $N_d$, it counts the number of rational degree $d$ contact curves in $\mathbb{P}^3$, parametrizing by stable maps, that are incident to $2d+1$ lines in general position.

# 3 Calculating Contact Invariants

For calculate $N_d$, we apply the localization formula for stable maps (3), obtaining

$$\int_{\overline{\mathcal{M}}_d} c_{2d-1}(\mathcal{E}_d) \cup H^{2d+1} = \sum_{\Gamma} \frac{1}{a(\Gamma)} \int_{\mathcal{M}_\Gamma} \frac{\iota_\Gamma^*(c_{2d-1}^T(\mathcal{E}_d) \cup H_T^{2d+1})}{\mathrm{Euler}^T(\mathcal{N}_\Gamma^{\mathrm{virt}})}.$$

More explicity:

**Theorem 3.1.** *If $\lambda_i$, $i = 1, \cdots, 4$, are the weights for the action $T = (\mathbb{C}^*)^4$ over $\mathbb{P}^3$, we have*

$$N_d = \sum_{\Gamma} \frac{c_{2d-1}^T(\mathcal{E}_d)(\Gamma) \cdot H_T(\Gamma)^{2d+1}}{a(\Gamma)\, \mathrm{Euler}^T(\mathcal{N}_\Gamma^{\mathrm{virt}})},$$

*in which*

$$H_T(\Gamma) := \sum_{\nu \in V} \mu(\nu)\lambda_i \tag{10}$$

*and*

$$c_{2d-1}^T(\mathcal{E}_d)(\Gamma) := \prod_{e \in E} \prod_{\alpha=1}^{2d_e-1} \frac{\alpha\lambda_i + (d_e - \alpha)\lambda_j}{d_e} \cdot \prod_{\nu \in V} (2\lambda_{i_\nu})^{\mathrm{val}\,\nu - 1}, \tag{11}$$

*where $\mathrm{val}\,\nu$ is the number of adjacent edges to $\nu$ and $\mu(\nu)$ is the sum of weights of adjacent edges to $\nu$.*

The expression for $N_d$ follow from localization formula and $\mathrm{Euler}^T(\mathcal{N}_\Gamma^{\mathrm{virt}})$ is calculated by (4) and (5). It remains verify the expressions (10) and (11).

*Proof:* First, we fix general notations. Let $(z_0 : \cdots : z_3)$ be homogeneous coordinates for $\mathbb{P}^3$ and $q_i$ the base points. Let $\lambda_i$ be the weights of T-action over $\mathbb{P}^3$, in correspondence to each $z_i$.

Now, let $f : C \to \mathbb{P}^3$ be a (representative of) rational degree $d$ stable map of $\overline{M}_{0,0}(\mathbb{P}^3, d)$. Let $C = C' \cup \mathbb{P}^1$ be a decomposition for $C$, where $C'$ is a connected tree of projective lines, $\mathbb{P}^1$ a twig of $C$. Let $N$ be the intersection node between $C'$ and $\mathbb{P}^1$. Let $d_1$ be the degree of $f|_{\mathbb{P}^1}$ and $d' = d - d_1$ the degree of $f|_{C'}$. Let $(s : t)$ be homogeneous coordinates for $\mathbb{P}^1$, with $V(s) \mapsto N$.

For equation (10), we consider the formula

$$H = \pi_*(h^2) = c_1(\pi_*(\mathcal{O}(2h))) - 2c_1(\pi_*(\mathcal{O}(h))), \tag{12}$$

obtained by application of GRR theorem and techniques for treating exact sequences at nodal curves, based in [HM98], §3.E. (see [Amo14] §3.1.1 for complete calculations). Now, we proceed by mathematical induction over the number of components in C.

For base case, consider $C = \mathbb{P}^1$ and f a degree d map. If the image of f is the coordinate line $q_i q_j$, then $z_i = s^d$ e $z_j = t^d$ and the weights of s and t are $\lambda_i/d$ and $\lambda_j/d$. Thus, for $H^0(\mathbb{P}^1, \mathcal{O}_{\mathbb{P}^1}(kd))$, with generators

$$s^\alpha t^{kd-\alpha}\ ,\ 0 \leq \alpha \leq kd,$$

the weights are

$$\frac{\alpha\lambda_i + (kd - \alpha)\lambda_j}{d}.$$

By (12), using $k = 1$ and 2, we obtain $H_T(\Gamma) = d(\lambda_i + \lambda_j)$.

For inductive step, consider the exact sequence for nodal curves (extending a sheaf by zero from inclusion $C' \hookrightarrow C$, [Har77] §II.1, exercise 19(c))

$$\mathcal{O}_{\mathbb{P}^1}(kd_1 - N) \rightarrowtail \mathcal{O}_C(kh) \twoheadrightarrow \mathcal{O}_{C'}(kh).$$

As $H^1(\mathbb{P}^1, \mathcal{O}_{\mathbb{P}^1}(kd_1 - N)) = 0$, by Serre's duality, we have the cohomological exact sequence

$$H^0(\mathbb{P}^1, \mathcal{O}_{\mathbb{P}^1}(kd_1 - N)) \rightarrowtail H^0(C, \mathcal{O}_C(kh)) \twoheadrightarrow H^0(C', \mathcal{O}_{C'}(kh)).$$

A basis for for $H^0(\mathbb{P}^1, \mathcal{O}_{\mathbb{P}^1}(kd_1 - N))$ is

$$s^\alpha t^{kd_1-\alpha},\ 1 \leq \alpha \leq kd_1$$

and the associated weights are

$$\frac{\alpha\lambda_i + (kd_1 - \alpha)\lambda_j}{d_1},\ 1 \leq \alpha \leq kd_1.$$

Thus, for $k = 1$ we obtain $c_1(\pi_*\mathcal{O}_{\mathbb{P}^1}(d_1 - N))$ and

$$\sum_{\alpha=1}^{d_1} \frac{\alpha\lambda_i + (d_1 - \alpha)\lambda_j}{d_1} = \frac{d_1 + 1}{2}\lambda_i + \frac{d_1 - 1}{2}\lambda_j.$$

Similarly, the weights for $c_1(\pi_*\mathcal{O}_{\mathbb{P}^1}(2d_1 - N))$, with $k = 2$, given by

$$\sum_{\alpha=1}^{2d_1} \frac{\alpha\lambda_i + (2d_1 - \alpha)\lambda_j}{d_1} = (2d + 1)\lambda_i + (2d - 1)\lambda_j.$$

So, the weight for $H_T(\Gamma)$ is

$$\begin{aligned}
H &= \sum \text{weight of } H^0(\mathcal{O}_C(2h)) - 2\sum \text{weight of } H^0(\mathcal{O}_C(h)) = \\
&= \Big[\sum \text{weight of } H^0(\mathcal{O}_{C'}(2h)) + \sum \text{weight of } H^0(\mathcal{O}_{\mathbb{P}^1}(2d_1 - N))\Big] + \\
&\quad -2\Big[\sum \text{weight of } H^0(\mathcal{O}_{C'}(h)) + \sum \text{weight of } H^0(\mathcal{O}_{\mathbb{P}^1}(d_1 - N))\Big] = \\
&= \Big[\sum \text{weight of } H^0(\mathcal{O}_{C'}(2h)) - 2\sum (\text{weight of } H^0(\mathcal{O}_{C'}(h))\Big] + \\
&\quad +\Big[\sum \text{weight of } H^0(\mathcal{O}_{\mathbb{P}^1}(2d_1 - N)) - 2\sum \text{weight of } H^0(\mathcal{O}_{\mathbb{P}^1}(d_1 - N))\Big] = \\
&= \sum_{\nu\in\Gamma_{C'}} \mu(\nu)\lambda_{i_\nu} + \Big[(2d_1+1)\lambda_i + (2d_1-1)\lambda_j - (d_1+1)\lambda_i - (d_1-1)\lambda_j\Big] = \\
&= \sum_{\nu\in\Gamma_{C'}} \mu(\nu)\lambda_{i_\nu} + [d_1\lambda_i + d_1\lambda_j] = \sum_{\nu\in\Gamma_C} \mu(\nu)\lambda_{i_\nu},
\end{aligned}$$

proving the formula (10).

For equation (11), write $\omega_C(k) := \omega_C \otimes f^*\mathcal{O}_{\mathbb{P}^3}(k)$, for simplicity. Again, we proceed by mathematical induction over the number of components in C.

For base case, consider the Euler exact sequence

$$\omega_{\mathbb{P}^1}(2d) \rightarrowtail \mathcal{O}_{\mathbb{P}^1}(2d-1) \otimes H^0(\mathbb{P}^1, \mathcal{O}_{\mathbb{P}^1}(1)) \twoheadrightarrow \mathcal{O}_{\mathbb{P}^1}(2d).$$

and its cohomological exact sequence

$$H^0(\mathbb{P}^1, \omega_{\mathbb{P}^1}(2d)) \rightarrowtail H^0(\mathbb{P}^1, \mathcal{O}_{\mathbb{P}^1}(2d-1)) \otimes H^0(\mathbb{P}^1, \mathcal{O}_{\mathbb{P}^1}(1)) \twoheadrightarrow H^0(\mathbb{P}^1, \mathcal{O}_{\mathbb{P}^1}(2d)).$$

The generators for right space are $s^\alpha t^{2d-\alpha}$, $0 \le \alpha \le 2d$ and for center space are $s^\alpha t^{2d-1-\alpha} \otimes s, s^\alpha t^{2d-1-\alpha} \otimes t$, $0 \le \alpha \le 2d-1$. So, the generators for left space are

$$s^\alpha t^{2d-1-\alpha} \otimes t\ ,\ 1 \le \alpha \le 2d-1$$

and its weights are

$$\frac{\alpha\lambda_i + (2d-\alpha)\lambda_j}{d}\ ,\ 1 \le \alpha \le 2d-1.$$

For inductive step, consider the excision exact sequence

$$\omega_{C'}(2) \rightarrowtail \omega_C(2) \twoheadrightarrow \omega_{\mathbb{P}^1}(2d_1 + N).$$

and its cohomological sequence

$$H^0(C', \omega_{C'}(2)) \rightarrowtail H^0(C, \omega_C(2)) \twoheadrightarrow H^0(\mathbb{P}^1, \omega_{\mathbb{P}^1}(2d_1 + N)),$$

exact by vanishing of $H^1(C, \omega_C(2))$ for any C tree of projective lines (see [Amo14] §2.4).

To calculate the weights of right space, consider the Euler sequence

$$\omega_{\mathbb{P}^1}(2d_1 + N) \rightarrowtail \mathcal{O}_{\mathbb{P}^1}(2d_1 - 1 + N) \otimes H^0(\mathbb{P}^1, \mathcal{O}_{\mathbb{P}^1}(1)) \twoheadrightarrow \mathcal{O}_{\mathbb{P}^1}(2d_1 + N).$$

and its cohomological sequence

$$H^0(\mathbb{P}^1, \omega_{\mathbb{P}^1}(2d_1 + N)) \rightarrowtail H^0(\mathbb{P}^1, \mathcal{O}_{\mathbb{P}^1}(2d_1 - 1 + N)) \otimes H^0(\mathbb{P}^1, \mathcal{O}_{\mathbb{P}^1}(1)) \twoheadrightarrow H^0(\mathbb{P}^1, \mathcal{O}_{\mathbb{P}^1}(2d_1 + N)).$$

As this space contains the degree $2d_1$ rational functions with simple pole at N, its generators are

$$\frac{s^{\alpha}t^{2d_1+1-\alpha}}{s},\ 0 \le \alpha \le 2d_1+1.$$

Also, the generators for middle space are

$$\frac{s^{\alpha}t^{2d_1-\alpha}}{s} \otimes s, \frac{s^{\alpha}t^{2d_1-\alpha}}{s} \otimes t,\ 0 \le \alpha \le 2d_1.$$

Thus, the generators for left space are

$$\frac{s^{\alpha}t^{2d_1-\alpha}}{s} \otimes s,\ 0 \le \alpha \le 2d_1-1$$

and their weights

$$\frac{\alpha\lambda_i + (2d_1-\alpha)\lambda_j}{d_1},\ 0 \le \alpha \le 2d_1-1.$$

For $\alpha = 0$, we have the weight $2\lambda_j$ associated to the vertex for N. For $\alpha \neq 0$ the weights given by degree $2d_1$ monomials in s and t, associated to the edge for $\mathbb{P}^1$.

Now, the weight for Chern class is obtained by

$$c^T_{2d-1}(\mathcal{E}_d) = c^T_{top}(\omega_{C'}(2)) \cdot c^T_{top}(\omega_{\mathbb{P}^1}(2d_1+N)).$$

The first factor is given by induction hypothesis and the second was calculated above. ■

## 3.1 Some virtual values

The invariants $N_1$ and $N_2$, for contact lines and conics, are known as particular case of involutive plane curves, by [LV11]. For $d = 1, 2$, the theorem (3.1) results

$$N_1 = \int_{\overline{\mathcal{M}}_1} c_1(\mathcal{E}_1) \cdot H^3 = 2 \qquad \text{and} \qquad N_2 = \int_{\overline{\mathcal{M}}_2} c_3(\mathcal{E}_2) \cdot H^5 = 40.$$

which coincide with previous known values.

For contact cubics and quartics, we have

$$N_3 = \int_{\overline{\mathcal{M}}_3} c_5(\mathcal{E}_3) \cdot H^7 = 4160 \qquad \text{and} \qquad N_4 = \int_{\overline{\mathcal{M}}_4} c_7(\mathcal{E}_4) \cdot H^9 = 1089024.$$

These results are the first relevant values for the contact invariant. For contact curves with degree $d > 4$, the numbers $N_d$ are calculated by theorem (3.1), but not efficiently. The direct and ingenuous algorithms for calculation $N_d$ depend, at first, on the combinatorial construction and enumeration of graph data, its isomorphism classes and automorphisms.

## 3.2 Enumerative questions

By [LV11], the numbers $N_1$ and $N_2$ are exactly the number of contact lines and conics meeting 3 and 5 lines, respectively. In this work, the number are directly calculated by Schubert calculus and combinatorial counting.

The enumerative interpretation for numbers $N_3$ and $N_4$ remains incomplete at present work. Some general informations about degree d contact curves come from the following proprieties:

- There are irreducible contact curves for each degree d ([Amo14], theorem 2.3.1, based in [Buc06], § 3.1).

- The multiple covers don't contribute in $N_d$, by dimensional counting ([Amo14], theorem 2.5.1, based in [Bry00], Proposition 5.1).

For contact cubic curves, the contact cuspidal cubic (ramified) curves are plane curves, so don't contribute in $N_3$. The number of reducible contact cubics, at each boundary component of $M_{0,0}(\mathbb{P}^3, 3)$, is calculated by combinatorial configuration of tree of projective lines for stable maps (see [Amo14], §3.3). These combinatorial configuration are listed at figure 3.2, and it's considered the number of marked points at each twig and the number of automorphisms. There are 560 type (a) curves, 840 for type (b) and 1680 for type (c). Thus, if there isn't multiplicity at these boundary components, it's expected 1080 irreducible contact cubic curves. So far, the multiplicities are not calculated yet.

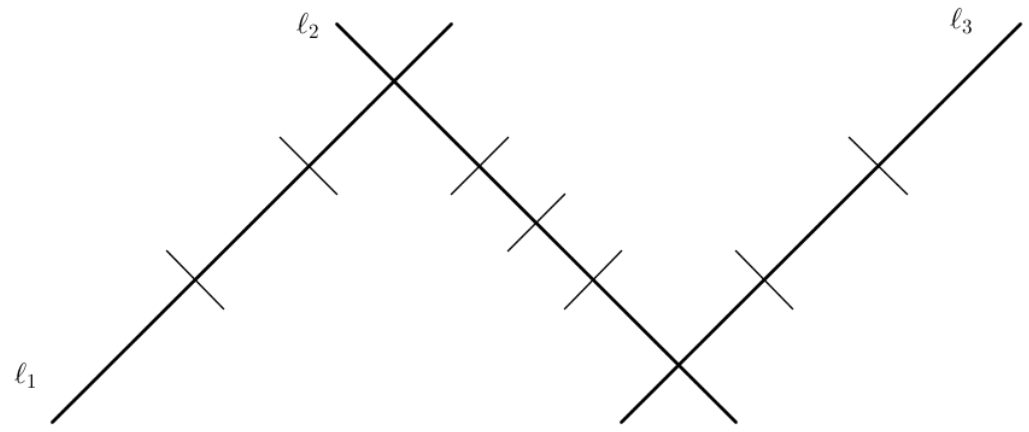

(a) $((2+3+2))$, Aut = 2

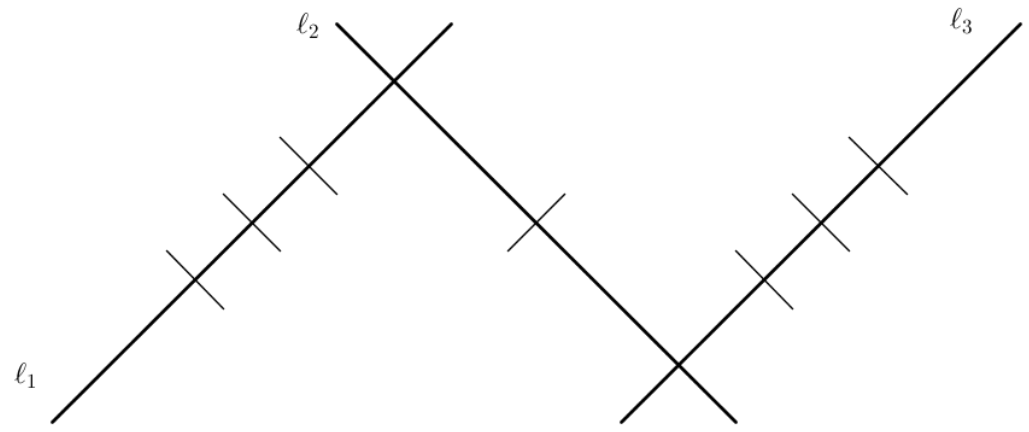

(b) $((3+1+3))$, Aut = 2

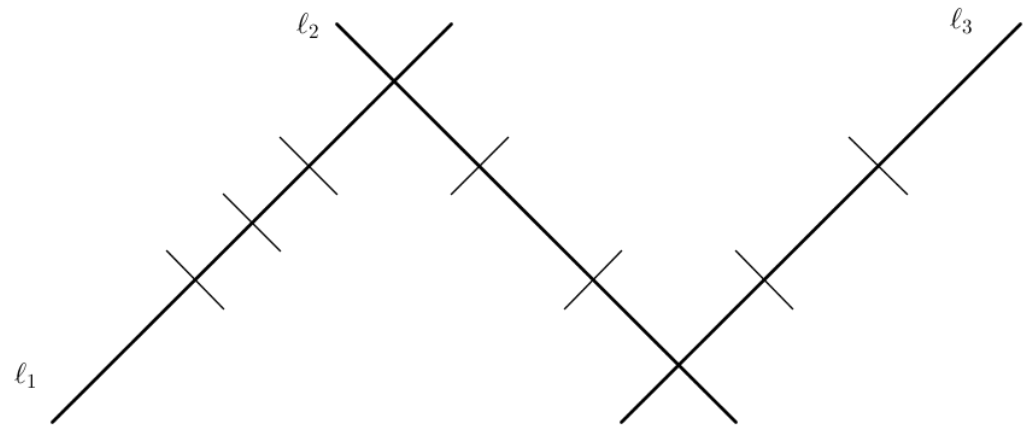

(c) $((3+2+2))$, Aut = 1

Fig. 2: Combinatorial configurations for contact tree of projective lines at boundary of $M_{0,0}(\mathbb{P}^3, 3)$
.

For contact quartic curves, the numbers of reducible contact quartics are calculated in analogue way, by combinatorial configurations of contact tree of projective lines at boundary of $M_{0,0}(\mathbb{P}^3, 4)$ (see [Amo14], §3.3). If there isn't multiplicity at these boundary components, the total of reducible contact quartic curves is 710080 and so, it's expected 378944 irreducible contact quartic curves.

As future work, the objective is complete these open questions using generators for Chow groups of stable maps ([Pan99]) and known tools in quantum cohomology ([FP97]).

Éden Santana Campos Amorim[1]

# Curvas de contato no espaço projetivo

*Contact Curves In Projective Space*

Tese submetida ao Corpo Docente do Programa de Pós-Graduação do Departamento de Matemática da Universidade Federal de Minas Gerais, como parte dos requisitos necessários para a obtenção do grau de Doutor em Matemática.
*PhD thesis submitted to Programa de Pós-Graduação of Departamento de Matemática da Universidade Federal de Minas Gerais.*

Orientador/Advisor: Israel Vainsencher[2]

Universidade Federal de Minas Gerais - UFMG
2014

[1]Departamento de Matemática - Centro Federal de Educação Tecnológica de Minas Gerais - CEFET-MG, Av Amazonas, 7675, 30510-000, Belo Horizonte, Minas Gerais, Brazil. email:eden.amorim@cefetmg.br

[2]Departamento de Matemática - Universidade Federal de Minas Gerais - UFMG, Av. Antônio Carlos, 6627, 31270-901, Belo Horizonte, Minas Gerais, Brazil. Caixa Postal 702. email: israel@mat.ufmg.br

## Resumo

O espaço projetivo de dimensão ímpar $\mathbb{P}^{2n-1}$ admite uma estrutura de contato proveniente da distribuição (não integrável) de hiperplanos determinada por uma forma simplética em $\mathbb{C}^{2n}$. Um objeto de interesse e tema de nosso estudo é o conjunto das curvas racionais de grau $d$ tangentes aos planos da distribuição de contato em $\mathbb{P}^3$. Tais curvas são chamadas curvas de contato, ou legendrianas. Para explorar a geometria de tais curvas, construimos o espaço de parâmetros $\mathcal{L}_d$ usando os mapas estáveis de Kontsevich com a estrutura do *stack* algébrico $\overline{\mathcal{M}}_{0,0}(\mathbb{P}^3, d)$. A teoria de interseção para *stacks* nos permite definir nesse espaço um invariante virtual $N_d$, associado ao número de curvas de contato de grau $d$ incidentes a $2d+1$ retas. Através de uma combinatória de grafos e partições oriunda da fórmula de localização de Bott determinamos uma fórmula geral para $N_d$. Calculamos explicitamente para curvas de grau até 4 - confirmando os casos já conhecidos de retas e cônicas de contato e apresentando os novos números associados a cúbicas e quárticas. Por fim, discutimos o significado enumerativo desses invariantes, ainda conjectural.

## Abstract

The odd dimensional projective space $\mathbb{P}^{2n-1}$ admits a contact structure arising from a non integrable distribution of hyperplanes determined by a symplectic form in $\mathbb{C}^{2n}$. Our object of interest is the set of rational curves of degree d which are tangent to that contact distribution in $\mathbb{P}^3$. Such curves are called contact curves or legendrian curves. To explore the geometry of contact curves, we construct the parameter space $\mathcal{L}_d$ using Kontsevich's stable maps, $\overline{\mathcal{M}}_{0,0}(\mathbb{P}^3, d)$, endowed with the structure of algebraic stack. The intersection theory on stacks allows us to define in that space the virtual invariant $N_d$, associated with the number of degree d contact curves incident to $2d+1$ lines. Using graph combinatorics and partitions originated from Bott's localization formula, we determine a general formula for $N_d$. We explicitly calculate it for contact curves up to degree 4 - confirming the known cases of contact lines and conics and introducing the new numbers for cubics and quartics. Finally, we discuss the enumerative significance of these invariants, still conjectural for $d > 4$.

# Sumário

# Lista de Figuras

# Capítulo 1

# Mapas estáveis

Em poucas palavras, os mapas estáveis fornecem um modo de estudar curvas através de suas parametrizações e estão na base de vários problemas enumerativos, incluindo os invariantes de Gromov-Witten, e teorias da física matemática. Introduzidos por Maxim Kontsevich, esses elementos e seu espaço de parâmetros guardam de modo geral mais informações que a estrutura de esquema algébrico suporta - como automorfismos dos elementos - e naturalmente exigem estruturas categóricas ou funtoriais, como os *stacks*. Nesta seção apresentaremos as definições e propriedades básicas dos espaços de mapas estáveis, mas suficientes como ferramentas para descrever e resolver problemas no estudo das curvas de contato no espaço projetivo. Em particular, estaremos restritos a mapas estáveis parametrizando curvas racionais em espaços projetivos $\mathbb{P}^r$, o caso mais bem descrito e explorado.

As seções estão estruturadas de modo a introduzir o conceito de estabilidade (1.1), o problema de *moduli* associado (1.2), assim como a extensão para o *stack* de mapas estáveis (1.3). Por fim, baseado na teoria de interseção para *stacks* e no trabalho de Kontsevich, introduzimos a fórmula de localização no anel de Chow equivariante do *stack* de mapas estáveis e descrevemos os mapas estáveis fixos pela ação de um toro (1.4).

A referência canônica usada é [16], ou também o texto mais introdutório [23]. Definições e propriedades gerais de um *stack* algébrico são deixadas para o apêndice (A). Para a fórmula de localização de Bott nos referimos ao próprio trabalho de Kontsevich [24], além de [19]. Outra referência que merece destaque é [11], reunindo todos esses conceitos e ferramentas e apresentando exemplos associados aos problemas da *Mirror Symmetry*.

## 1.1 Mapas estáveis

Considere $C$ um curva conexa, nodal, de gênero $g$. Tome $p_1, \cdots, p_n$ pontos não-singulares sobre essa curva. Chamamos esses pontos de **marcas** e dizemos que $C$ é uma **curva marcada**, denotada por $(C, p_1, \cdots, p_n)$. Um mapa $f : (C, p_1, \cdots, p_n) \to \mathbb{P}^r$ de grau fixo $d$ é por sua vez chamado **mapa marcado** (em $\mathbb{P}^r$) e por vezes denotado como $(C, p_1, \cdots, p_n; f)$ ou simplesmente por $f : C \to \mathbb{P}^r$ quando as marcas puderem ser omitidas. Nas definições seguintes, fixamos os inteiros $g$, $n$, $r$ e $d$.

Um **isomorfismo entre mapas marcados** $(C, p_1, \cdots, p_n; f)$ e $(C', p'_1, \cdots, p'_n, f')$ é um isomorfismo de curvas $\psi : C \to C'$ que, além de comutar com os mapas $f$ e $f'$, também preserva as marcas, ou seja, $\phi(p_i) = p'_i$.

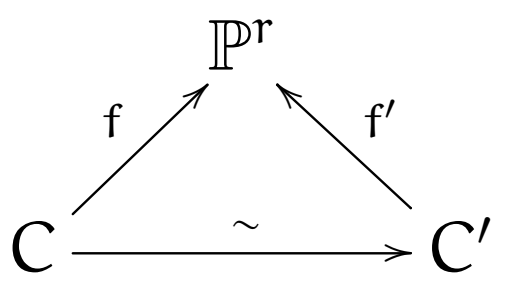

Observe que essa definição inclui automorfismos, ou seja, isomorfismos $\psi : C \to C$ tal que $f = f \circ \psi$ nas condições acima. Assim, ao tentar considerar o espaço quociente desses mapas marcados em relação a seus isomorfismos, estamos também quocientando esses automorfismos. De outro modo, estaremos interessados em estudar tais mapas como parametrizações de curvas em $\mathbb{P}^r$, a menos de suas reparametrizações.

Para criar um espaço de parâmetros compacto não podemos eliminar totalmente esses automorfismos, mas restringí-los quanto a cardinalidade. Daí surge o conceito de **estabilidade de Kontsevich**, que iremos nos referir simplesmente por estabilidade.

Diremos que $(C, p_1, \cdots, p_n; f)$ é um **mapa estável** se apresenta apenas uma quantidade finita de automorfismos. Mais precisamente, isso significa que:

(0) Componentes irredutíveis racionais contraídas por $f$ contêm ao menos 3 marcas ou nós;

(1) Componentes irredutíveis elípticas contraídas por $f$ contêm ao menos 1 marca.

Uma classe de isomorfismos que possa ser representada pelo mapa estável $(C, p_1, \cdots, p_n; f)$ será denotada por $[f]$.

Para a construção de um espaço de parâmetros para os mapas estáveis, precisamos também do conceito de famílias. Definimos uma **família de mapas estáveis** para $\mathbb{P}^r$ com $n$ marcas, gênero $g$ e grau $d$ sobre um esquema $B$ por um morfismo plano $\pi : \mathcal{C} \to B$ com seções $\sigma_1, \cdots, \sigma_n$ e um morfismo $\varphi : \mathcal{C} \to \mathbb{P}^r$. Essa estrutura é tal que, para cada $b \in B$ temos $(\mathcal{C}_b, \sigma_1(b), \cdots, \sigma_n(b); \varphi|_{\mathcal{C}_b})$ como mapa estável (para $\mathbb{P}^r$ com $n$ marcas, de gênero $g$ e grau $d$).

$$\begin{array}{ccc} \mathcal{C} & \xrightarrow{\varphi} & \mathbb{P}^r \\ {\scriptstyle \sigma_i} \uparrow \downarrow {\scriptstyle \pi} & & \\ B & & \end{array} \tag{1.1}$$

Definimos também um **isomorfismo de famílias de mapas estáveis** pela comutação dos diagramas abaixo, juntamente com o requerimento da preservação das seções (ou marcas ao longo das fibras).

$$\begin{array}{ccccc} & & \mathbb{P}^r & & \\ & {\scriptstyle \varphi} \nearrow & & \nwarrow {\scriptstyle \varphi'} & \\ \mathcal{C} & & \xrightarrow{\sim} & & \mathcal{C}' \\ {\scriptstyle \sigma_i} \uparrow \downarrow {\scriptstyle \pi} & & & & {\scriptstyle \pi'} \downarrow \uparrow {\scriptstyle \sigma'_i} \\ B & & = = = = & & B \end{array}$$

Em particular, a noção de mapa estável ou de isomorfismo de mapas estáveis pode ser recuperada a partir de uma família com base $B = \bullet$.

## 1.2 O *moduli*

Em busca de um espaço de parâmetros para esses mapas estáveis, definimos o funtor contravariante

$$\overline{\mathcal{M}}_{g,n}(\mathbb{P}^r, d) : (Sch/\mathbb{C}) \to (Set)$$

associando cada esquema $B$ ao conjunto das classes de isomorfismos de famílias de mapas estáveis para $\mathbb{P}^r$ com $n$ marcas de gênero $g$ e grau $d$ sobre $B$; morfismos entre esquemas são levados em diagramas cartesianos de famílias. Observe que $\overline{\mathcal{M}}_{g,n}(\mathbb{P}^r, d)(B)$ pode ser tomado como uma categoria na qual todo morfismo é um isomorfismo, caracterizando um grupoide. De fato, podemos associar a esse funtor (na verdade um 2-funtor em grupoides) uma subcategoria fibrada por grupoides (ou CFG) sobre $(Sch/\mathbb{C})$ como apresentado no apêndice (A).

De acordo com a construção de Kontsevich, esse funtor não é representável em geral, ou seja, não existe um espaço de *moduli* fino para mapas estáveis. Porém, temos um esquema que nos dá um espaço de *moduli* grosseiro para esse funtor, denotado por $\overline{M}_{g,n}(\mathbb{P}^r, d)$.

**Teorema 1.2.1.** *Existe um esquema projetivo $\overline{M}_{g,n}(\mathbb{P}^r, d)$ que é um espaço de* moduli *grosseiro parametrizando classes de isomorfismo de mapas estáveis para $\mathbb{P}^r$ com $n$ marcas de gênero $g$ e grau $d$.*

Em particular, dada uma família $\mathcal{C} \to B$, sempre temos um morfismo $B \to \overline{M}_{g,n}(\mathbb{P}^r, d)$, mas não necessariamente existe uma família $U \to \overline{M}_{g,n}(\mathbb{P}^r, d)$ que descreve $\mathcal{C} \to B$ como sua imagem inversa.

Os próximos parágrafos fazem uma síntese dos elementos na construção do espaço de mapas estáveis para curvas racionais sem marcas, $\overline{M}_{0,0}(\mathbb{P}^r, d)$, como forma de fixar notações e investigar melhor a estrutura desse espaço. Para construções mais gerais, citamos [16].

Consideramos inicialmente o espaço de parametrizações $W(\mathbb{P}^r, d)$ dos mapas estáveis $f : \mathbb{P}^1 \to \mathbb{P}^r$, cujas coordenadas são polinômios homogêneos de grau $d$ que não se anulam simultaneamente. Temos assim $W(\mathbb{P}^r, d) \subset \mathbb{P}(H^0(\mathbb{P}^1, \mathcal{O}_{\mathbb{P}^1}(d))^{\oplus(r+1)})$. Em particular, $\dim W(\mathbb{P}^r, d) = (r+1)(d+1) - 1$. Na intenção de parametrizar curvas racionais em $\mathbb{P}^r$ através desses mapas, devemos observar que esse espaço contém 'redundâncias': há vários mapas parametrizando uma mesma curva, além de recobrimentos, permitindo curvas com estruturas não reduzidas. Além disso, não é compacto: há presença de famílias a 1-parâmetro incompletas (vide [23], 2.2).

Consideramos um subespaço $W^*(\mathbb{P}^r, d)$ formado pelos mapas livres de automorfismos, ou seja, que não são recobrimentos. Tomando o quociente pelas reparametrizações de curvas, dada por automorfismos no domínio $\mathbb{P}^1$, podemos definir

$$M^*_{0,0}(\mathbb{P}^r, d) := W^*(\mathbb{P}^r, d)/\mathrm{Aut}(\mathbb{P}^1),$$

o qual é um espaço de *moduli* fino para classes de isomorfismos de tais mapas. Com a presença de automorfismos do mapa, a condição de estabilidade sobre esses mapas irredutíveis consegue garantir que

$$M_{0,0}(\mathbb{P}^r, d) := W(\mathbb{P}^r, d)/\mathrm{Aut}(\mathbb{P}^1),$$

é apenas um espaço de *moduli* grosseiro. A falta de uma família universal é discutida na próxima seção, (1.2.1).

A condição de estabilidade também provê uma compactificação de $M_{0,0}(\mathbb{P}^r, d)$ que adiciona os mapas com domínio redutível. Esse espaço compacto é $\overline{M}_{0,0}(\mathbb{P}^r, d)$. Para um número $n$ qualquer de marcas, temos

**Teorema 1.2.2.** $\overline{M}_{0,n}(\mathbb{P}^r, d)$ *é uma variedade normal irredutível, localmente isomorfa a um quociente de uma variedade lisa pela ação de um grupo finito. Ela contém a subvariedade lisa* $\overline{M}^*_{0,n}(\mathbb{P}^r, d)$ *como aberto denso* e $\dim \overline{M}_{0,n}(\mathbb{P}^r, d) = (r+1)(d+1) - 1 + n - 3$.

Alguns desses espaços são classicamente bem conhecidos. Por exemplo, $\overline{M}_{0,0}(\mathbb{P}^r, 1) = \overline{M}^*_{0,0}(\mathbb{P}^r, 1)$ são as Grassmannianas e $\overline{M}_{0,1}(\mathbb{P}^r, 1)$ seus fibrados universais, caso no qual temos um espaço de *moduli* fino. Já $\overline{M}_{0,0}(\mathbb{P}^2, 2)$ é isomorfo ao espaço das cônicas completas. Quando $d = 0$ e $n \geq 3$ temos a identificação $\overline{M}_{0,n}(\mathbb{P}^r, 0) = \overline{M}_{0,n}$, espaço das **curvas racionais estáveis**, uma variedade projetiva lisa e um *moduli* fino para tais curvas (vide [22]).

### 1.2.1 Morfismos Estruturais

Apesar de não termos em geral famílias universais para $\overline{M}_{0,n}(\mathbb{P}^r, d)$, temos alguns morfismos que desempenham papéis análogos:

- Avaliação: Para cada marca $p_i$, definimos o morfismo (plano) de avaliação $\nu_i$ que toma a imagem dessa marca pelo mapa.

$$\begin{array}{rrcl} \nu_i : & \overline{M}_{0,n}(\mathbb{P}^r, d) & \to & \mathbb{P}^r \\ & (C; p_1, \cdots, p_n; f) & \mapsto & f(p_i) \end{array}$$

- Esquecimento: Para cada marca $p_i$, definimos um morfismo $\pi_i$ que retira ou 'esquece' essa marca da curva domínio e a estabiliza em seguida (contraindo componentes instáveis, se necessário).

$$\begin{array}{rrcl} \pi_i : & \overline{M}_{0,n+1}(\mathbb{P}^r, d) & \to & \overline{M}_{0,n}(\mathbb{P}^r, d) \\ & (C; p_1, \cdots, p_i, \cdots, p_{n+1}; f) & \mapsto & (C; p_1, \cdots, \hat{p_i}, \cdots, p_{n+1}; f) \end{array}$$

Através desses morfismos, podemos de fato definir famílias universais para os casos nos quais $\overline{M}_{0,n}(\mathbb{P}^r, d)$ é um espaço de *moduli* fino (Grassmannianas, por exemplo). Porém, mesmo perdendo a propriedade universal no caso geral, tais mapas continuam bem definidos, constituindo famílias como abaixo:

$$\begin{array}{ccc} \overline{M}_{0,n+1}(\mathbb{P}^r, d) & \xrightarrow{\nu_{n+1}} & \mathbb{P}^r \\ \sigma_i \uparrow \downarrow \pi_{n+1} & & \\ \overline{M}_{0,n}(\mathbb{P}^r, d) & & \end{array} \quad (1.2)$$

A condição da universalidade dessa família está relacionada à presença de automorfismos não triviais. Quando livre de automorfismos, dada $[f] \in \overline{M}^*_{0,n}(\mathbb{P}^r, d)$ representado pelo mapa estável $f : C \to \mathbb{P}^r$, a fibra $F_f := \pi^{-1}_{n+1}([f])$ pode ser identificada canonicamente com a própria curva C. Por sua vez, a restrição $\nu_{n+1}|_{F_f}$ é identificada com o próprio morfismo f.

Com a presença de automorfismos, a fibra $F_f$ pode ser não reduzida, de certa forma 'absorvendo' o grau de recobrimento de f, deixando $\nu_{n+1}|_{F_f}$ como uma bijeção. Um exemplo ilustrativo: o recobrimento $(x : y) \mapsto (x^2 : y^2 : 0)$ (exemplos no cap.2 de [23]].

As seções $\sigma_i$ podem ser descritas da seguinte forma: dada classe [f], a imagem por $\sigma_i$ adiciona uma nova marca $p_{n+1}$ em C que coincide com a marca $p_i$ e em seguida a estabiliza a curva. A estabilização ocorre inserindo $p_{n+1}$ e $p_i$ em uma nova componente $\mathbb{P}^1$ de grau zero, cuja interseção com a curva C original ocorre no ponto antes ocupado pela marca $p_i$. Em particular, a curva domínio resultante é redutível e está necessariamente incluída no bordo de $\overline{M}_{0,n}(\mathbb{P}^r, d)$, descrito adiante.

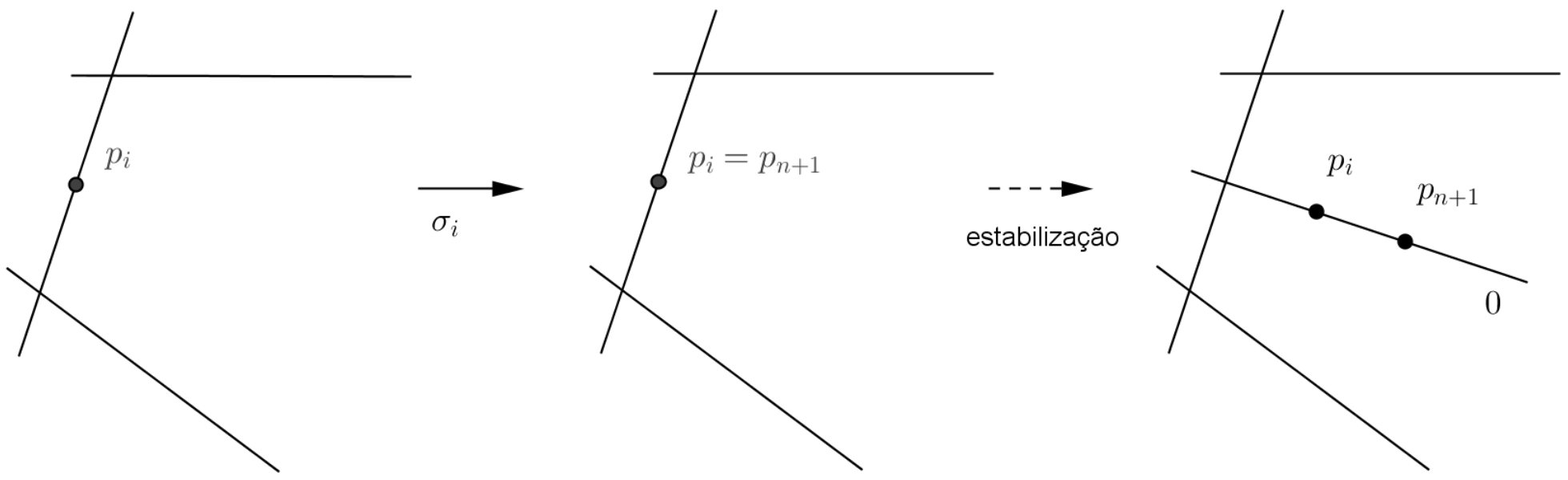

Figura 1.1: Ação da seção $\sigma_{n+1}$ sobre a curva domínio

### 1.2.2 Bordo

O espaço $\overline{M} = \overline{M}_{0,n}(\mathbb{P}^r, d)$ é construído como uma compactificação do espaço $M = M_{0,n}(\mathbb{P}^r, d)$, cujos mapas têm como domínio apenas curvas irredutíveis, ou simplesmente $\mathbb{P}^1$. O conjunto $\Delta := \overline{M} \setminus M$, é portanto o bordo de $\overline{M}$, formado por mapas cujo domínio são curvas racionais redutíveis, também chamadas de **árvore de retas projetivas**. De modo geral, os mapas estáveis do bordo são aqueles obtidos por limites de famílias de mapas em $M_{0,n}(\mathbb{P}^r, d)$.

A decomposição das curvas domínio dos mapas estáveis no bordo gera uma decomposição e uma estrutura recursiva em $\Delta$. De fato, o bordo se divide em várias componentes divisoriais $\Delta(d_A, d_B; A, B)$, cada uma delas associada a uma partição do grau $d_A + d_B = d$ e do conjunto $A \sqcup B = \{p_1, \cdots, p_n\}$ das $n$ marcas. Explicitamente: um nó $p$ determina, como ponto de interseção, duas componentes $C_A$ e $C_B$, com conjunto de marcas $A \cup \{p\}$ e $B \cup \{p\}$ respectivamente e as restrições $f|_{C_A}$ e $f|_{C_B}$ com graus $d_A$ e $d_B$. Assim, com exceção dos casos simétricos, temos um isomorfismo

$$\Delta(d_A, d_B; A, B) \simeq \overline{M}_{0,|A|+1}(\mathbb{P}^r, d_A) \times_{\mathbb{P}^r} \overline{M}_{0,|B|+1}(\mathbb{P}^r, d_B)$$

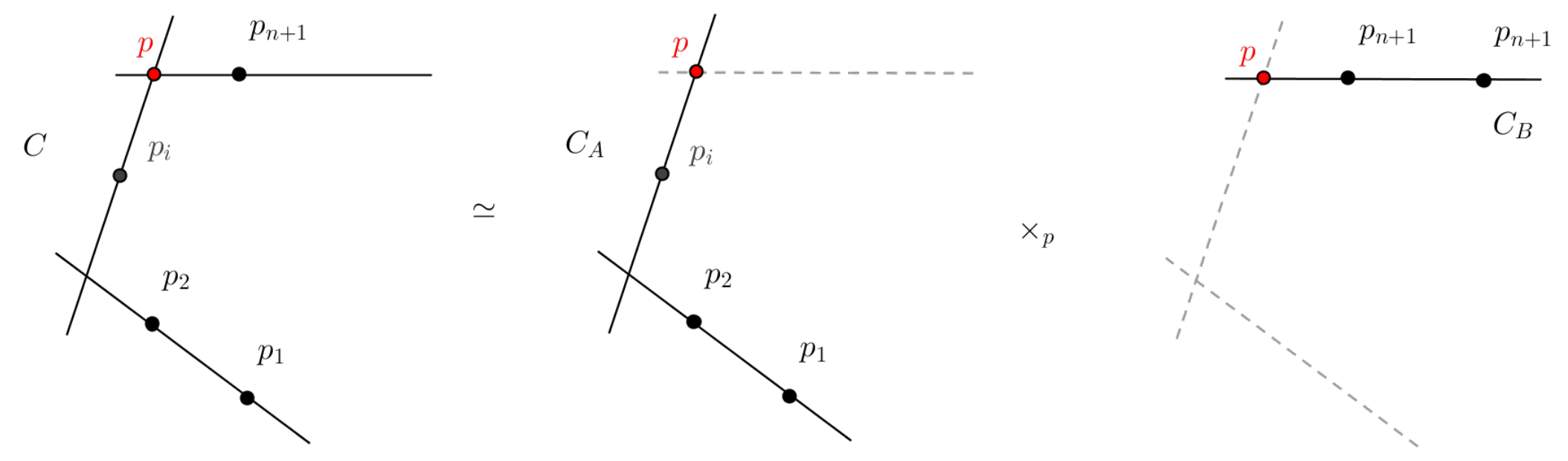

Figura 1.2: Decomposição no bordo

Ou seja, o bordo de $\overline{M}_{0,n}(\mathbb{P}^r, d)$ pode ser descrito por mapas estáveis com mais marcas e graus menores. Dessa estrutura recursiva é que surgem, por exemplo, as equações recursivas para os invariantes de Gromov-Witten, de acordo com o trabalho de Kontsevich e Manin [25].

## 1.3 O *stack*

Como apresentado pelo teorema (1.2.1), o esquema $\overline{M}_{0,n}(\mathbb{P}^r, d)$ apresenta em geral singularidades quocientes. Além disso, como espaço de *moduli* é apenas grosseiro. Isso se deve ao fato de que as famílias de mapas estáveis guardam 'multiplicidades' dadas por seus automorfismos, que não são explicitadas no esquema $\overline{M}_{0,n}(\mathbb{P}^r, d)$. Para estudarmos os mapas estáveis com todas suas informações, precisamos de uma extensão do conceito de esquema. Para isso - uma explicação pouco rigorosa - trocamos nosso conjunto de pontos por uma categoria. Essa categoria deve ser capaz de codificar não só a informação das famílias de mapas estáveis (objetos), como também seus automorfismos (morfismos). Além disso, deve possui propriedades 'algébricas', que a torne próxima da estrutura de um esquema. Tal extensão é conhecida como *stack algébrico*. Mais precisamente, o conceito de *stack* parte do estudo do funtor do problema de *moduli*, como descrito mais formalmente no apêndice A.

**Teorema 1.3.1.** *O funtor* $\overline{\mathcal{M}}_{0,n}(\mathbb{P}^r, d)$ *é um* stack *algébrico próprio e liso sobre* $\mathbb{C}$.

De acordo com o apresentado no apêndice, a descrição de um *stack* como categoria é equivalente e pode parecer mais natural:

- objetos: classes de isomorfismos de famílias planas de mapas estáveis $\pi : C \to B$, como definido por (1.1);
- morfismos: diagramas cartesianos entre famílias de mapas estáveis.

Essa categoria é equipada com um funtor contravariante de projeção

$$p : \overline{\mathcal{M}}_{0,n}(\mathbb{P}^r, d) \to (Sch/\mathbb{C}),$$

que leva uma família de mapas estáveis $\pi : C \to B$ em sua base B e um diagrama cartesiano de mapas estáveis no morfismo de base que o constitui. Além do mais, é uma CFG: a fibra de p sobre um esquema B (e a identidade $id_B$) é uma subcategoria $\overline{\mathcal{M}}_B$ na qual todos morfismos são de fato isomorfismos. Isso nos permite construir uma ideia intuitiva sobre a estrutura de *stacks* para mapas estáveis: sobre um ponto $\bullet$, a fibra $\overline{\mathcal{M}}_\bullet$ é o mapa $C \to \mathbb{P}^r$ com seus automorfismos.

O que o torna $\overline{\mathcal{M}}_{0,n}(\mathbb{P}^r, d)$ um *stack* algébrico (no sentido de Deligne-Mumford) é o fato dessa categoria apresentar características análogas às de um esquema, graças a propriedades de sua diagonal e a existência de um atlas a partir de um esquema. Na prática, podemos realizar as principais propriedades de $\overline{\mathcal{M}}_{0,n}(\mathbb{P}^r, d)$, desde as topológicas, produtos fibrados e construções de feixes, como no caso de esquemas lisos. Ainda mais, podemos construir uma teoria de interseção para $\overline{\mathcal{M}}_{0,n}(\mathbb{P}^r, d)$, como melhor descrito na próxima seção.

Antes, devemos entender a relação precisa entre $\overline{M}_{0,n}(\mathbb{P}^r, d)$ e $\overline{\mathcal{M}}_{0,n}(\mathbb{P}^r, d)$. Assim como todo esquema, podemos realizar $\overline{M}_{0,n}(\mathbb{P}^r, d)$ como um *stack* através de seu funtor de pontos, que continuaremos denotando da mesma forma - um abuso bem justificável. Por definição de um espaço de *moduli* para um *stack*, temos um morfismos de *stacks* (uma transformação natural)

$$m : \overline{\mathcal{M}}_{0,n}(\mathbb{P}^r, d) \to \overline{M}_{0,n}(\mathbb{P}^r, d). \tag{1.3}$$

Essa transformação pode ser descrita do seguinte modo: uma classe de famílias de mapas estáveis $\pi : \mathcal{C} \to B$ tem como imagem seu morfismo classificante $B \to \overline{M}_{0,n}(\mathbb{P}^r, d)$.

## 1.4 Teoria de interseção

O grupo de Chow $A_*(\overline{\mathcal{M}}_{0,n}(\mathbb{P}^r, d))_{\mathbb{Q}}$ é construído de forma análoga ao caso de esquemas lisos, porém definido sobre os racionais. A construção é suas propriedades são descritas brevemente no apêndice A ou com mais detalhes em [29]. Em particular, temos a classe fundamental $[\overline{\mathcal{M}}_{0,n}(\mathbb{P}^r, d)]$.

O grupo de Chow acima é um módulo sobre o anel de Chow bivariante $A^*(\overline{\mathcal{M}}_{0,n}(\mathbb{P}^r, d))_{\mathbb{Q}}$ e ainda vale a dualidade $A^*(\overline{\mathcal{M}}_{0,n}(\mathbb{P}^r, d))_{\mathbb{Q}} \to A_*(\overline{\mathcal{M}}_{0,n}(\mathbb{P}^r, d))_{\mathbb{Q}}$ dada pela interseção com a classe fundamental.

Como consequência das proposições A.6.1 e A.6.2, se consideramos o morfismo do *moduli* $m : \overline{\mathcal{M}}_{0,n}(\mathbb{P}^r, d) \to \overline{M}_{0,n}(\mathbb{P}^r, d)$, temos isomorfismos $m_* : A_*(\overline{\mathcal{M}}_{0,n}(\mathbb{P}^r, d))_{\mathbb{Q}} \to A_*(\overline{M}_{0,n}(\mathbb{P}^r, d))_{\mathbb{Q}}$ e $m^* : A^*(\overline{M}_{0,n}(\mathbb{P}^r, d))_{\mathbb{Q}} \to A^*(\overline{\mathcal{M}}_{0,n}(\mathbb{P}^r, d))_{\mathbb{Q}}$. Em particular, segue a dualidade $A^*(\overline{M}_{0,n}(\mathbb{P}^r, d))_{\mathbb{Q}} \to A_*(\overline{M}_{0,n}(\mathbb{P}^r, d))_{\mathbb{Q}}$ também sobre o espaço de moduli dos mapas estáveis. Para o nosso caso, com $g = 0$ e sendo $\mathbb{P}^r$ espaço convexo, vale que $[\overline{M}_{0,n}(\mathbb{P}^r, d)]^v = [\overline{M}_{0,n}(\mathbb{P}^r, d)]$ ([11], 7.1.5).

Com isso, $A_*(\overline{M}_{0,n}(\mathbb{P}^r, d))_{\mathbb{Q}}$ apresenta a maioria das propriedades formais que conhecemos para a teoria de interseção para esquemas não singulares, uma vez que herda a estrutura do grupo de Chow de um *stack* liso. O fato do isomorfismo ser sobre o corpo dos racionais $\mathbb{Q}$, pode ser pensando, grosso modo, para descontar as repetições geradas pelos automorfismos dos objetos. Uma manifestação disso fica explícita na fórmula de localização para *stacks*, apresentada a seguir.

### 1.4.1 Fórmula de localização

Baseado na construção do apêndice A, considere o toro $T = (\mathbb{C}^*)^{r+1}$ agindo sobre $\mathbb{P}^r$. Isso induz uma ação em $\overline{\mathcal{M}}_{0,n}(\mathbb{P}^r, d)$ e podemos definir o anel de Chow T-equivariante, a coeficientes racionais, $A_*^T(\overline{\mathcal{M}}_{0,n}(\mathbb{P}^r, d))_{\mathbb{Q}}$. Além disso, se $\lambda_1, \cdots, \lambda_{r+1}$ denotam os pesos da ação de T sobre $\mathbb{P}^r$, o anel de Chow T-equivariante de um ponto será $\mathbb{C}[\lambda_1, \cdots, \lambda_{r+1}]$. Se $\mathcal{R}_T = \mathbb{C}(\lambda_1, \cdots, \lambda_{r+1})$ é seu corpo de frações, definimos a localização do anel de Chow T-equivariante $A_*^T(\overline{\mathcal{M}}_{0,n}(\mathbb{P}^r, d))_{\mathbb{Q}} \otimes \mathcal{R}_T$. A fórmula de localização para *stacks*, também mostrada no apêndice A, permite calcular o grau de uma classe $\alpha$ em $A_*^T(\overline{\mathcal{M}}_{0,n}(\mathbb{P}^r, d))_{\mathbb{Q}} \otimes \mathcal{R}_T$ a partir das componentes do lugar dos pontos fixos pela ação induzida de T sobre $\overline{\mathcal{M}}_{0,n}(\mathbb{P}^r, d)$. De acordo com [24] e [19], as componentes podem ser parametrizadas por grafos $\Gamma$, que serão descritos melhor adiante. Denotamos por $\iota_\Gamma : \overline{\mathcal{M}}_\Gamma \hookrightarrow \overline{\mathcal{M}}_{0,n}(\mathbb{P}^r, d)$ a inclusão da componente associada ao grafo $\Gamma$. Com as notações e definições acima explicitadas, temos

**Teorema 1.4.1** (Fórmula de localização para mapas estáveis). *Se* $\alpha \in A_*^T(\overline{\mathcal{M}}_{0,n}(\mathbb{P}^r, d))_{\mathbb{Q}} \otimes \mathcal{R}_T$, *então*

$$\int_{\overline{\mathcal{M}}_{0,n}(\mathbb{P}^r,d)} \alpha = \sum_{\Gamma} \frac{1}{a(\Gamma)} \int_{\overline{\mathcal{M}}_\Gamma^T} \frac{\iota_\Gamma^* \alpha}{\mathrm{Euler}^T(\mathcal{N}_\Gamma)},$$

*onde* $a(\Gamma)$ *é a ordem do grupo de automorfismos agindo sobre a componente* $\overline{\mathcal{M}}_\Gamma$ *e* $\mathcal{N}_\Gamma$ *seu fibrado normal.*

A partir de agora, vamos descrever melhor a estrutura combinatória que rege o lugar dos mapas estáveis fixos em $\overline{\mathcal{M}}_{0,n}(\mathbb{P}^r, d)$ pela ação de T.

Começamos com uma ação sobre $\mathbb{P}^r$ que fixa seus pontos base,

$$q_0 = (1 : 0 : \cdots : 0), \cdots, q_r = (0 : 0 : \cdots : 1).$$

A partir dos mapas de avaliação, temos uma ação induzida sobre $\overline{M}_{0,n}(\mathbb{P}^r, d)$ com a seguinte caracterização: uma classe representada pelo mapa estável $(C, p_1, \cdots, p_n; f)$ é fixada pela ação se, e somente se,

- Uma componente de C na qual $f|_C$ tem grau positivo é mapeada em uma das retas coordenadas entre os pontos $q_i$;
- Uma componente contraída em um ponto deve ter como imagem algum $q_i$;
- As marcas, nós e pontos de ramificação também têm como imagem algum dos pontos $q_i$.

Cada componente de pontos fixos pode ser identificada por dados combinatórios condensados em uma estrutura de grafo ponderado colorido (para definições, ver apêndice C) $\Gamma = (V, i, S; E, \delta)$:

(V) Um vértice $v$ é uma componente conexa $C_v$ do conjunto $f^{-1}(q_0, \cdots, q_r)$, ou seja, corresponde a um nó em C ou uma componente de grau 0;

  (i) A cada vértice $v$ associamos um rótulo $i(v) = i_v \in \{0, \cdots, r\}$ de modo que $f(C_v) = q_{i_v}$;

  (S) A cada vértice $v$ associamos um subconjunto $S(v) = S_v$ de $\{p_1, \cdots, p_n\}$ formado pelos pontos $p_i$ que estão em $C_v$;

(E) Uma aresta $e$ corresponde a uma componente $C_e$ da curva que seja mapeada em alguma reta coordenada;

  (δ) A cada aresta $e$ associamos um peso definido pelo pelo grau: $\delta(e) := \deg f|_{C_e} = d_e$

Por estarmos restritos a curvas racionais, o grafo $\Gamma$ é necessariamente uma árvore. Independentemente disso, por construção não apresenta laços: dada uma aresta $e = (v, w)$, $C_e$ tem como imagem a reta coordenada definida pelos pontos $q_{i_v}$ e $q_{i_w}$, com $i_v \neq i_w$. Isso também justifica que $i$ pode ser tomada como uma coloração de $\Gamma$. Também, os pesos formam uma partição do grau do mapa estável, $d = \sum_e d_e$, assim como os subconjuntos de marcas são uma partição do conjunto das marcas, $\{p_1, \cdots, p_n\} = \bigsqcup_v S_v$. Essas condições definem unicamente a estrutura do grafo $\Gamma$.

Por outro lado, dada uma árvore $\Gamma$ como acima, podemos construir não apenas um mapa estável, mas um *substack* $\overline{\mathcal{M}}_\Gamma \subset \overline{\mathcal{M}}_{0,n}(\mathbb{P}^r, d)$ formado pelos mapas estáveis com árvore $\Gamma$. Para entender a estrutura desse *substack*, tome um elemento $(C, p_1, \cdots, p_n; f) \in \overline{\mathcal{M}}_\Gamma$. A cada vértice $v$ com $C_v$ uma curva, associamos o número $n(v) = |S_v| + val(v)$, que conta os pontos especiais sobre a componente: as $|S_v|$ marcas e os $val(v)$ nós. Em particular, $C_v$ é uma curva se, e somente se, $n(v) \geq 3$. Portanto $(C_v, S_v)$ é uma curva estável e $(C_v, S_v; f|_{C_v})$ um mapa estável de grau zero. Podemos então construir um morfismo

$$\prod_{v \mid \dim C_v = 1} \overline{M}_{0,n(v)} \to \overline{\mathcal{M}}_\Gamma \tag{1.4}$$

tal que, dada a aresta $e = (v, w)$, tomamos $C_e \cong \mathbb{P}^1$ com nós em $0$ e $\infty$ para $C_v$ e $C_w$, $f(C_e)$ a reta coordenada entre $q_{i_v}$ e $q_{i_w}$ e $f|_{C_e}$ de grau $d_e$. Com isso produzimos um mapa estável em $\overline{\mathcal{M}}_\Gamma$. Se não há componentes com dimensão positiva, tomamos o domínio do mapa como um ponto.

Os mapas construídos dessa forma trazem naturalmente automorfismos: aqueles provenientes da própria estrutura do grafo e aqueles provenientes do recobrimento das retas coordenadas, no caso das componentes de grau positivo. O grupo de automorfismos do mapa estável fixo, $A(\Gamma)$, pode ser descrito pela sequência exata

$$0 \longrightarrow \prod_{e \in E} \mathbb{Z}_{d_e} \longrightarrow A(\Gamma) \longrightarrow Aut(\Gamma) \longrightarrow 0$$

onde $Aut(\Gamma)$ denota os automorfismos de $\Gamma$ como grafo ponderado colorido (preservando não só incidência, mas também pesos e cores). Com essas informações, temos o primeiro ingrediente para calcular explicitamente o grau de uma classe $\alpha$ através da fórmula de localização (1.4.1).

Explicitando os elementos: a enumeração dos grafos ponderados rotulados coloridos $\Gamma$ a partir de um inteiro $d$, conjunto $S$ e $r+1$ cores; a ordem do grupo de automorfismos $A(\Gamma)$; a expressão em pesos da ação para a classe equivariante $\iota_\Gamma^* \alpha$; a expressão em pesos da ação para a classe de Euler do fibrado normal $N_\Gamma$.

A classe de Euler do fibrado normal é descrita explicitamente em [24] através da contribuição de dois fatores, associados aos vértices e arestas respectivamente. Temos

$$\frac{1}{\mathrm{Euler}^T(\mathcal{N}_\Gamma)} = V(\Gamma)E(\Gamma)$$

onde

$$V(\Gamma)=\prod_{v\in V}\left(\left(\prod_{j\neq i_v}(\lambda_{i_v}-\lambda_j)\right)^{\mathrm{val}(v)-1}\left(\sum_{e=(v,w)}\frac{d_e}{\lambda_{i_v}-\lambda_{i_w}}\right)^{\mathrm{val}(v)-3}\prod_{e=(v,w)}\frac{d_e}{\lambda_{i_v}-\lambda_{i_w}}\right) \tag{1.5}$$

e

$$E(\Gamma)=\prod_{e=(v,w)\in E}\left(\frac{(-1)^{d_e}(\frac{d_e}{\lambda_{i_v}-\lambda_{i_w}})^{2d_e}}{d_e!^2}\prod_{k\neq i_v,i_w}\prod_{\alpha=0}^{d_e}\frac{1}{\frac{\alpha\lambda_{i_v}+(d_e-\alpha)\lambda_{i_w}}{d_e}-\lambda_k}\right) \tag{1.6}$$

Vamos descrever a enumeração das árvores $\Gamma$ através de exemplos no caso sem marcas e em $\mathbb{P}^3$ - nossa real situação de interesse. Consideremos uma ação em $\mathbb{P}^3$ induzindo uma ação em $\overline{M}_{0,0}(\mathbb{P}^3,3)$. Nesse caso, os pontos fixos são isolados, uma vez que não temos componentes $C_v$ de dimensão positiva (de acordo com a descrição do morfismo (1.4). Um mapa estável fixo típico pode ser representado pelo diagrama da figura (1.3), onde identificamos cada um dos 4 pontos base de $\mathbb{P}^3$ com uma cor. Decodificamos todas essas informações geométricas no grafo também apresentado na figura.

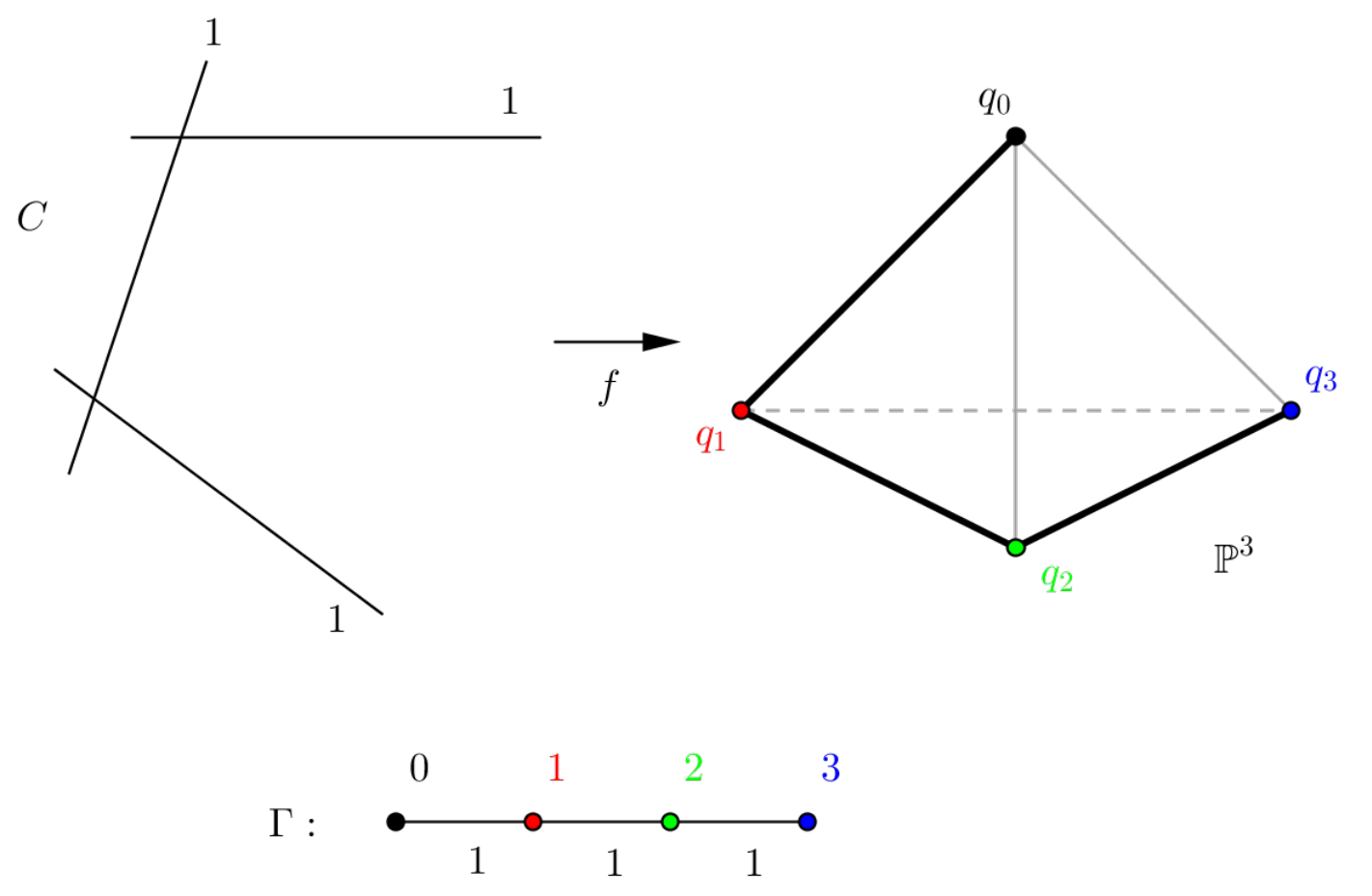

Figura 1.3: Mapa estável fixo de grau 3

Para o caso de um mapa estável fixo com componente de grau zero, temos a estrutura como mostrada no exemplo da figura (1.4).

O primeiro caso em que não teremos todos pontos fixos isolados é para $d = 4$. Nesse caso, há um tipo combinatório de grafo (desconsiderando coloração) que apresenta uma componente de grau zero, como mostrado na figura (1.5). Nesse caso, a componente de pontos fixos é isomorfa a $\overline{\mathrm{M}}_{0,4} \simeq \mathbb{P}^1$, que corresponde à razão cruzada entre os 4 pontos especiais (nós) na componente de grau zero. Na figura, a reta $\mathrm{q}_0\mathrm{q}_1$ é recoberta por duas componentes, indicado pelo $2\times$ na imagem e representado no grafo pelas duas arestas cujos vértices têm cores iguais.

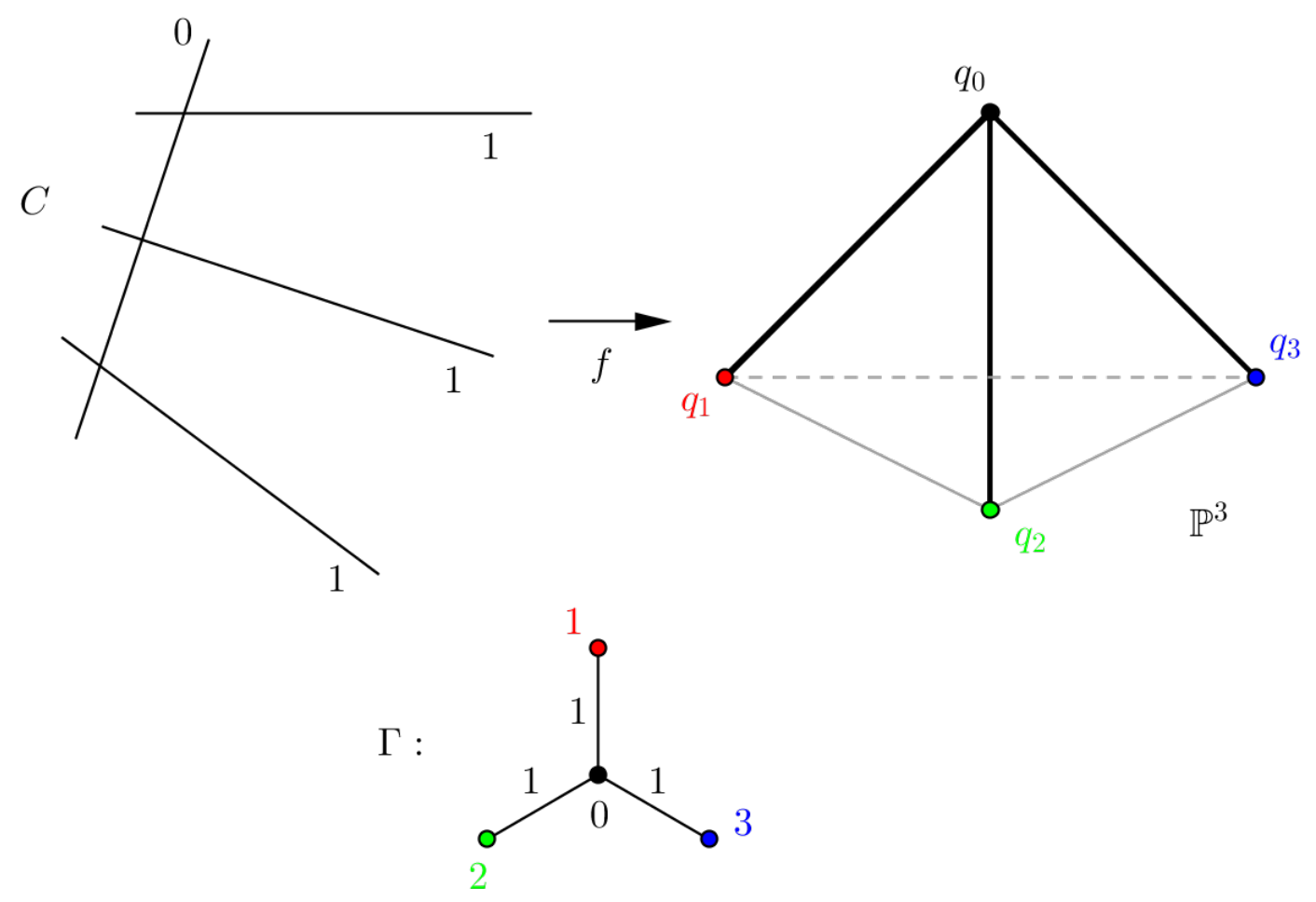

Figura 1.4: Mapa estável fixo de grau 3 com componente contraída

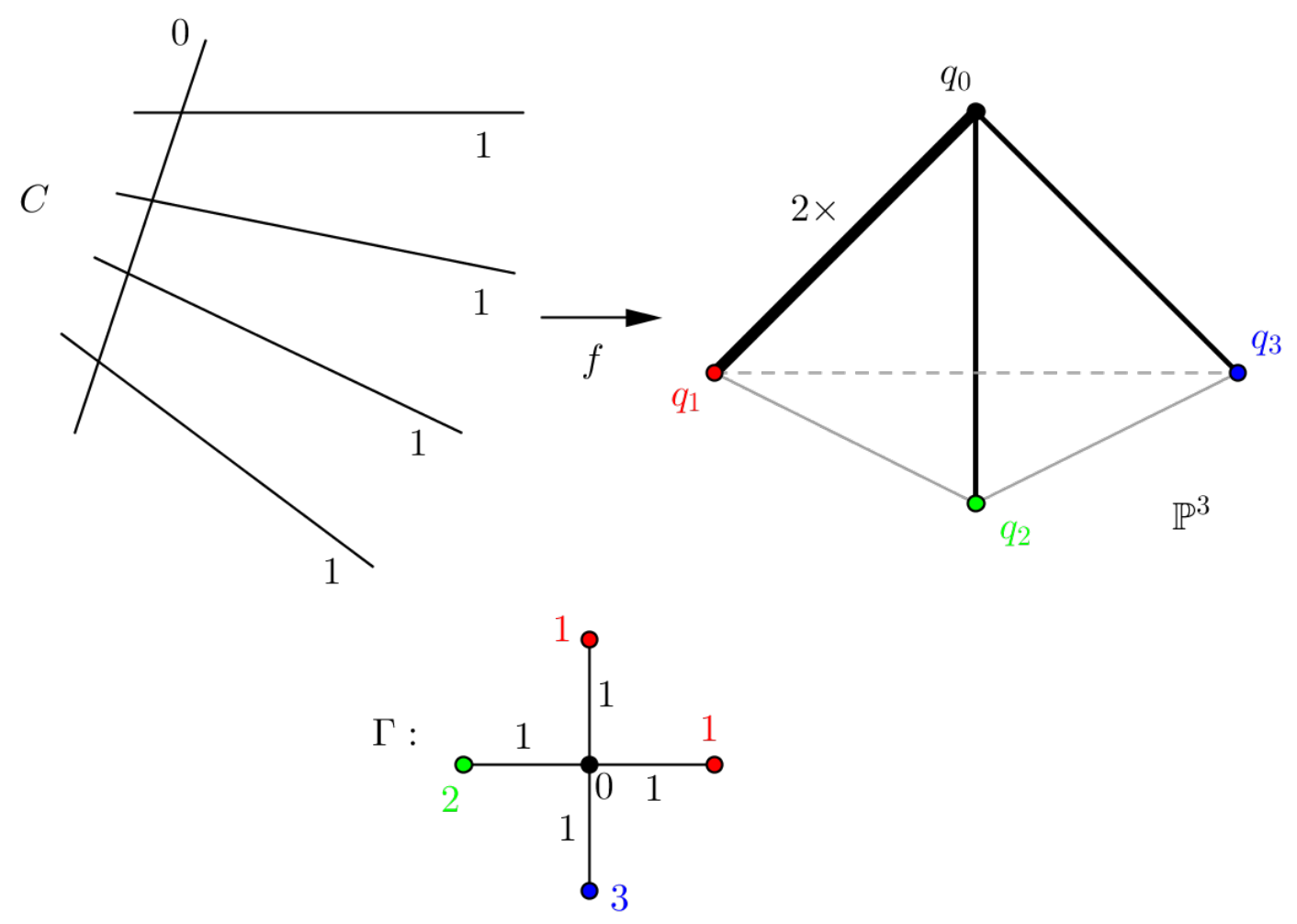

Figura 1.5: Mapa estável fixo de grau 4 de uma componente de dimensão positiva

# Capítulo 2

# Curvas de contato no espaço projetivo

Seguindo a proposta apresentada em [26], iremos determinar um espaço de parâmetros para a classe das curvas de contato racionais no espaço projetivo $\mathbb{P}^3$, assim como definir alguns de seus invariantes enumerativos. Também faz parte dessa proposta obter tal espaço de parâmetros como subespaço do *stack* $\overline{\mathcal{M}}_{0,0}(\mathbb{P}^3, d)$, para cada grau $d$.

Primeiramente, definimos nosso principal objeto de estudo: curvas de contato. Essas curvas aparecem em espaços de dimensão ímpar com a chamada estrutura de contato, proveniente de um espaço simplético. Os conceitos e nomenclaturas para as estruturas simplética e de contato são apresentados em (2.1), tendo como referência básica o texto [4].

O primeiro caso tratado é o de curvas de contato planas, descrito em [26]. Em seguida apresentamos uma descrição mais geral através da caracterização legendriana dessas curvas, como em [10]. Apresentamos ainda a caracterização da estrutura de contato e das curvas de contato via fibrações, ainda seguindo o projeto em [26].

Finalmente, em (2.4) iniciamos nosso estudo principal, definindo o espaço de parâmetros $\mathcal{L}_d$ e os invariantes $N_d$ associados às curvas de contato como mapas estáveis de grau $d$. Também, fazemos uma releitura dos exemplos da subseção anterior frente essa nova ótica, discutindo a relação entre a dimensão esperada e a dimensão efetiva de $\mathcal{L}_d$, assim como o significado enumerativo dos invariantes $N_d$.

## 2.1 Estrutura de contato

### 2.1.1 A estrutura simplética

Em um espaço vetorial complexo $V$, dizemos que uma forma bilinear $\omega : V \times V \to \mathbb{C}$ é **simplética** se

- é anti-simétrica: $\omega(v, w) = -\omega(w, v)$
- é não degenerada: para todo $v$, existe $w$ tal que $\omega(v, w) \neq 0$

Observe que se um espaço vetorial possui forma simplética, ele deve ter dimensão par. Portanto, consideremos daqui em diante o espaço $V = \mathbb{C}^{2r}$, $r$ inteiro. Se tomamos em $\mathbb{C}^{2r}$ coordenadas $x = (x_1, \cdots, x_r, y_1, \cdots, y_r)$, uma 2-forma simplética, única a menos de uma mudança de coordenadas simpléticas, é

$$\omega = \sum_{i=1}^{r} dx_i \wedge dy_i.$$

Um espaço vetorial munido de uma forma simplética é chamado **espaço simplético**. Dado um subespaço linear $W$, definimos o seu complementar relativo à forma $\omega$ por

$$W^\perp = \{v \in \mathbb{C}^{2r} \mid \omega(w, v) = 0 \ \forall w \in W\}.$$

Claramente, temos as relações $(W^\perp)^\perp = W$ e $\dim W + \dim W^\perp = 2r$. Do fato dessa forma ser simplética, podemos ainda definir alguns subespaços especiais:

- $W$ é simplético se $W \cap W^\perp = \{0\}$
- $W$ é isotrópico se $W \subseteq W^\perp$
- $W$ é coisotrópico se $W^\perp \subseteq W$, ou se $W^\perp$ é isotrópico.
- $W$ é lagrangeano se $W^\perp = W$, ou $W$ é isotrópico e coisotrópico.

Observe que se $W$ é isotrópico, temos $\dim W \leq r$; se é coisotrópico, $\dim W \geq r$. Em particular, se $W$ é lagrangeano, $\dim W = r$. Outro modo de dizer: $W$ é lagrangeano se é um subespaço maximal tal que $\omega|_W \equiv 0$.

Também temos uma classificação das variedades em um espaço simplético: uma subvariedade é chamada isotrópica, coisotrópica ou lagrangeana se o são seus espaços tangentes em cada ponto não singular.

### 2.1.2 A distribuição de contato

Proveniente da estrutura simplética em um espaço vetorial de dimensão par $(\mathbb{C}^{2r}, \omega)$, temos uma estrutura análoga para o espaço projetivo de dimensão ímpar $\mathbb{P}^{2r-1} = \mathbb{P}(\mathbb{C}^{2r})$. Para construir essa estrutura, considere o isomorfismo

$$\begin{array}{rcl} \mathbb{C}^{2r} & \to & (\mathbb{C}^{2r})^\vee \\ v & \mapsto & \omega(v, \cdot) \end{array}$$

Definimos então uma 1-forma diferencial, também denotada por $\omega$, tal que em cada ponto $p \in \mathbb{P}^{2r-1}$ seja $\omega_p := \omega(p, \cdot)$, definida a menos de múltiplo escalar. Essa forma diferencial define uma distribuição de hiperplanos em $T\mathbb{P}^{2r-1}$ determinados por $p \mapsto \operatorname{Nuc} \omega_p$. Observe que $\omega_p(p) = 0$, o que nos permite realizar $\omega$ como uma distribuição de hiperplanos em $\mathbb{P}^{2r-1}$. Vale também que $\omega \wedge (d\omega)^{\wedge r-1}$ é não degenerada. A uma distribuição com essa última propriedade damos o nome de **distribuição de contato** e dizemos que ela define uma **estrutura de contato** em $\mathbb{P}^{2r-1}$. Observe que a condição de integrabilidade para uma distribuição de codimensão 1 dada por uma 1-forma $\omega$ em geral é $\omega \wedge d\omega = 0$. Por isso, uma distribuição de contato é tida como de 'não integrabilidade máxima'.

Estamos interessados em analisar a condição de tangência a essa distribuição.

**Definição 2.1.1.** *Seja $X$ uma subvariedade de $\mathbb{P}^{2r-1}$ com sua estrutura de contato $\omega$. Dizemos que $X$ é uma **subvariedade de contato** se é tangente à distruibuição de contato em cada um de seus pontos não singulares. Isto é, se*

$$T_pX \subset \operatorname{Nuc} \omega_p, \ \forall p \in X \setminus \operatorname{Sing} X$$

Uma subvariedade de contato $X$ tem no máximo dimensão $r-1$. De fato, a 2-forma $d\omega_p$ dá uma estrutura simplética ao hiperplano $\operatorname{Nuc} \omega_p$, no qual $T_pX$ é então um subespaço isotrópico. Se $\dim X = r-1$, dizemos que $X$ é **legendriana**. Outra definição: dizemos que uma variedade $X \in \mathbb{P}^{2r-1}$ é **involutiva** se seu tangente é coisotrópico, tendo no mínimo dimensão $r-1$.

Em particular - e na verdade nosso caso de interesse - uma curva de contato em $\mathbb{P}^3$ é legendriana (e portanto, também involutiva).

## 2.2 Curvas de contato planas

O resultado a seguir ([26], proposição 2.2), caracteriza as curvas de contato planas em $\mathbb{P}^3$ através da involutividade.

**Proposição 2.2.1.** *Seja* $X$ *uma curva involutiva plana de grau* $d \geq 2$ *em* $\mathbb{P}^3$*. Então* $X$ *é a união de* $d$ *retas passando pelo ponto simplético do plano definido por* $X$*.*

Considere o espaço que parametriza as curvas de contato planas de grau d. Segue da proposição acima que temos $3 + d$ parâmetros para descrever tal espaço: 3 para a escolha do plano de contato e d escolhas de retas passando pelo ponto simplético e contidas no plano de contato. Dado que a condição de incidência a uma reta em posição geral no espaço é divisorial, um invariante de interesse seria o que conta o número de curvas de contato planas incidentes a $3 + d$ retas genéricas. Esse invariante é calculado em [26]:

**Proposição 2.2.2** (Levcovitz-Vainsencher). *O número de curvas planas involutivas de grau* $d \geq 2$ *encontrando* $3 + d$ *retas em posição geral de* $\mathbb{P}^3$ *é* $20d\binom{d+3}{5}$*.*

Dessas proposições já somos capazes de obter um invariante enumerativo para retas e cônicas de contato:

1. Retas de contato: Podemos parametrizar o conjunto das retas de contato em $\mathbb{P}^3$ por uma variedade de dimensão 3 e grau 2, como uma seção hiperplana na quádrica de Plücker. Assim, o número de retas de contato incidente a 3 retas em posição geral de $\mathbb{P}^3$ é 2. O espaço dessas retas é a Grassmanniana Lagrangeana de retas em $\mathbb{P}^3$ [3].

2. Cônicas de contato: Pela proposição (2.2.1), toda cônica de contato é um par de retas passando pelo ponto simplético do plano da distribuição de contato que a contém. Assim, a dimensão da variedade de cônicas de contato é $3 + 2 = 5$. O número de tais cônicas incidentes a 5 retas em posição geral então é 40.

## 2.3 Caracterização legendriana

Para estudar o caso não plano, $d \geq 3$, é necessário extrair mais informações da estrutura de contato em $\mathbb{P}^3$. O caso de cúbicas já nos indica isso. Pelas proposições da seção anterior, temos dimensão 6 para conjunto das cúbicas de contato planas e 360 delas incidentes a 6 retas em posição geral. Mas dentre as cúbicas já temos o primeiro caso em que as proposições anteriores não cobrem, como a cúbica reversa. De fato, a cúbica reversa de parametrização $f(s:t) = (s^3 : \frac{1}{3}t^3 : st^2 : s^2t)$ é uma curva legendriana - e portanto de contato. Para verificar isso, de acordo com o exemplo em [10], escolhemos $\omega$ apropriadamente como definida pela matriz

$$\begin{bmatrix} 0 & 1 & 0 & 0 \\ -1 & 0 & 0 & 0 \\ 0 & 0 & 0 & 1 \\ 0 & 0 & -1 & 0 \end{bmatrix}. \tag{2.1}$$

Tomamos os vetores geradores do espaço tangente ao cone em $f(s:t)$

$$f_s = (3s^2 : 0 : t^2 : 2st) \quad \text{e} \quad f_t = (0 : t^2 : 2st : s^2),$$

e calculamos

$$\omega(f_s, f_t) = 3s^2 \cdot t^2 - 2st \cdot 2st + t^2 \cdot s^2 = 0,$$

para todo $(s:t)$, mostrando que a cúbica é de contato.

De modo mais geral, em [10], seção 3.1, é verificado que, em relação à forma de contato dada pela matriz (2.1), as curvas de parametrização

$$(s:t) \mapsto \left( s^{k+l} : \frac{k-l}{k+l} t^{k+l} : s^l t^k : s^k t^l \right), \tag{2.2}$$

com $k, l$ coprimos, são curvas de contato. Temos então um gerador de exemplos para curvas de contato não planas de grau arbitrário. Tais exemplos seguem do seguinte resultado:

**Teorema 2.3.1** (Buczynski). *Uma curva racional* $X \subset \mathbb{P}^3$ *com parametrização*

$$f(s:t) = \big(f_0(s,t) : f_1(s,t) : f_2(s,t) : f_3(s,t)\big)$$

*é legendriana (em relação a forma (2.1)) se e somente se* $f_0 \equiv 0$ *e* $X$ *é uma reta legendriana, ou as funções* $f_i(1,t)$ *satisfazem a equação diferencial*

$$f_1' \cdot f_0 - f_1 \cdot f_0' = f_2' \cdot f_3 - f_2 \cdot f_3' \tag{2.3}$$

*Prova:*

Seja $W_d = W(\mathbb{P}^3, d) \subset \mathbb{P}\big(H^0(\mathbb{P}^1, \mathcal{O}_{\mathbb{P}^1}(d))^{\oplus 4}\big)$ o espaço de parametrizações para curvas de grau d em $\mathbb{P}^3$ (1.2), ou seja, mapas $f(s:t) = \big(f_0(s,t) : f_1(s,t) : f_2(s,t) : f_3(s,t)\big)$ definidos a menos de um múltiplo com coordenadas $f_i$ polinômios homogêneos em $s, t$ que não se anulam simultaneamente. No aberto onde $s \neq 0$, tomamos $f_i(1,t)$ polinômios de grau menor ou igual a d. O espaço tangente sobre pontos não singulares do cone afim $\big(\alpha f_0(1,t), \alpha f_1(1,t), \alpha f_2(1,t), \alpha f_3(1,t)\big)$ é gerado pelos vetores $(f_0, f_1, f_2, f_3)$ e $(f_0', f_1', f_2', f_3')$. Assim,

$$\begin{bmatrix} f_0 & f_1 & f_2 & f_3 \end{bmatrix} \begin{bmatrix} 0 & 1 & 0 & 0 \\ -1 & 0 & 0 & 0 \\ 0 & 0 & 0 & 1 \\ 0 & 0 & -1 & 0 \end{bmatrix} \begin{bmatrix} f_0' \\ f_1' \\ f_2' \\ f_3' \end{bmatrix} = 0 \Rightarrow f_1' \cdot f_0 - f_1 \cdot f_0' = f_2' \cdot f_3 - f_2 \cdot f_3'.$$

Em particular, o exemplo (2.2) é produzido tomando-se $f_2 = t^l$ e $f_3 = t^k$.

■

Por essa equação diferencial é fácil perceber que a tangência de uma curva em um plano de contato implica necessariamente a sua osculação. De fato, diferenciando a equação (2.3) obtemos

$$f_1'' \cdot f_0 - f_1 \cdot f_0'' = f_2'' \cdot f_3 - f_2 \cdot f_3'' \tag{2.4}$$

devido à simetria da expressão. Dessa forma, seja $p = f(1, a) = (a_0 : a_1 : a_2 : a_3)$ ponto de uma curva de contato e $D_p : a_1 z_0 - a_0 z_1 + a_3 z_2 - a_2 z_3 = 0$ o plano da distribuição de contato no ponto p. Os outros pontos de interseção da curva com o plano $D_p$ são as raízes do polinômio de grau d

$$a_1 f_0(1,t) - a_0 f_1(1,t) + a_3 f_2(1,t) - a_2 f_3(1,t) = 0.$$

Temos então que $t = a$ é uma raiz de multiplicidade pelo menos 3, devido as equações diferenciais (2.3) e (2.4) acima. Assim

**Corolário 2.3.2.** *Uma curva de contato em $\mathbb{P}^3$ é osculadora à distribuição de contato.*

Em geral, osculação é o contato de ordem máxima que podemos garantir para todo ponto de uma curva de contato. Como exemplo, para a quártica de contato dada por (2.2) com $k = 3$ e $l = 1$, o polinômio da interseção entre a curva e o plano da distribuição de contato tem fatoração $(t-a)^3(t+a)$. Observe ainda que nesse caso teremos hiperosculação para os pontos $(1:0)$ e $(0:1)$ de $\mathbb{P}^1$.

Relações dessa descrição com os mapas estáveis são explicitadas na próxima seção.

## 2.4 Curvas de contato em famílias de mapas estáveis

Retomamos momentaneamente a estrutura de contato em $\mathbb{P}^{2r-1}$. Uma outra descrição para a estrutura pode ser dada em termos de fibrados através do morfismo (ainda denotado por $\omega$)

$$\begin{array}{rccl} \omega: & T\mathbb{P}^{2r-1} & \to & \mathcal{O}_{\mathbb{P}^{2r-1}}(2) \\ & (p,x) & \mapsto & \omega_p(x). \end{array}$$

Se $\mathcal{D}$ denota o núcleo dessa mapa, a fibra sobre $p$ pode ser identificada com o (espaço tangente do) hiperplano da distribuição de contato em $p$. Obtemos a sequência exata abaixo:

$$\mathcal{D} \rightarrowtail T\mathbb{P}^{2r-1} \xrightarrow{\omega} \mathcal{O}_{\mathbb{P}^{2r-1}}(2)\,. \tag{2.5}$$

Esse fibrado está associado ao fibrado de correlação nula, como definido em [5], §7.

Suponha que temos uma curva racional lisa $X$ descrita por uma parametrização $f : \mathbb{P}^1 \to \mathbb{P}^{2r-1}$. A condição de contato então significa que para cada $p = f(t) \in C$, a forma simplética $\omega_p$ deve se anular na derivada da parametrização $f$. Mais formalmente,

$$f^*\omega \circ df = 0.$$

Essa informação está contida no diagrama abaixo, proveniente da sequência exata (2.5):

$$\begin{array}{ccccc} f^*\mathcal{D} & \rightarrowtail & f^*T\mathbb{P}^{2r-1} & \xrightarrow{f^*\omega} & f^*\mathcal{O}_{\mathbb{P}^{2r-1}}(2), \\ & \nwarrow & \uparrow{\scriptstyle df} & \nearrow & \\ & & T\mathbb{P}^1 & & \end{array} \tag{2.6}$$

onde a condição de contato se traduz na fatoração de $df$ por $f^*\mathcal{D}$, que por sua vez implica na anulação de $T\mathbb{P}^1 \to f^*\mathcal{O}_{\mathbb{P}^{2r-1}}(2)$.

Observe que podemos estender essa caracterização facilmente para $f : C \to \mathbb{P}^{2r-1}$, com $C$ árvore de retas projetivas. Apesar de não termos um fibrado tangente definido nesse caso, por conta dos nós, podemos substituir pelo dual do dualizante $\omega_C$ (vide [20], cap.3A). De fato, $\omega_C$ é um fibrado sobre C e sobre um ponto liso coincide com o cotangente no ponto. Como a definição (2.1.1) se restringe à parte não singular da curva, a condição geral de contato é a anulação de $\omega_C^\vee \to f^*\mathcal{O}_{\mathbb{P}^{2r-1}}(2)$.

Nosso objetivo a partir de agora será descrever espaços que parametrizam as curvas de contato em $\mathbb{P}^3$ ($r = 2$), sem marcas, estratificando pelo grau da curva, para em seguida, obter a dimensão desses espaços e um invariantes de interseção a retas como os adiantados na seção (2.2). Por mera simplificação de notação, denotaremos por $\overline{\mathcal{M}}_d$ e $\overline{M}_d$ o *stack* e o espaço de *moduli* associados a esses mapas estáveis e $\overline{\mathcal{M}}_{d,n}$ e $\overline{M}_{d,n}$ para os casos com marcas[1].

Um curva de grau d em $\mathbb{P}^3$ será identificada a partir de sua parametrização $f : C \to \mathbb{P}^3$ como mapa estável, representante da classe $[f] \in \overline{M}_d$. Vamos então descrever como a condição de contato se traduz em termos dos mapas estáveis.

Considere a família de mapas estáveis (sem marcas)

$$\begin{array}{ccc} \mathcal{C} & \xrightarrow{\varphi} & \mathbb{P}^3 \\ \downarrow{\scriptstyle \pi} & & \\ B & & \end{array} \tag{2.7}$$

Um membro dessa família é uma curva de contato se, para o parâmetro $b \in B$ correspondente, temos um diagrama como em (2.6). Para identificar esses parâmetros, tome o dualizante relativo de $\pi$, $\omega_\pi \to \mathcal{C}$, e o mapa $\omega_\pi^\vee \to \varphi^*T\mathbb{P}^3$. Esse mapa é tal que, restrito à fibra de $\pi$ sobre b parametrizando mapa estável livre de automorfismos, é identificado com a derivada do mapa estável $\varphi_b$, ou seja, com $d\varphi_b : T_bC \to \varphi_b^*T\mathbb{P}^3$ sobre a parte lisa da curva. Assim, temos o diagrama

$$\begin{array}{ccccc} \varphi^*\mathcal{D} & \rightarrowtail & \varphi^*T\mathbb{P}^3 & \twoheadrightarrow & \varphi^*\mathcal{O}_{\mathbb{P}^3}(2), \\ & \nwarrow & \uparrow & \nearrow{\scriptstyle s} & \\ & & \omega_\pi^\vee, & & \end{array} \tag{2.8}$$

no qual agora, a condição de contato de $\varphi_b$ é traduzida pela anulação de s ao longo da fibra de $\pi$ sobre b. De outro modo, considerando a

[1]Em particular, para $d = 0$ a notação $\overline{M}_{0,n}$ representa o espaço de curvas racionais estáveis, coincidindo com a notação usual, com a ressalva que na notação original o "0" se remete ao gênero da curva e não ao grau.

seção

$$\mathcal{O} \to \omega_\pi \otimes \varphi^* \mathcal{O}_{\mathbb{P}^3}(2),$$

o lugar das curvas de contato nessa família será o esquema de zeros da seção adjunta, ainda denotada por $s$ (vide [2], (2.3))

$$\mathcal{O} \to \pi_*(\omega_\pi \otimes \varphi^* \mathcal{O}_{\mathbb{P}^3}(2)).$$

Devemos verificar que $\pi_*(\omega_\pi \otimes \varphi^* \mathcal{O}_{\mathbb{P}^3}(2))$ é de fato um fibrado vetorial. Por definição (vide [21], cap. 8), temos que a fibra sobre um mapa estável $(C; f)$ é o espaço $H^0(\pi^{-1}(C; f), \omega_\pi \otimes \varphi^* \mathcal{O}_{\mathbb{P}^3}(2)|_{\pi^{-1}(C;f)})$ que por sua vez é identificado com $H^0(C, \omega_C \otimes f^* \mathcal{O}_{\mathbb{P}^3}(2))$, pela observação na seção (1.2.1). Assim, para verificar que temos um fibrado vetorial resta verificar que $H^1(C, \omega_C \otimes f^* \mathcal{O}_{\mathbb{P}^3}(2)) = 0$. Para simplificar a notação, a partir de agora escreveremos $\omega_C(k) := \omega_C \otimes f^* \mathcal{O}_{\mathbb{P}^3}(k)$ sobre a curva $C$.

**Lema 2.4.1.** $H^1(C, \omega_C(2)) = 0$*, para* $C$ *árvore de retas projetivas.*

*Prova:* Usaremos indução sobre o número de componentes da curva. Primeiramente, para $C \simeq \mathbb{P}^1$, da dualidade de Serre temos $H^1(\mathbb{P}^1, \omega_{\mathbb{P}^1}(2d)) = H^0(\mathbb{P}^1, \mathcal{O}_{\mathbb{P}^1}(-2d))^\vee = 0$, $\forall d > 0$.

Para o passo indutivo, considere um mapa estável $(C; f)$ com $C$ árvore de retas projetivas e normalização $C' \cup \mathbb{P}^1 \to C$. A decomposição é tomada de modo que existam pontos $P \in C'$, $Q \in \mathbb{P}^1$ mapeados sobre um nó $N \in C$. Considerando as restrições de $f$ a cada componente, seja $d_1$ o grau de $f$ restrito à componente $\mathbb{P}^1$. Assim, de acordo com [21], II.1, podemos construir a sequência exata de 'excisão' para a curva nodal

$$\cdots \longrightarrow H^1(C', \omega_{C'}(2)) \longrightarrow H^1(C, \omega_C(2)) \longrightarrow H^1(\mathbb{P}^1, \omega_{\mathbb{P}^1}(2d_1 + N)) \longrightarrow \cdots,$$

onde, por hipótese de indução vale $H^1(C', \omega_{C'}(2)) = 0$ e, por dualidade de Serre, $H^1(\mathbb{P}^1, \omega_{\mathbb{P}^1}(2d_1 + N)) = H^1(\mathbb{P}^1, \mathcal{O}_{\mathbb{P}^1}(-2d_1 - N))^\vee = 0$, de onde podemos concluir que $H^1(C, \omega_C(2)) = 0$. ■

Vamos explorar as fibras no caso livre de automorfismos. Se considerarmos $C \cong \mathbb{P}^1$, obtemos $H^0(\mathbb{P}^1, \omega_{\mathbb{P}^1} \otimes f^* \mathcal{O}_{\mathbb{P}^3}(2)) = H^0(\mathbb{P}^1, \mathcal{O}_{\mathbb{P}^1}(2d-2))$, pois $d$ é o grau de $f$ e $\omega_{\mathbb{P}^1} = \mathcal{O}_{\mathbb{P}^1}(-2)$, sendo portanto um espaço vetorial de dimensão $2d - 1$. Portanto, ao menos no aberto $M_d^*$, temos que a dimensão é $\dim M_d^* - (2d - 1) = 2d + 1$. De fato, essa é a dimensão esperada para nosso espaço de parâmetros.

Para lidar com casos com automorfismos, uma vez tendo sua descrição em famílias, podemos tomar essas construções sobre *stack* de mapas estáveis $\overline{\mathcal{M}}_d$, como será descrito na próxima seção. Antes, vamos retomar à seção (2.3), onde descrevemos uma equação diferencial para parametrizações de curvas de contato. Observe que tal descrição independe de reparametrizações $\widetilde{f} = f \circ \psi$; como consequência independe do representante de uma classe de famílias $[f] \in M_d$. Além do mais, como a condição de contato se restringe a cada componente, temos que o teorema (2.3.1) pode ser estendido para mapas estáveis em $\overline{M}_d$, com curva domínio sendo árvore de retas projetivas.

De forma mais precisa, se $W_d := W(\mathbb{P}^3, d)$, considere a família universal de parametrizações

$$\begin{array}{ccc} W_d \times \mathbb{P}^1 & \xrightarrow{\varphi} & \mathbb{P}^3 \\ \downarrow \pi & & \\ W_d & & \end{array}$$

onde o mapa $\varphi$ é uma avaliação $(f, p) \mapsto f(p)$ e $\pi$ é a projeção na primeira coordenada. Temos também o mapa classificante $W_d \to M_d \subset \overline{M}_d$.

Como as curvas de domínio nesse caso são simplesmente $\mathbb{P}^1$, podemos tomar a construção de uma seção como construída em (2.8) usando o tangente relativo $T_\pi$ e sua condição de anulação: $\varphi^*\omega \circ df(p) = 0$, para todo $p \in \mathbb{P}^1$. A construção da equação diferencial no teorema (2.3.1) descreve então localmente, ao menos como conjunto, o lugar de zeros dessa seção.

A independência dessas equações, assim como as componentes e suas multiplicidades do nosso espaço de parâmetros podem ser verificadas nos casos particulares de cúbicas e quárticas por meio de computações com SINGULAR, o qual apresentaremos nas seções 3.3 e 3.4.

## 2.5 Um invariante virtual para as curvas de contato

Para uma construção mais geral, vamos descrever o espaço de parâmetros procurado em termos de *stacks* de mapas estáveis.

Consideramos agora a família universal (em *stacks*)

$$\begin{array}{ccc} \overline{\mathcal{M}}_{d,1} & \xrightarrow{\nu} & \mathbb{P}^3 \\ \big\downarrow{\scriptstyle \pi} & & \\ \overline{\mathcal{M}}_d & & \end{array}$$

Definimos $\mathcal{E}_d$ o fibrado vetorial de posto $2d-1$ sobre $\overline{\mathcal{M}}_d$, tal que que, para cada família como em (2.7) construímos o fibrado $\pi_*(\omega_\pi \otimes \varphi^*\mathcal{O}_{\mathbb{P}^3}(2))$ compatível com mudanças de bases (de acordo com a definição no apêndice A.5). Escrevemos $\mathcal{E}_d = \pi_*(\omega_\pi \otimes \nu^*\mathcal{O}_{\mathbb{P}^3}(2))$. Esse é de fato um fibrado vetorial pelo lema (2.4.1), cuja fibra sobre um mapa estável $(C; f)$ é $H^0(C, \omega_C \otimes f^*\mathcal{O}_{\mathbb{P}^3}(2))$.

Temos então uma seção $\bar{s}$ de $\mathcal{E}_d$ associada a $s$ em (2.8), através de mudança de base pelo do mapa de moduli (1.3). Definimos $\mathcal{L}_d = \mathcal{Z}(\bar{s})$, lugar de zeros da seção. Seja $\mathcal{L}_d \to L_d$ o espaço de *moduli* associado.

Como *substack*, $\mathcal{L}_d$ é associado ao funtor $\mathcal{L}_d$ que, a cada esquema $B$, associa as classes de isomorfismo de famílias de mapas estáveis sobre $B$ para as quais a seção $s$ em (2.8) se anula ao longo das fibras.

Seguindo a construção apresentada no apêndice A.6, podemos definir a **classe fundamental virtual** de $\mathcal{L}_d$. Em nosso caso, se $\iota : \mathcal{L}_d \hookrightarrow \overline{\mathcal{M}}_d$, temos que

$$\iota_*([\mathcal{L}_d]^{\mathrm{virt}}) = c_{2d-1}(\mathcal{E}_d) \cap [\overline{\mathcal{M}}_d] \in A_{2d+1}(\overline{\mathcal{M}}_d)_{\mathbb{Q}},$$

uma vez que $\overline{\mathcal{M}}_d$ tem dimensão $4d$ e $\mathcal{E}_d$ tem posto $2d-1$. Em particular, a **dimensão esperada** de $\mathcal{L}_d$ é $2d+1$.

A construção de $\mathcal{L}_d$ feita acima é central em nosso estudo. Vamos enunciar esse resultado através do

**Teorema 2.5.1.** *O espaço das curvas de contato racionais de grau* $d$ *em* $\mathbb{P}^3$ *pode ser parametrizado por um* stack $\iota : \mathcal{L}_d \hookrightarrow \overline{\mathcal{M}}_{0,0}(\mathbb{P}^3, d)$, *definido como lugar de zeros de uma seção do fibrado vetorial* $\mathcal{E}_d = \pi_*(\omega_\pi \otimes \nu^*\mathcal{O}_{\mathbb{P}^3}(2))$. *Em particular,* $\mathcal{L}_d$ *tem dimensão esperada* $2d+1$ *e sua classe fundamental virtual é tal que*

$$\iota_*([\mathcal{L}_d]^{\mathrm{virt}}) = c_{2d-1}(\mathcal{E}_d) \cap [\overline{\mathcal{M}}_d] \in A_{2d+1}(\overline{\mathcal{M}}_d)_{\mathbb{Q}}.$$

A dimensão é apenas dita esperada, no sentido discutido em [14], como um limite inferior para a dimensão ao longo do esquema. Ainda devemos verificar as componentes irredutíveis desse espaço, suas dimensões e multiplicidades, para que possamos tirar conclusões com significado enumerativo. Nos casos específicos tratados nesse trabalho, com $d = 1, 2, 3, 4$, verificamos alguns desses aspectos.

Pela discussão na seção 2.2, vamos definir um invariante (virtual) natural para nosso estudo:

**Definição 2.5.1.** *O invariante virtual relacionado às curvas de contato de grau* $d$ *em* $\mathbb{P}^3$ *incidentes a* $2d+1$ *retas em posição geral de* $\mathbb{P}^3$ *é dado por*

$$N_d := \int_{\mathcal{L}_d} H^{2d+1} = \int_{\overline{\mathcal{M}}_d} c_{2d-1}(\mathcal{E}_d) \cup H^{2d+1}, \tag{2.9}$$

*onde* $H$ *é o divisor de incidência a reta em* $\mathbb{P}^3$*.*

Assim, onde pudermos garantir o significado enumerativo desse invariante, $N_d$ conta exatamente o número de curvas de contato de grau $d$ em $\mathbb{P}^3$ parametrizadas por mapas estáveis e incidentes a $2d+1$ retas em posição geral de $\mathbb{P}^3$. Nos casos de retas e cônicas de contato podemos verificar que isso acontece.

### 2.5.1 De volta às curvas de contato planas

Vamos retomar retas e cônicas de contato da seção (2.2) à luz da descrição por *stacks* de mapas estáveis.

1. Retas de contato: sendo o próprio $\overline{M}_1$ isomorfo à grassmanniana $G(2, 4)$, nossa descrição segue em acordo com o já apresentado: $L_1$ coincide com a grassmanniana lagrangeana, sendo lisa de dimensão pura 3 em $\overline{M}_1$. Calculamos $N_1 = 2$, com o significado enumerativo devido: o número de retas de contato incidentes a 3 retas em posição geral no espaço.

2. Cônicas de contato: pelo apresentado anteriormente, temos que $\dim L_2 = 5$. Porém, nesse caso é interessante diferenciarmos dois tipos: o par de retas concorrentes distintas, ou uma reta dupla. Diferenciar a situação onde aparecem recobrimentos duplos é um dos pontos chaves para entender a estrutura de $L_2$ - e também para graus maiores.

Ao menos como conjunto, podemos identificar duas componentes de $L_2$ que iremos denotar por partições de $d = 2$ - como será melhor explicado na seção (3.3):

- Para (2): temos o conjunto dos recobrimentos duplos $f : \mathbb{P}^1 \to \mathbb{P}^3$. Considere a fatoração $f = h \circ g$, onde $g : \mathbb{P}^1 \to \mathbb{P}^3 \in L_1$ e $h : \mathbb{P}^1 \to \mathbb{P}^1 \in \overline{M}_{0,0}(\mathbb{P}^1, 2)$. Sabendo que $\dim L_1 = 3$ e $\dim \overline{M}_{0,0}(\mathbb{P}^1, 2) = 2$, calculamos dimensão 5. Essa fatoração revela a presença de automorfismo não trivial desses mapas estáveis. Observe que há automorfismo também no caso limite no bordo de $\overline{M}_{0,0}(\mathbb{P}^1, 2)$, em que a curva domínio do mapa estável é redutível.

- Para (1+1): o conjunto dos mapas estáveis $f : \ell_1 \cup \ell_2 \to \mathbb{P}^3$ formando pares de retas de contato, que está contido no bordo $\Delta(1,1)$ de $\overline{M}_2$. A partir do mapa $\delta : \overline{M}_{1,1} \times_{\mathbb{P}^3} \overline{M}_{1,1} \to \Delta(1,1)$ podemos verificar que a dimensão dessa componente é 5 (observe que $\delta$ não é um isomorfismo nesse caso, mas genericamente $2:1$).

  Essa componente também inclui o caso em que as duas retas de contato coincidem (pela diagonal em $L_1 \times L_1$). Isso revela a presença de automorfismo não trivial nos mapas estáveis desse limite, não apenas pela simetria em $\delta$, mas pela imagem desses mapas estáveis.

Porém, observe que mesmo contando com a dimensão esperada, a componente (2) não contribui no cálculo de $N_2$. De fato, sendo a imagem dos mapas estáveis nessa componente uma reta de contato, o produto com o fator $H^5$, a condição de incidência às 5 retas em posição geral, é nula. Com isso, podemos considerar apenas (1+1), concluindo que $N_2 = 40$, pelos cálculos anteriores, com significado enumerativo devido.

Na seção (3.2) recuperaremos os valores de $N_1$ e $N_2$ diretamente da expressão (2.9).

### 2.5.2 Recobrimentos e ramificações

Como indicado no exemplo das cônicas de contato na seção anterior, recobrimentos são casos em que naturalmente aparecem automorfismos. Um outro caso que requer cuidados são mapas reduzidos (irredutíveis?) mas com ramificação. Faremos agora uma análise mais geral sobre a estrutura das ramificações e recobrimentos.

Iniciamos com um resultado que garante que a dimensão do espaço para curvas de contato sem ramificação coincide com a dimensão esperada.

**Teorema 2.5.2** (Bryant). *Seja* $f : \mathbb{P}^1 \to C \subset \mathbb{P}^3$ *curva de contato não ramificada de grau* $d$, *cuja restrição do fibrado da estrutura de contato* *é* $\mathcal{O}_C(-2d)$. *O espaço de* moduli *para essas curvas é um subespaço liso de dimensão* $2d+1$.

*Prova:* Veja [9], proposição 5.1, para um enunciado mais geral, do qual nosso caso se restringe a "$k = 2d - 1$". ■

Esse resultado ainda pode ser estendido para curvas de domínio $C$ nodais.

**Corolário 2.5.3.** *Seja* $f : C \to \mathbb{P}^3$ *uma curva de contato de grau* $d$ *com* $C$ *árvore de retas projetivas tal que a restrição de* $f$ *a cada componente satisfaz condições como no teorema. O espaço para essas curvas tem dimensão* $2d+1$.

*Prova:* Procedemos por indução no número de componentes de $C$. Sem perda de generalidade, tome $f : C = C' \cup \mathbb{P}^1 \to \mathbb{P}^3$, com $d'$ e $d_1$ os graus de $f|_{C'}$ e $f|_{\mathbb{P}^1}$, respectivamente. Por hipótese de indução, a restrição a $f|_{C'}$ está num espaço de dimensão igual a $2d'+1$ e, pelo caso base, a componente $\mathbb{P}^1$ em dimensão $2d_1+1$. Considerando a condição de incidência pelo nó entre $C$ e $\mathbb{P}^1$, temos então dimensão $(2d'+1)+(2d_1+1)-1 = 2d+1$. ■

Para recobrimentos, podemos verificar o seguinte resultado.

**Proposição 2.5.1.** *Seja* $f : C \to \mathbb{P}^3$ *uma curva de contato de grau* $d \geq 2$ *com [imagem não ramificada e] alguma componente sendo um recobrimento de grau* $a > 1$ *(e nenhum outro tipo de ramificação). Então a dimensão de seu espaço de parâmetros é menor ou igual a* $2d+1$, *com a igualdade ocorrendo para recobrimento de retas.*

*Prova:* Consideremos inicialmente o caso em que o mapa estável $f : \mathbb{P}^1 \to \mathbb{P}^3$ de grau $d$ se fatora como $f = h \circ g$, com $g : \mathbb{P}^1 \to \mathbb{P}^1 \in \overline{M}_{0,0}(\mathbb{P}^1, a)$, $a > 1$ e $h : \mathbb{P}^1 \to \mathbb{P}^3$ mapa estável não ramificado de grau $b$, com imagem reduzida, sendo portanto $d = a \cdot b$. Contando dimensões, temos que os mapas $h$ são parametrizados por um espaço de dimensão $2b+1$, por (2.5.2), enquanto que $\dim \overline{M}_{0,0}(\mathbb{P}^1, a) = 2a-2$. Assim, mapas como $f$ estão parametrizados por um espaço de dimensão $2(a+b)-1$. Porém,

$$2(a+b)-1 \leq 2d+1$$

com igualdade ocorrendo apenas quando $b = 1$, ou seja, um recobrimento múltiplo de grau $d$ de uma reta.

Para $f : C \to \mathbb{P}^3$, C árvore de retas projetivas, a mesma conclusão segue por indução no número de componentes de C, com notação como no corolário (2.5.3). Por hipótese de indução, a restrição a $f|_{C'}$ está num espaço de dimensão menor ou igual a $2d' + 1$: se é não ramificado conta com a dimensão exata, se é recobrimento, conta com dimensão menor ou igual a $2d' + 1$. Análogo vale para a componente $\mathbb{P}^1$, em dimensão menor ou igual a $2d_1+1$. Considerando a condição de incidência pelo nó entre as duas componentes, temos então dimensão máxima $(2d' + 1) + (2d_1 + 1) - 1 = 2d + 1$.

∎

Segue dessa descrição:

**Corolário 2.5.2.** *Curvas de contato como descritas na proposição (2.5.1) não contribuem enumerativamente em* $N_d$.

*Prova:*
Nos casos da proposição, apenas recobrimentos múltiplos de grau d de retas de contato contam com a dimensão esperada $2d + 1$ e portanto os demais casos não contribuem no cálculo de $N_d$, por questão de dimensão. Porém, os recobrimentos de retas de contato ainda não contribuem no cálculo de $N_d$, dessa vez devido a condição de incidência a $2d + 1$ retas em posição geral.

∎

Dessa forma, além dos mapas $\mathbb{P}^1 \to \mathbb{P}^3$ não ramificados de grau d, a proposição anterior afirma que os mapas $C = \bigcup \ell_i \to \mathbb{P}^3$, com as restrições a cada $\ell_i$ sendo mapas não ramificados, também não contribuem no cálculo de $N_d$. Em particular, os resultados apresentados garantem que $L_3$, o espaço das cúbicas racionais de contato, tem de fato a dimensão esperada. O caso de uma cúbica cuspidal remete-nos a cúbicas planas; as de contato consistem na união de retas. Com a presença de ramificações, ainda não temos caracterização da dimensão para $d \geq 4$.

# Capítulo 3

# Cálculo dos invariantes virtuais

A partir da fórmula de localização em $\overline{\mathcal{M}}_d$ e a caracterização do espaço das curvas de contato $\mathcal{L}_d$, somos capazes de calcular $N_d$ através de fórmulas gerais. Explicitamos essas fórmulas na seção (3.1), descrevendo as classes equivariantes da condição de incidência a retas e da condição de contato, além de exemplos diretos de aplicação. Até esse ponto de estudo, apesar dessas fórmulas serem explícitas, os cálculos associados demandam intrincadas construções combinatórias, o que limitam o encontro de uma fórmula geral fechada ou recursiva, ou ao menos mais eficiente para $N_d$. Optamos então em tratar alguns casos iniciais. Na seção (3.2) aplicamos a fórmula obtida na subseção anterior nos casos já conhecidos de retas e cônicas de contato. Os números inéditos são obtidos nas seções seguintes: em (3.3) obtemos o invariante $N_3$ para cúbicas de contato, assim como uma descrição geométrica do significado desse invariante; em (3.4) o mesmo para quárticas. Finalmente, em (3.5) fazemos considerações para o problema do cálculo em graus maiores com as ferramentas aqui apresentadas.

As demais referências são as básicas no estudo de curvas e cohomologia, como [21], [14], [20].

## 3.1 As classes equivariantes

Podemos determinar o invariante virtual $N_d$ através da fórmula de localização para *stacks* de mapas estáveis, 1.4.1. Seja $T = (\mathbb{C}^*)^4$ um toro agindo sobre $\mathbb{P}^3$, induzindo uma ação em $\overline{M}_d$. Então, pela fórmula de localização, com inclusões $\iota : \mathcal{M}_\Gamma \hookrightarrow \overline{\mathcal{M}}_d$, a integral dada em (2.9) é

$$\int_{\overline{\mathcal{M}}_d} c_{2d-1}(\mathcal{E}_d) \cup H^{2d+1} = \sum_\Gamma \frac{1}{a(\Gamma)} \int_{\mathcal{M}_\Gamma} \frac{\iota_\Gamma^*(c_{2d-1}^T(\mathcal{E}_d) \cup H_T^{2d+1})}{\mathrm{Euler}^T(\mathcal{N}_\Gamma^{virt})}.$$

onde $\Gamma$ varia entre os grafos ponderados coloridos que parametrizam as componentes de mapas estáveis fixos pela ação de $T$, $\mathcal{N}_\Gamma^{virt}$ seu fibrado normal virtual e $a(\Gamma)$ a ordem do grupo de automorfismos de uma mapa estável associado ao grafo $\Gamma$.

Para cada grafo $\Gamma$, a classe de Euler no denominador é determinada pelos fatores $V(\Gamma)$ e $E(\Gamma)$, dados pelas expressões (1.5) e (1.6). O que nos resta é determinar as classes equivariantes no numerador, correspondentes à condição de contato, dada por $c_{2d-1}(\mathcal{E}_d)$, e às condições de incidência a retas, $H^{2d+1}$.

Resumindo o que será obtido nas seções posteriores: para cada $\Gamma$ parametrizando as componentes $\overline{M}_\Gamma$ do lugar de pontos fixos e $\lambda_i$ os pesos associados a ação de $T$ sobre $\mathbb{P}^3$, temos

$$H_T(\Gamma) := \sum_{v \in V} \mu(v)\lambda_i,$$

$$c_{2d-1}^T(\mathcal{E}_d)(\Gamma) := \prod_{e \in E} \prod_{\alpha=1}^{2d_e-1} \frac{\alpha\lambda_i + (d_e - \alpha)\lambda_j}{d_e} \cdot \prod_{v \in V} (2\lambda_{i_v})^{\mathrm{val}\, v - 1}.$$

Assim, temos o seguinte resultado:

**Proposição 3.1.1.** *Uma fórmula aberta para* $N_d$ *é dada por*

$$N_d = \sum_\Gamma \frac{c_{2d-1}^T(\mathcal{E}_d)(\Gamma) \cdot H_T(\Gamma)^{2d+1}}{a(\Gamma)\,\mathrm{Euler}^T(\mathcal{N}_\Gamma^{virt})}, \tag{3.1}$$

### 3.1.1 Incidência a retas

Seja $h = c_1(\nu^*\mathcal{O}_{\mathbb{P}^3}(1))$ a classe da incidência de um mapa estável a um plano em posição geral de $\mathbb{P}^3$ no ponto marcado. A classe em $\overline{\mathcal{M}}_d$ para a incidência de uma mapa estável a uma reta em posição geral de $\mathbb{P}^3$, em um ponto qualquer, pode então ser definida por $H := \pi_*(h^2)$. Ainda podemos escrever $H$ de modo mais conveniente, para efeito de cálculos, através do teorema de Grothendieck-Riemann-Roch (GRR) [14].

**Lema 3.1.1.** *A classe* $H$ *da incidência de um mapa estável a uma reta em* $\mathbb{P}^3$ *pode ser escrita como*

$$H = \pi_*(h^2) = c_1(\pi_*(\mathcal{O}(2h))) - 2c_1(\pi_*(\mathcal{O}(h))).$$

*Prova:* Considere o fibrado linear $\mathcal{O}(ah)$ em $\overline{\mathcal{M}}_{d,1}$, com $a$ um inteiro não-negativo. Então, o fibrado $\pi_*\mathcal{O}(ah)$ em $\overline{\mathcal{M}}_d$ tem sua fibra sobre o mapa estável $(C; f)$ identificada com $H^0(C, \mathcal{O}(ah)|_C) = H^0(C, f^*\mathcal{O}_{\mathbb{P}^3}(a))$. O objetivo é escrever $H$ em função de $\pi_*\mathcal{O}(ah)$ para valores convenientes de $a$. Para isso, começamos explicitando o GRR para nosso caso, com o mapa de esquecimento e o fibrado $\mathcal{O}(ah)$:

$$\operatorname{ch}(\pi_!\mathcal{O}(ah)) = \pi_*(\operatorname{ch}(\mathcal{O}(ah)) \cdot \operatorname{td}(T_\pi)).$$

Podemos fazer uma primeira simplificação nessa expressão, a partir do seguinte

**Fato 1.** $H^i(C, \mathcal{O}_C(ah)) = 0$ *para todos* $i \geq 1$, $a \geq 0$ *e* $C$ *árvore de retas projetivas.*

*Prova:* Primeiramente, temos a afirmação válida para qualquer $i > 1$, por questão de dimensão. Assim, basta calcular para $i = 1$, o que faremos por indução no número de componentes de $C$.

Se temos mapa estável $(\mathbb{P}^1; f)$ de grau $d$, é direto da dualidade de Serre que

$$H^1(\mathbb{P}^1, \mathcal{O}_{\mathbb{P}^1}(ad)) = H^0(\mathbb{P}^1, \mathcal{O}_{\mathbb{P}^1}(-ad-2))^\vee = 0.$$

Assim como no lema (2.4.1), considere um mapa estável $(C; f)$ com $C$ árvore de retas projetivas e decomposição $C' \cup \mathbb{P}^1 \to C$, tal que existam pontos $P \in C'$, $Q \in \mathbb{P}^1$ mapeados sobre um nó $N \in C$. Seja $d_1$ o grau de $f$ restrito à componente $\mathbb{P}^1$. Assim, podemos construir a sequência exata de 'excisão' para a curva nodal (veja o apêndice B.1 para construção geral)

$$\mathcal{O}_C(ad_1)_{\mathbb{P}^1} \otimes \mathcal{O}_C(-N) \rightarrowtail \mathcal{O}_C(ah) \twoheadrightarrow \mathcal{O}_C(ah)_{C'} \tag{3.2}$$

e a sequência exata longa de cohomologia

$$\cdots \longrightarrow H^1(\mathbb{P}^1, \mathcal{O}_{\mathbb{P}^1}(ad_1 - Q)) \longrightarrow H^1(C, \mathcal{O}_C(ah)) \longrightarrow H^1(C', \mathcal{O}_{C'}(ah)) \longrightarrow 0.$$

Novamente pela dualidade de Serre, temos $H^1(\mathbb{P}^1, \mathcal{O}_{\mathbb{P}^1}(ad_1 - Q)) = H^0(\mathbb{P}^1, \mathcal{O}_{\mathbb{P}^1}(-ad_1 - 2 + Q)) = 0$. Também, pela hipótese de indução, vale $H^1(C', \mathcal{O}_{C'}(ah)) = 0$, já que $C'$ tem menos componentes que $C$. Logo, $H^1(C, \mathcal{O}_C(ah)) = 0$. ∎

Mediante o fato acima, como $H^1(C, \mathcal{O}_C(ah)) = 0$, temos que $R^i\pi_*$ se anulam para $i \geq 1$, e portanto obtemos

$$ch(\pi_*\mathcal{O}(ah)) = \pi_*(ch(\mathcal{O}(ah)) \cdot td(T_\pi))$$

O produto do lado direito da equação acima pode ser calculado explicitamente como (vide [20], 3.E)

$$\left(1 + ah + \frac{a^2h^2}{2}\right)\left(1 - \frac{\kappa_\pi}{2} + \frac{\kappa_\pi^2 + \eta_\pi}{12} + \cdots\right),$$

onde $\kappa_\pi = c_1(\omega_\pi)$ e $\eta_\pi = c_2(\Omega_\pi)$, a classe do lugar de nós nas fibras de $\pi$.

Mas $\pi_*\left(\frac{\kappa_\pi^2 + \eta_\pi}{12} + \cdots\right) = 0$, verificado pelo cálculo de $ch(\pi_*\mathcal{O})$:

$$ch(\pi_*\mathcal{O}) = \pi_*(ch(\mathcal{O}) \cdot td(T_\pi)) = \pi_*\left(1 - \frac{\kappa_\pi}{2} + \frac{\kappa_\pi^2 + \eta_\pi}{12} + \cdots\right)$$

e portanto $0 = c_1(\pi_*\mathcal{O}) = \pi_*\left(\frac{\kappa_\pi^2 + \eta_\pi}{12}\right)$, assim como para as classes de ordem superior.

Com isso, concluímos que

$$c_1(\pi_*(\mathcal{O}(ah))) = \pi_*\left(\frac{a^2h^2}{2} - \frac{ah\kappa_\pi}{2}\right).$$

Comparando a expressão acima para $a = 1$ e $a = 2$ obtemos

$$H = \pi_*(h^2) = c_1(\pi_*(\mathcal{O}(2h))) - 2c_1(\pi_*(\mathcal{O}(h))).$$

■

Com isso, a parte de incidência do numerador em 3.1 se resume a uma soma de pesos da ação induzida sobre $H^0(C, f^*\mathcal{O}_{\mathbb{P}^3}(1))$ e $H^0(C, f^*\mathcal{O}_{\mathbb{P}^3}(2))$.

**Exemplo 1.1.1**: Como ilustração, vamos obter o número de cônicas de $\mathbb{P}^3$ incidentes a 8 retas em posição geral, o qual sabemos ser 92 ([13] ou [15], exemplo 3.2.22). Corresponde ao invariante de Gromov-Witten,

$$I(h^2 \cdots h^2) = \int_{\overline{\mathcal{M}}_2} H^8$$

que por localização é calculado como

$$\sum_\Gamma \frac{1}{a(\Gamma)} \int_{\overline{\mathcal{M}}_\Gamma} \frac{i_\Gamma^*(H_T^8)}{\mathrm{Euler}^T(\mathcal{N}_\Gamma^{virt})}.$$

Primeiramente, determinamos os tipos combinatórios dos grafos $\Gamma$ que parametrizam as componentes de pontos fixos em $\overline{\mathcal{M}}_2$ pela ação usual. Eles são mostrados na tabela abaixo, junto com o número de automorfismos e classes de colorações.

| | | |
|---|---|---|
| $i \overset{2}{\text{---}} j$ | $i \overset{1}{\text{---}} j \overset{1}{\text{---}} k$ | $i \overset{1}{\text{---}} j \overset{1}{\text{---}} i$ |
| $a(\Gamma) = 2$ | $a(\Gamma) = 1$ | $a(\Gamma) = 2$ |
| $6\times$ | $12\times$ | $12\times$ |

Apenas nesse exemplo, vamos determinar os pesos usando as sequências exatas usadas nessa seção.

- Recobrimento duplo:

  Para esse caso obtemos os pesos diretamente, com $C = \mathbb{P}^1$ e $f$ mapa de grau 2. Se $s, t$ são coordenadas homogêneas para $\mathbb{P}^1$ e $z_0, \cdots, z_3$ coordenadas homogêneas para $\mathbb{P}^3$, $f$ pode ser tomada por $z_i = s^2$ e $z_j = t^2$. Assim, se $\lambda_i$ é o peso para a coordenada $z_i$, temos que os pesos para $s$ e $t$ são $\lambda_i/2$ e $\lambda_j/2$ respectivamente. Portanto,

  Para $H^0(\mathbb{P}^1, \mathcal{O}_{\mathbb{P}^1}(2))$, com geradores $s^2, st, t^2$, os pesos são $\frac{2\lambda_i}{2}$, $\frac{\lambda_i + \lambda_j}{2}$ e $\frac{2\lambda_j}{2}$;

  Para $H^0(\mathbb{P}^1, \mathcal{O}_{\mathbb{P}^1}(4))$, com geradores $s^4, s^3t, s^2t^2, st^3, t^4$, os pesos são $\frac{4\lambda_i}{2}$, $\frac{3\lambda_i + \lambda_j}{2}$, $\frac{2\lambda_i + 2\lambda_j}{2}$, $\frac{\lambda_i + 3\lambda_j}{2}$ e $\frac{4\lambda_j}{2}$

  e portanto, pelo lema (3.1.1), $H = 2\lambda_i + 2\lambda_j$.

- Duas retas simples:

  Agora, sejam $X = \mathbb{P}^1$ e $Y = \mathbb{P}^1$ as duas componentes de $C$, com nó $N$. Como coordenadas, sejam $s, t$ para $X$ e $u, v$ para $Y$, de modo que $N$ seja a imagem de $\mathcal{V}(s)$ e $\mathcal{V}(u)$ pela normalização. Suponha também que $f(N) = q_j$. Como o grau do mapa $f$ é linear restrito a qualquer componente, temos $z_i = s$ e $z_j = t$ para $X$; $z_k = u$ e $z_j = v$ para $Y$. Portanto, $s$ tem peso $\lambda_i$, $t$ e $v$ têm peso $\lambda_j$ e $u$ tem peso $\lambda_k$. Tomando cohomologias da sequência exata (3.2) aplicada nesse caso, temos

$$H^0(X, \mathcal{O}_X(a - N)) \rightarrowtail H^0(C, \mathcal{O}_C(ah)) \twoheadrightarrow H^0(Y, \mathcal{O}_Y(a)).$$

Uma base para o espaço da direita são os monômios em $u, v$ de grau $a$. Já para o espaço da esquerda, temos os monômios em $s, t$ de grau $a$ que se anulam em $\mathcal{V}(s)$, ou seja, múltiplos de $s$.

Para $a = 1$, os pesos para esses espaços são respectivamente $\lambda_j, \lambda_k$ e $\lambda_i$. Para $a = 2$, os pesos são $2\lambda_j, \lambda_j + \lambda_k, 2\lambda_k$ e $\lambda_i + \lambda_j, 2\lambda_i$. Dessa forma, a contribuição da classe $H$ é $\lambda_i + 2\lambda_j + \lambda_k$.

Observe que todo os cálculos são análogos para o caso $i = k$, resultando em $2\lambda_i + 2\lambda_j$.

Por fim, a soma das 30 funções racionais nas variáveis $\lambda_0, \cdots, \lambda_3$ correspondentes aos 30 pontos fixos (isolados) resulta em... 92. Cálculos podem ser verificados através do algoritmo apresentado no apêndice E.

□

De modo geral, vale a seguinte fórmula:

**Proposição 3.1.2.** *A classe equivariante para a condição de incidência a retas na componente $\mathcal{M}_\Gamma$ é dada por*

$$H_T(\Gamma) = \sum_{v \in V} \mu(v)\lambda_{i_v},$$

*ou seja, cada vértice $v \in \Gamma$ contribui com uma parcela constituída pelo seu peso $\lambda_{i_v}$ contado com multiplicidade $\mu(v)$, o somatório dos graus das arestas incidentes a $v$.*

*Prova:* O caso base pode ser tomado como no caso de recobrimento do exemplo (3.1.1), assim como sua notação, resultando em $d(\lambda_i+\lambda_j)$. Dado um mapa estável $(C; f) \in \overline{M}_d$, seguimos por indução, considerando que o resultado é válido para qualquer mapa estável com menor número de componentes na curva domínio que em C.

Tome uma decomposição para a curva C, como no lema (2.4.1), do tipo $C' \cup \mathbb{P}^1$. Temos $C'$ árvore de retas projetivas, conexa, e $\mathbb{P}^1$ um galho de C. Seja N o nó da interseção entre $C'$ e $\mathbb{P}^1$. Seja $d_1$ o grau do mapa restrito à componente $\mathbb{P}^1$ e $d' = d - d_1$ o grau do mapa restrito a $C'$. Tome $(s : t)$ coordenadas para $\mathbb{P}^1$, com $V(s) \mapsto N$. Se $\mathbb{P}^1$ é mapeada na reta coordenada $q_i q_j$, com N mapeado em $q_j$, então $z_i = s^{d_1}$ e $z_j = t^{d_1}$. Assim, os pesos associados às variáveis $s$ e $t$ são respectivamente $\lambda_i/d_1$ e $\lambda_j/d_1$.

Pelo lema (3.1.1), precisamos então calcular explicitamente os pesos para $H^0(\mathcal{O}_C(h))$ e $H^0(\mathcal{O}_C(2h))$. Para isso, considere a sequência exata

$$\mathcal{O}_{\mathbb{P}^1}(d_1 - N) \rightarrowtail \mathcal{O}_C(h) \twoheadrightarrow \mathcal{O}_{C'}(h).$$

Pela dualidade de Serre, temos $H^1(\mathbb{P}^1, \mathcal{O}_{\mathbb{P}^1}(d_1 - N)) = H^0(\mathbb{P}^1, \omega_{\mathbb{P}^1}(-d_1 + N))^\vee = H^0(\mathbb{P}^1, \mathcal{O}_{\mathbb{P}^1}(-2 - d_1 + N))^\vee = 0$ para qualquer $d_1 \geq 0$. Desse modo, obtemos uma sequência exata de cohomologias

$$H^0(\mathbb{P}^1, \mathcal{O}_{\mathbb{P}^1}(d_1 - N)) \rightarrowtail H^0(C, \mathcal{O}_C(h)) \twoheadrightarrow H^0(C', \mathcal{O}_{C'}(h)).$$

Os pesos de $H^0(\mathbb{P}^1, \mathcal{O}_{\mathbb{P}^1}(d_1 - N))$ são determinados pela base formada pelos monômios de grau $d_1$ em $s, t$ múltiplos de $s$

$$s^\alpha t^{d_1 - \alpha},\ 1 \leq \alpha \leq d_1$$

e portanto associado aos pesos

$$\frac{\alpha\lambda_i + (d_1 - \alpha)\lambda_j}{d_1},\ 1 \leq \alpha \leq d_1.$$

Assim, para $c_1(\pi_* \mathcal{O}_{\mathbb{P}^1}(d_1 - N))$ temos

$$\begin{aligned} \sum_{\alpha=1}^{d_1} \frac{\alpha\lambda_i + (d_1 - \alpha)\lambda_j}{d_1} &= \frac{1}{d_1}\left(\left(\frac{d_1(d_1+1)}{2}\right)\lambda_i + \left(\frac{d_1(d_1-1)}{2}\right)\lambda_j\right) = \\ &= \frac{d_1+1}{2}\lambda_i + \frac{d_1-1}{2}\lambda_j. \end{aligned}$$

Analogamente, podemos obter uma sequência exata de cohomologias

$$H^0(\mathbb{P}^1, \mathcal{O}_{\mathbb{P}^1}(2d_1 - N)) \rightarrowtail H^0(C, \mathcal{O}_C(2h)) \twoheadrightarrow H^0(C', \mathcal{O}_{C'}(2h)).$$

O peso para $c_1(\pi_* \mathcal{O}_{\mathbb{P}^1}(2d_1 - N))$ é dado por

$$\sum_{\alpha=1}^{2d_1} \frac{\alpha\lambda_i + (2d_1 - \alpha)\lambda_j}{d_1} = (2d+1)\lambda_i + (2d-1)\lambda_j.$$

Com isso, temos que o peso para a classe H será

$$
\begin{aligned}
H &= \sum \text{pesos de } H^0(\mathcal{O}_C(2h)) - 2\sum \text{pesos de } H^0(\mathcal{O}_C(h)) = \\
&= \Big[\sum \text{pesos de } H^0(\mathcal{O}_{C'}(2h)) + \sum \text{pesos de } H^0(\mathcal{O}_{\mathbb{P}^1}(2d_1 - N))\Big] + \\
&\quad -2\Big[\sum \text{pesos de } H^0(\mathcal{O}_{C'}(h)) + \sum \text{pesos de } H^0(\mathcal{O}_{\mathbb{P}^1}(d_1 - N))\Big] = \\
&= \Big[\sum \text{pesos de } H^0(\mathcal{O}_{C'}(2h)) - 2\sum (\text{pesos de } H^0(\mathcal{O}_{C'}(h))\Big] + \\
&\quad \Big[\sum \text{pesos de } H^0(\mathcal{O}_{\mathbb{P}^1}(2d_1 - N)) - 2\sum \text{pesos de } H^0(\mathcal{O}_{\mathbb{P}^1}(d_1 - N))\Big] = \\
&= \sum_{v \in \Gamma_{C'}} \mu(v)\lambda_{i_v} + \Big[(2d_1+1)\lambda_i + (2d_1-1)\lambda_j - (d_1+1)\lambda_i - (d_1-1)\lambda_j\Big] = \\
&= \sum_{v \in \Gamma_{C'}} \mu(v)\lambda_{i_v} + [d_1\lambda_i + d_1\lambda_j] = \sum_{v \in \Gamma_C} \mu(v)\lambda_{i_v},
\end{aligned}
$$

ou seja, uma vez que a contribuição da componente $C'$ provém da hipótese de indução, a contribuição da componente $\mathbb{P}^1$ consiste em adicionar parcelas da forma $d_1\lambda_{i_v}$ correspondente a seus vértices. ■

### 3.1.2 Condição de contato

Vamos obter agora os pesos correspondentes à classe de Euler do fibrado que descreve a condição de contato em termos de mapas estáveis, $\mathcal{E}_d$. Lembramos que, por simplificação, denotamos $\omega_C(k) := \omega_C \otimes f^*\mathcal{O}_{\mathbb{P}^3}(k)$.

**Exemplo 1.2.2**: Como ilustração, descrevemos os pesos para $\mathcal{E}_1$. Como a fibra de $\mathcal{E}_1$ sobre $(C, f)$ é $H^0(C, \omega_C(2)) = H^0(\mathbb{P}^1, \omega_{\mathbb{P}^1}(2))$, vamos usar a sequência de Euler, tensorizada por $\mathcal{O}_{\mathbb{P}^1}(2)$

$$\omega_{\mathbb{P}^1}(2) \rightarrowtail \mathcal{O}_{\mathbb{P}^1}(1) \otimes H^0(\mathbb{P}^1, \mathcal{O}_{\mathbb{P}^1}(1)) \twoheadrightarrow \mathcal{O}_{\mathbb{P}^1}(2).$$

para obter os pesos. Tomando cohomologias, teremos uma sequência exata com $H^0$, uma vez que, pela dualidade de Serre, $H^1(\mathbb{P}^1, \omega_{\mathbb{P}^1}(2)) = H^0(\mathbb{P}^1, \mathcal{O}_{\mathbb{P}^1}(-2))^\vee = 0$:

$$H^0(\mathbb{P}^1, \omega_{\mathbb{P}^1}(2)) \rightarrowtail H^0(\mathbb{P}^1, \mathcal{O}_{\mathbb{P}^1}(1)) \otimes H^0(\mathbb{P}^1, \mathcal{O}_{\mathbb{P}^1}(1)) \twoheadrightarrow H^0(\mathbb{P}^1, \mathcal{O}_{\mathbb{P}^1}(2)).$$

Como os geradores para o espaço da direita são $s^2, st, t^2$ e os do centro são $s \otimes s, s \otimes t, t \otimes s, t \otimes t$, o mapa sobrejetor identifica $s \otimes s \mapsto s^2, s \otimes t \mapsto$ $st, t\otimes t \mapsto t^2$. Assim, a classe equivariante para $c_1(\mathcal{E}_1)$, correspondente a $t\otimes s$, será então $\lambda_i + \lambda_j$. □

Para o caso de cônicas, ainda podemos proceder de modo análogo.

**Exemplo 1.2.3**: Retomemos a notação e os mapas fixos descritos no exemplo (3.1.1). Vamos obter, como antes, os pesos correspondentes $\mathcal{E}_2$ para cada uma dessas configurações.

- Recobrimento duplo:

  Devemos obter os pesos de $H^0(C, \omega_C(2))) = H^0(\mathbb{P}^1, \omega_{\mathbb{P}^1}(4))$. Novamente, a partir da sequência de Euler temos

$$H^0(\mathbb{P}^1, \omega_{\mathbb{P}^1}(4)) \rightarrowtail H^0(\mathbb{P}^1, \mathcal{O}_{\mathbb{P}^1}(3)) \otimes H^0(\mathbb{P}^1, \mathcal{O}_{\mathbb{P}^1}(1)) \twoheadrightarrow H^0(\mathbb{P}^1, \mathcal{O}_{\mathbb{P}^1}(4)).$$

  Usando coordenadas $(s:t)$, temos que os pesos para $H^0(\mathbb{P}^1, \omega_{\mathbb{P}^1}(4))$ correspondem aos geradores da forma

$$s^{3-\alpha}t^{\alpha}\otimes s \quad \text{e} \quad s^{3-\alpha}t^{\alpha}\otimes t, \quad \text{com} \quad 0\le \alpha\le 3$$

  Podemos tomá-los como sendo $s^{3-\alpha}t^{\alpha}\otimes t$, $0\le\alpha\le 2$, nos dando $(3\lambda_i+\lambda_j)/2, (2\lambda_i+2\lambda_j)/2, (\lambda_i+3\lambda_j)/2$. Assim, a classe de Euler equivariante nesse caso é

$$\frac{3\lambda_i+\lambda_j}{2}\cdot\frac{2\lambda_i+2\lambda_j}{2}\cdot\frac{\lambda_i+3\lambda_j}{2}.$$

- Duas retas simples:

  Ainda nas considerações do exemplo anterior, a partir de uma sequência exata de excisão como em (3.2) aplicada a esse caso obtemos a sequência exata de cohomologia

$$H^0(X, \omega_{\mathbb{P}^1}(2)) \rightarrowtail H^0(C, \omega_C(2)) \twoheadrightarrow H^0(Y, \omega_{\mathbb{P}^1}(2+N)).$$

  que novamente é exata pela dualidade de Serre. Os pesos procurados então serão os pesos provenientes dos dois espaços das extremidades, correspondentes a cada componente de C.

  Para a componente $\mathbb{P}^1$ do espaço à esquerda, tomamos a partir da sequência de Euler obtendo

$$\omega_{\mathbb{P}^1}(2) \rightarrowtail \mathcal{O}_{\mathbb{P}^1}(1)\otimes H^0(\mathbb{P}^1, \mathcal{O}_{\mathbb{P}^1}(1)) \twoheadrightarrow \mathcal{O}_{\mathbb{P}^1}(2).$$

e produzimos a sequência, novamente exata, nas cohomologias

$$H^0(\mathbb{P}^1, \omega_{\mathbb{P}^1}(2)) \rightarrowtail H^0(\mathbb{P}^1, \mathcal{O}_{\mathbb{P}^1}(1)) \otimes H^0(\mathbb{P}^1, \mathcal{O}_{\mathbb{P}^1}(1)) \twoheadrightarrow H^0(\mathbb{P}^1, \mathcal{O}_{\mathbb{P}^1}(2)).$$

Geradores para o espaço da direita são os monômios $s^2, st, t^2$. Já para o espaço do centro os geradores são $s \otimes s, s \otimes t, t \otimes s, t \otimes t$. Assim, pela sobrejeção temos $s \otimes s \mapsto s^2$, $s \otimes t \mapsto st$ e $t \otimes t \mapsto t^2$. Portanto um gerador para o espaço da esquerda é $t \otimes s$ cujo peso relacionado é $\lambda_j + \lambda_i$.

Agora, para o espaço $H^0(\mathbb{P}^1, \omega_{\mathbb{P}^1}(2+N))$, tomamos novamente a sequência de Euler

$$\omega_{\mathbb{P}^1}(2+N) \rightarrowtail \mathcal{O}_{\mathbb{P}^1}(1+N) \otimes H^0(\mathbb{P}^1, \mathcal{O}_{\mathbb{P}^1}(1)) \twoheadrightarrow \mathcal{O}_{\mathbb{P}^1}(2+N).$$

e produzimos a sequência exata nas cohomologias

$$H^0(\mathbb{P}^1, \omega_{\mathbb{P}^1}(2+N)) \rightarrowtail H^0(\mathbb{P}^1, \mathcal{O}_{\mathbb{P}^1}(1+N)) \otimes H^0(\mathbb{P}^1, \mathcal{O}_{\mathbb{P}^1}(1)) \twoheadrightarrow H^0(\mathbb{P}^1, \mathcal{O}_{\mathbb{P}^1}(2+N)$$

Os geradores para o espaço da direita correspondentes às funções racionais de grau 2 com pólo em $N = \mathcal{V}(u)$, ou seja

$$\frac{u^3}{u}, \frac{u^2v}{u}, \frac{uv^2}{u}, \frac{v^3}{u}.$$

Analogamente, para o espaço do centro os geradores são

$$\frac{u^2}{u} \otimes u, \ \frac{uv}{u} \otimes u, \ \frac{v^2}{u} \otimes u, \ \frac{u^2}{u} \otimes v, \ \frac{uv}{u} \otimes v, \ \frac{v^2}{u} \otimes v.$$

A sobrejeção identifica

$$\frac{u^2}{u} \otimes u \mapsto \frac{u^3}{u}, \ \frac{uv}{u} \otimes u \mapsto \frac{u^2v}{u}, \ \frac{v^2}{u} \otimes u \mapsto \frac{uv^2}{u}, \ \frac{v^2}{u} \otimes v \mapsto \frac{v^3}{u}.$$

Assim, os geradores para o espaço da esquerda são

$$\frac{u^2}{u} \otimes v, \ \frac{uv}{u} \otimes v$$

cujos peso relacionado são respectivamente $\lambda_k + \lambda_j$ e $2\lambda_j$.

Reunindo os pesos obtidos, concluímos que a classe equivariante para $c_3(\mathcal{E}_2)$ é

$$(\lambda_i + \lambda_j)(\lambda_j + \lambda_k)(2\lambda_j)$$

Para o caso $i = k$ os cálculos são análogos, obtendo o peso $(\lambda_i + \lambda_j)^2(2\lambda_j)$.

□

Os exemplos acima já são suficientes para generalizar a expressão dos pesos para qualquer d. De fato, temos

**Proposição 3.1.3.** *A classe equivariante para a condição de contato na componente $\mathcal{M}_\Gamma$ é dada por*

$$c^T_{2d-1}(\mathcal{E}_d)(\Gamma) = \prod_{e\in E}\prod_{\alpha=1}^{2d_e-1}\frac{\alpha\lambda_i+(2d_e-\alpha)\lambda_j}{d_e}\cdot\prod_{\nu\in V}(2\lambda_{i_\nu})^{\mathrm{val}\,\nu-1}$$

*ou seja, cada vértice $\nu$ de $\Gamma$ contribui com um fator $(2\lambda_{i_\nu})^{\mathrm{val}\,\nu-1}$ (em particular, vértices com apenas uma aresta adjacente não contam) e cada aresta* e *contribui com um fator associado a cada monômio misto de grau $2d_e$ nas coordenadas de sua reta imagem.*

*Prova:* Tomamos a mesma notação do caso de incidência a retas e usamos indução sobre o número de componentes de C, partindo agora da sequência exata

$$\omega_{C'}(2)\rightarrowtail\omega_C(2)\twoheadrightarrow\omega_{\mathbb{P}^1}(2d+N).$$

Para o caso base, procedemos como nos exemplos acima. Pela sequência de Euler dualizada

$$\omega_{\mathbb{P}^1}(2d)\rightarrowtail\mathcal{O}_{\mathbb{P}^1}(2d-1)\otimes H^0(\mathbb{P}^1,\mathcal{O}_{\mathbb{P}^1}(1))\twoheadrightarrow\mathcal{O}_{\mathbb{P}^1}(2d).$$

produzimos a sequência exata nas cohomologias

$$H^0(\mathbb{P}^1,\omega_{\mathbb{P}^1}(2d))\rightarrowtail H^0(\mathbb{P}^1,\mathcal{O}_{\mathbb{P}^1}(2d-1))\otimes H^0(\mathbb{P}^1,\mathcal{O}_{\mathbb{P}^1}(1))\twoheadrightarrow H^0(\mathbb{P}^1,\mathcal{O}_{\mathbb{P}^1}(2d))$$

uma vez, que, pela dualidade de Serre, $H^1(\mathbb{P}^1,\omega_{\mathbb{P}^1}(2d)) = H^1(\mathbb{P}^1,\mathcal{O}_{\mathbb{P}^1}(-2d))^\vee = 0$

Dado que o espaço à direita é formado pelos monômios de grau 2d, seus geradores são da forma

$$s^\alpha t^{2d-\alpha}, 0\leq\alpha\leq 2d.$$

Por sua vez, os geradores do espaço do meio são da forma

$$s^\alpha t^{2d-1-\alpha}\otimes s, s^\alpha t^{2d-1-\alpha}\otimes t,\ 0\leq\alpha\leq 2d-1.$$

O mapa sobrejetor então identifica

$$s^\alpha t^{2d-1-\alpha}\otimes s\mapsto s^{\alpha+1}t^{2d-1-\alpha},\ 0\leq\alpha\leq 2d-1,$$

e

$$t^{2d-1} \otimes t \mapsto t^{2d}.$$

Dessa forma, os geradores para o espaço à esquerda são

$$s^{\alpha}t^{2d-1-\alpha} \otimes t,\ 1 \leq \alpha \leq 2d-1$$

cujos pesos associados são

$$\frac{\alpha\lambda_i + (2d-\alpha)\lambda_j}{d},\ 1 \leq \alpha \leq 2d-1.$$

Agora, por hipótese de indução, considere que vale a fórmula em (3.1.3) de pesos para $C'$. A partir da sequência exata

$$\omega_{C'}(2) \rightarrowtail \omega_C(2) \twoheadrightarrow \omega_{\mathbb{P}^1}(2d_1+N).$$

obtemos uma sequência exata de cohomologias

$$H^0(C', \omega_{C'}(2)) \rightarrowtail H^0(C, \omega_C(2)) \twoheadrightarrow H^0(\mathbb{P}^1, \omega_{\mathbb{P}^1}(2d_1+N)),$$

exatidão garantida pelo lema (2.4.1). Resta então calcular os pesos para o espaço à direita. Para isso, tome a sequência exata

$$\omega_{\mathbb{P}^1}(2d_1+N) \rightarrowtail \mathcal{O}_{\mathbb{P}^1}(2d_1-1+N) \otimes H^0(\mathbb{P}^1, \mathcal{O}_{\mathbb{P}^1}(1)) \twoheadrightarrow \mathcal{O}_{\mathbb{P}^1}(2d_1+N).$$

da qual produzimos a sequência exata nas cohomologias

$$H^0(\mathbb{P}^1, \omega_{\mathbb{P}^1}(2d_1+N)) \rightarrowtail H^0(\mathbb{P}^1, \mathcal{O}_{\mathbb{P}^1}(2d_1-1+N)) \otimes H^0(\mathbb{P}^1, \mathcal{O}_{\mathbb{P}^1}(1)) \twoheadrightarrow H^0(\mathbb{P}^1, \mathcal{O}_{\mathbb{P}^1}(2d_1$$

pois $H^1(\mathbb{P}^1, \omega_{\mathbb{P}^1}(2d_1+N)) = H^1(\mathbb{P}^1, \mathcal{O}_{\mathbb{P}^1}(-2d_1-N))^{\vee} = 0$

Dado que o espaço à direita é formado pelas funções racionais de grau $2d_1$ com pólo simples em $N$, seus geradores são da forma

$$\frac{s^{\alpha}t^{2d_1+1-\alpha}}{s},\ 0 \leq \alpha \leq 2d_1+1.$$

Por sua vez, os geradores do espaço do meio são da forma

$$\frac{s^{\alpha}t^{2d_1-\alpha}}{s} \otimes s, \frac{s^{\alpha}t^{2d_1-\alpha}}{s} \otimes t,\ 0 \leq \alpha \leq 2d_1.$$

O mapa sobrejetor então identifica

$$\frac{s^{2d_1}}{s} \otimes s \mapsto \frac{s^{2d_1+1}}{s}$$

e

$$\frac{s^\alpha t^{2d_1-\alpha}}{s} \otimes t \mapsto \frac{s^\alpha t^{2d_1+1-\alpha}}{s},\ 0 \le \alpha \le 2d_1.$$

Dessa forma, os geradores para o espaço à esquerda são

$$\frac{s^\alpha t^{2d_1-\alpha}}{s} \otimes s,\ 0 \le \alpha \le 2d_1 - 1$$

cujos pesos associados são

$$\frac{\alpha\lambda_i + (2d_1 - \alpha)\lambda_j}{d_1},\ 0 \le \alpha \le 2d_1 - 1.$$

Dessa forma temos, para $\alpha = 0$, o peso $2\lambda_j$, e para $\alpha \neq 0$ os pesos associados aos monômios mistos de grau $2d_1$.

Com esses pesos, podemos determinar os pesos para a classe de Chern:

$$c^T_{2d-1}(\mathcal{E}_d) = c^T_{top}(\omega_{C'}(2)) \cdot c^T_{top}(\omega_{\mathbb{P}^1}(2d_1 + N)).$$

O primeiro fator, por hipótese de indução, contribui com o produto dos pesos associados aos monômios mistos em cada aresta de $C'$ e com o produto dos pesos $(2\lambda_{i_\nu})^{\mathrm{val}\,\nu-1}$ para cada vértice $\nu$. Por sua vez, o segundo fator contribui com o produto dos pesos associados aos monômios mistos de sua única aresta e adiciona um fator $2\lambda_{i_\nu}$ correspondente ao vértice do nó N. Com isso, o produto final dos pesos apresenta a forma em (3.1.3), como queríamos verificar. ∎

## 3.2 Retas e cônicas

A partir dessas expressões somos capazes de obter $N_d$, desde que possamos somar sobre todas árvores ponderadas coloridas. Baseado nos exemplos das seções anteriores e dos algoritmos apresentados no apêndice E vamos calcular $N_1$ e $N_2$ explicitamente, onde a enumeração das árvores é facilmente construída.

Recolhendo as informações dos exemplos (3.1.2):

$$N_1 = \int_{\overline{\mathcal{M}}_1} c_1(\mathcal{E}_1) \cdot H^3 = \sum_{i,j,k,l \in \{0,1,2,3\} \text{distintos}} \frac{(\lambda_i + \lambda_j)^4}{(\lambda_i - \lambda_k)(\lambda_i - \lambda_l)(\lambda_j - \lambda_k)(\lambda_j - \lambda_l)} = 2,$$

concordando com o que já havíamos obtido anteriormente.

A partir dos exemplos (3.1.1) e (3.1.2), obtemos

$$N_2 = \int_{\overline{\mathcal{M}}_2} c_3(\mathcal{E}_2) \cdot H^5 = 40 \tag{3.3}$$

Esse valor também está em concordância com a proposição (2.2.2) para curvas de contato planas.

## 3.3 Cúbicas

Como adiantado, o primeiro caso revelante do nosso problema é obter o invariante para as cúbicas de contato em $\mathbb{P}^3$. Primeiramente, como um 'processo de verificação', usando a fórmula de localização e a tabela no apêndice D somos capazes de obter o número de cúbicas em $\mathbb{P}^3$ incidentes a 12 retas no espaço, já apresentado em [13]:

$$\int_{\overline{\mathcal{M}}_3} H^{12} = 80160.$$

Observe que ainda para $d = 3$ temos apenas pontos fixos isolados, uma vez que para o caso de valência 3 temos a componente de mapas fixos isomorfa a $\overline{M}_{0,3} = \bullet$.

Para calcular o invariante $N_3$ pela fórmula de localização, tomamos novamente a tabela em (D) e aplicamos o algoritmo correspondente no apêndice E. Obtemos assim

$$N_3 = \int_{\overline{\mathcal{M}}_3} c_5(\mathcal{E}_3) \cdot H^7 = 4160. \tag{3.4}$$

Vamos fazer uma análise mais minuciosa do significado desse número, descrevendo $L_3$. Já sabemos que as cúbicas de contato planas não contribuem enumerativamente, pois a dimensão de seu espaço de parâmetros é 6. Dentre as espaciais, a discussão sobre ramificação e recobrimentos da seção (2.5.2), já nos permite identificar, como conjunto, as configurações de mapas estáveis $f : C \to \mathbb{P}^3$ possíveis. Entendemos por **configuração de contato** um tipo das imagens genéricas de tais mapas estáveis, identificando o número de componentes de C e seus graus de restrição. Isso pode ser relacionado com as partições do inteiro $d = 3$:

$$\underbrace{\bigcup C_i}_{C} \mapsto \Big(\underbrace{\sum \deg f|_{C_i}}_{3}\Big).$$

Nesse caso diremos que temos uma configuração $(\sum \deg f|_{C_i})$. Desse modo há apenas duas configurações que contribuem na contagem de $N_3$:

- Para $(3)$: mapas estáveis $\mathbb{P}^1 \to \mathbb{P}^3$, cuja única possibilidade são as cúbicas reversas, automaticamente sem ramificações.
- Para $(1+1+1)$: mapas estáveis $\ell_1 \cup t \cup \ell_2 \to \mathbb{P}^3$, árvore de três ramos $\mathbb{P}^1$ com dois nós. Observe que esse é o único caso redutível válido. De fato, uma configuração $(2+1)$, cônica e reta de contato, reduzem automaticamente para $(1+1+1)$, pela proposição (2.2.1). Também, se as retas se encontram em um único nó, a condição de contato faria as três serem necessariamente contidas no plano da distribuição de contato associado ao nó.

  O conjunto associado à configuração $(1+1+1)$ está contido em $\Delta(1,1,1)$, o 'bordo do bordo' de $\overline{M}_3$, identificado pelo isomorfismo com $(\overline{M}_{1,1} \times_{\mathbb{P}^3} \overline{M}_{1,2}) \times_{\mathbb{P}^3} \overline{M}_{1,1}$. Com isso, temos indicada uma estrutura recursiva nessa componente.

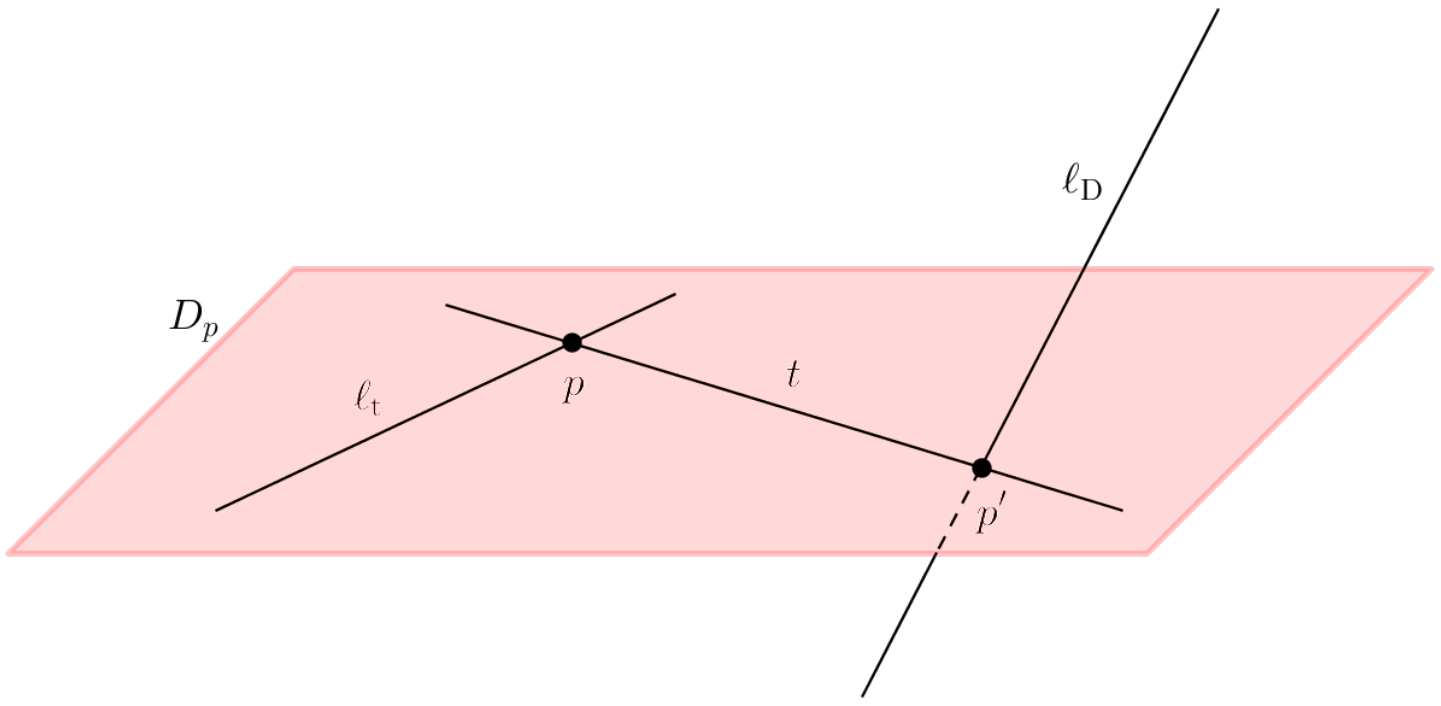

Figura 3.1: Curva típica de configuração $(1+1+1)$, "Z"

Como esquema, confirmamos a estrutura da componente $(3)$ a partir da descrição dada pela equação diferencial 2.3 e com o uso do SINGULAR. A decomposição primária das equações locais mostram duas componentes irredutíveis e reduzidas para mapas $\mathbb{P}^1 \to \mathbb{P}^3$ de codimensão 5: uma referente aos recobrimentos triplos de uma reta de contato, outra referente às cúbicas reversas (veja apêndice E). A primeira, como já observado, apesar de apresentar a dimensão esperada, não contribuiu enumerativamente em $N_3$. A segunda componente corresponde à configuração $(3)$.

Com base nessas considerações, podemos decompor $N_3$ em uma parcela $N_3^{(3)}$ e outra parcela $N_3^{(1,1,1)}$ correspondendo às cúbicas redutíveis de contato, incidentes a 7 retas. Podemos ainda identificar $N_3^{(1,1,1)}$ por configurações de incidência, ao menos de forma combinatória sem considerar possíveis multiplicidades como esquema, utilizando a estrutura recursiva do bordo de $\overline{M}_3$ e o que conhecemos sobre retas e cônicas de contato.

Para contabilizar $N_3^{(1,1,1)}$, fazemos nova observação sobre a natureza das configurações de contato. Considere 7 retas em posição geral no espaço, $h_1, \cdots, h_7$. Impor a incidência da curva $f : C \to \mathbb{P}^3$ a essas 7 retas determina elementos no bordo de $\overline{M}_{0,7}$, o espaço de parâmetro para curvas estáveis racionais com 7 marcas. De fato, de $f : C \to \mathbb{P}^3$, tome $(C, p_1, \cdots, p_7)$ uma curva marcada, onde os pontos $p_i$ são determinados pela interseção da imagem de C com as 7 retas em posição geral no espaço. A estabilidade é garantida pela condição de contato e generalidade: cada componente representa uma curva de contato e portanto tem pelo menos três pontos especiais, entre as marcas $p_i$ e nós com as outras componentes de contato. Isso pois $\dim L_d$ é um limite superior para o número de marcas, por condição de generalidade, sendo o menor valor possível $\dim L_1 = 3$. Dessa forma, a contagem de curvas em cada configuração pode ser associada a partições do conjunto de marcas $\{p_1, \cdots, p_7\}$ em tantas parcelas quanto o número de componentes das curvas analisadas, restritas às condições de contato:

$$\Big(\underbrace{\bigcup C_i}_{C}, \underbrace{\bigcup A_i}_{\{p_i, \cdots, p_7\}}\Big) \mapsto \Big(\Big(\underbrace{\sum |A_i|}_{7}\Big)\Big).$$

Denominamos **configuração de incidência** os mapas estáveis descritos por essa distribuição de marcas. A notação em parênteses duplos é apenas para diferenciar do caso de partição do grau. Observe que essa partição vale inclusive para C irredutível. Além disso, devem ser considerados os automorfismos associados às simetrias dessas configurações.

Com base nisso, vamos construir uma cúbica redutível $C = \ell_1 \cup t \cup \ell_2$, como na figura (3.3), onde t é a componente 'central' na qual $\ell_1$ e $\ell_2$ incidem em pontos necessariamente distintos. Consideramos então 3 casos baseado nas partições possíveis dessas 7 retas dentre as 3 componentes da cúbica, sendo que cada componente pode ser incidente a até 3 dentre as retas $h_i$:

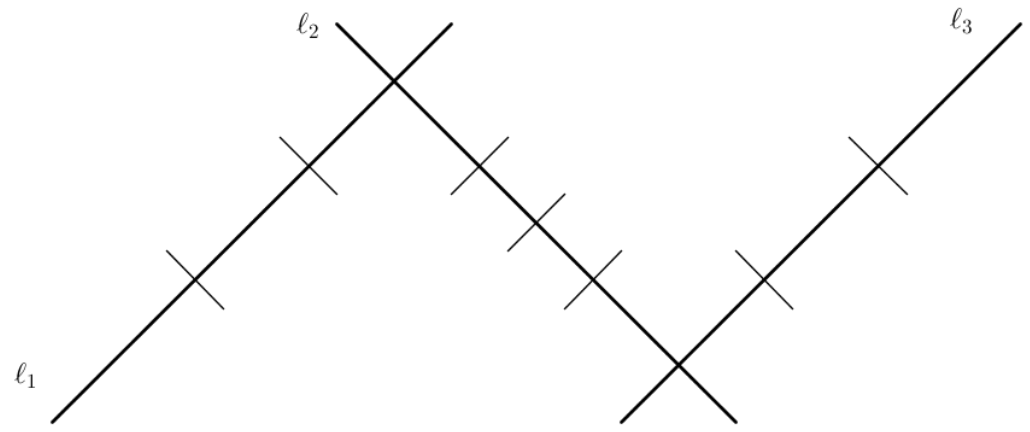

(a) $((2+3+2))$, Aut = 2

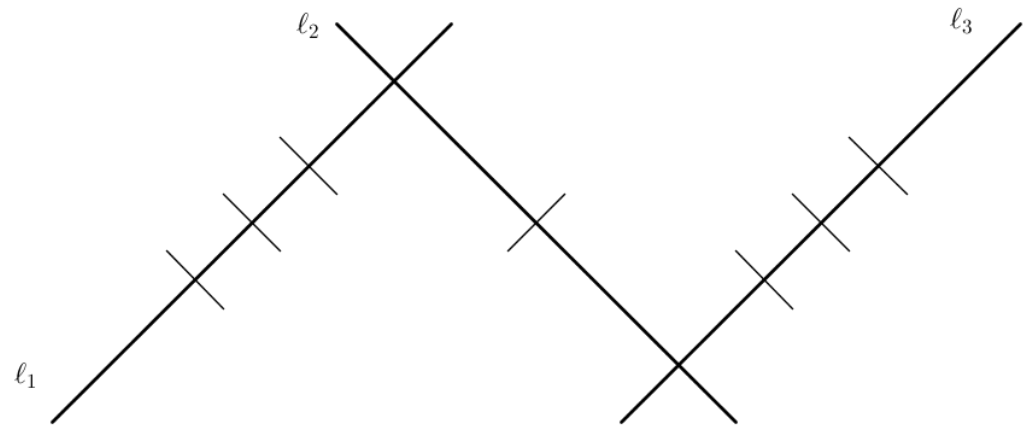

(b) $((3+1+3))$, Aut = 2

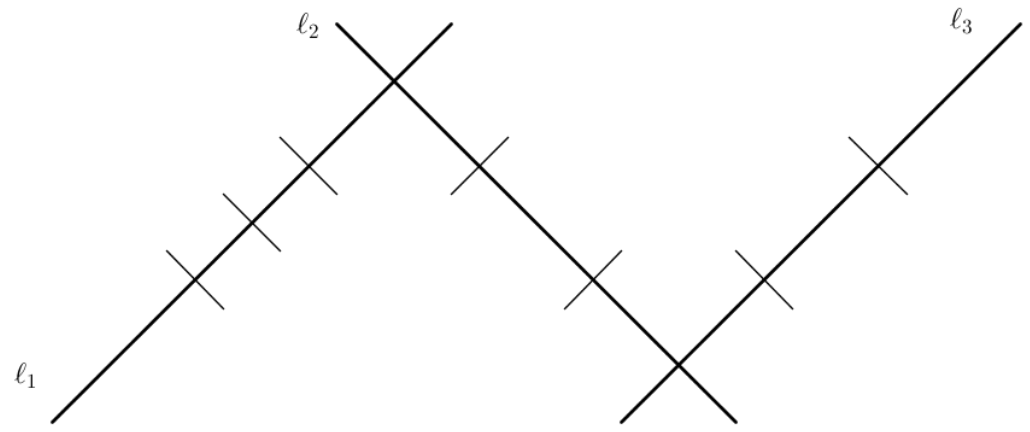

(c) $((3+2+2))$, Aut = 1

Figura 3.2: Distribuição das retas $h_i$ na configuração $(1+1+1)$

- $((3{+}1{+}3))$: Tome 3 dentre as $h_i$ e produza $N_1 = 2$ retas de contato, escolhendo uma delas como $\ell_1$. Tome 3 dentre as 4 $h_i$ restantes e produza outras 2 retas de contato, escolhendo uma delas como $\ell_2$. Com a $h_i$ restante e $\ell_1, \ell_2$, produzimos finalmente mais 2 retas de contato, escolhendo uma delas para ser $t$. Devemos levar em conta a simetria dessa configuração, descontando a duplicidade. Dessa forma, obtemos
$$\frac{1}{2} \cdot \binom{7}{3} \cdot \binom{N_1}{1} \cdot \binom{4}{3} \cdot \binom{N_1}{1} \cdot \binom{1}{1} \cdot \binom{N_1}{1} = 560$$
configurações de 3 retas de contato.

- $((2+3+2))$: Tome 3 dentre as $h_i$ e produza 2 retas de contato, escolhendo uma delas para ser $t$. Junto a ela, tome mais 2 dentre

as 4 $h_i$ restantes produzindo outras 2 retas de contato, escolhendo uma delas para ser $\ell_1$. Novamente, com $t$ e as outras 2 $h_i$ restantes, produza mais um par de reta de contato, escolhendo uma delas para ser $\ell_2$. Novamente devemos descontar a duplicidade da simetria. Obtemos assim mais

$$\frac{1}{2}\cdot\binom{7}{3}\cdot 2\cdot\binom{4}{2}\cdot 2\cdot\binom{2}{2}\cdot 2 = 840$$

configurações de 3 retas de contato.

- $((3+2+2))$: Tome 3 dentre as $h_i$ e produza 2 retas de contato, escolhendo uma como $\ell_1$. Junte a $\ell_1$ mais 2 dentre as 4 $h_i$ restantes e produza outras 2 retas de contato, escolhendo uma como $t$. Com as 2 $h_i$ restantes e $t$, produza finalmente mais 2 retas de contato, escolhendo uma delas para ser $\ell_2$. Dessa forma, obtemos mais

$$\binom{7}{3}\cdot 2\cdot\binom{4}{2}\cdot 2\binom{2}{2}\cdot 2 = 1680$$

configurações de 3 retas de contato.

Observe que as duas últimas configurações de incidência podem ser recuperadas pela construção de cônicas de contato por 5 retas no espaço: $(((3+2)+2)) = ((5+2))$ e $(((2+3)+2)) = ((5+2)) \sim ((2+5)) = ((2+(3+2)))$. Porém, a primeira é inerente ao caso de cúbicas. Por esse motivo, escolhemos enumerar as configurações de incidência a partir de retas de contato, como faremos também para quárticas.

Os casos listados acima contabilizam 3080. Dessa forma, se conjecturarmos que essas configurações contam todas com multiplicidade 1, teríamos

$$N_3^{(3)} = N_3 - N_3^{(1,1,1)} = 1080$$

cúbicas reversas de contato incidentes a 7 retas em posição geral no espaço.

## 3.4 Quárticas

Para as quárticas de contato encontramos componentes de dimensão positiva no cálculo de $N_4$ pela fórmula de localização. De fato, a partir desse caso há sempre grafos de mapas estáveis com valência maior que 3, criando um isomorfismos $\mathcal{M}_\Gamma \simeq \overline{M}_{0,n}$. Independente disso, os cálculos dos pesos abrangem inclusive componentes de dimensão positiva. Sendo assim, o cálculo dos invariantes seguem o mesmo padrão dos anteriores.

A partir da enumeração dos grafos associados às partições do grau $d = 4$, mostrados na tabela do apêndice D, aplicamos à fórmula de localização novamente a um caso de verificação, o número de quárticas incidentes a 16 retas em posição geral no espaço:

$$\int_{\overline{\mathcal{M}}_4} H^{16} = 383306880,$$

número também conhecido e apresentado em [13].

Voltando a nossa caso de interesse, calculamos

$$N_4 = \int_{\overline{\mathcal{M}}_4} c_7(\mathcal{E}_4) \cdot H^9 = 1089024. \tag{3.5}$$

Assim como no caso de cúbicas e usando notação análoga, vamos analisar melhor o conteúdo desse invariante identificando configurações de quárticas pelas partições de $d = 4$.

- Para $(4)$: temos as quárticas dadas por mapas estáveis $\mathbb{P}^1 \to \mathbb{P}^3$ irredutíveis, que incluem quárticas lisas, cuspidais ou nodais.
- Para $(3+1)$: parametrizadas por $c \cup \ell \to \mathbb{P}^3$, curva nodal composta por $c$ cúbica reversa (irredutível) de contato e $\ell$ reta de contato;
- Para $(1+1+1+1)$: com $\ell_1 \cup \ell_2 \cup \ell_3 \cup \ell_4 \to \mathbb{P}^3$, retas de contato com três nós, 3 a 3 não coplanares. As retas devem ser não coplanares 3 a 3, caso contrário já teríamos redução de dimensão pela cúbica plana. A configuração de quatro retas com três ou quatro delas incidentes em um mesmo ponto também reduz ao caso plano. Por fim, a configuração $(2+2)$, par de cônicas de contato, se reduz à configuração $(1+1+1+1)$. Ainda assim temos duas subconfigurações para essa partição, como representadas nas figuras (3.4).

Primeiramente, vamos discutir a componente de mapas $\mathbb{P}^1 \to \mathbb{P}^3$, $(4)$. Novamente usando o SINGULAR, confirmamos a dimensão esperada 9 em $L_4$ para mapas de $W_4$, usando a descrição em 2.3. Propomos também efetuar a decomposição em $W_4$, mas diferente do caso anterior, a decomposição primária não foi computacionalmente efetiva. Porém,

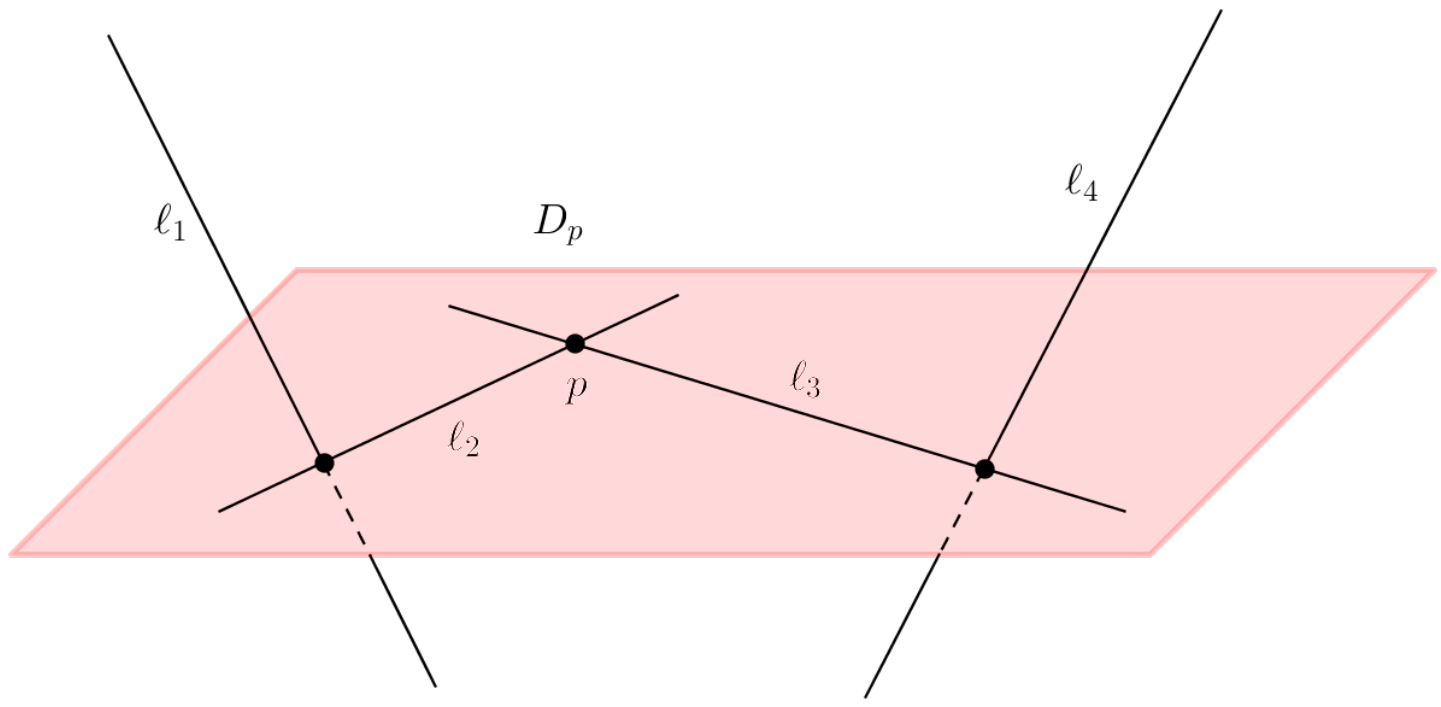

(a) Subconfiguração "W"

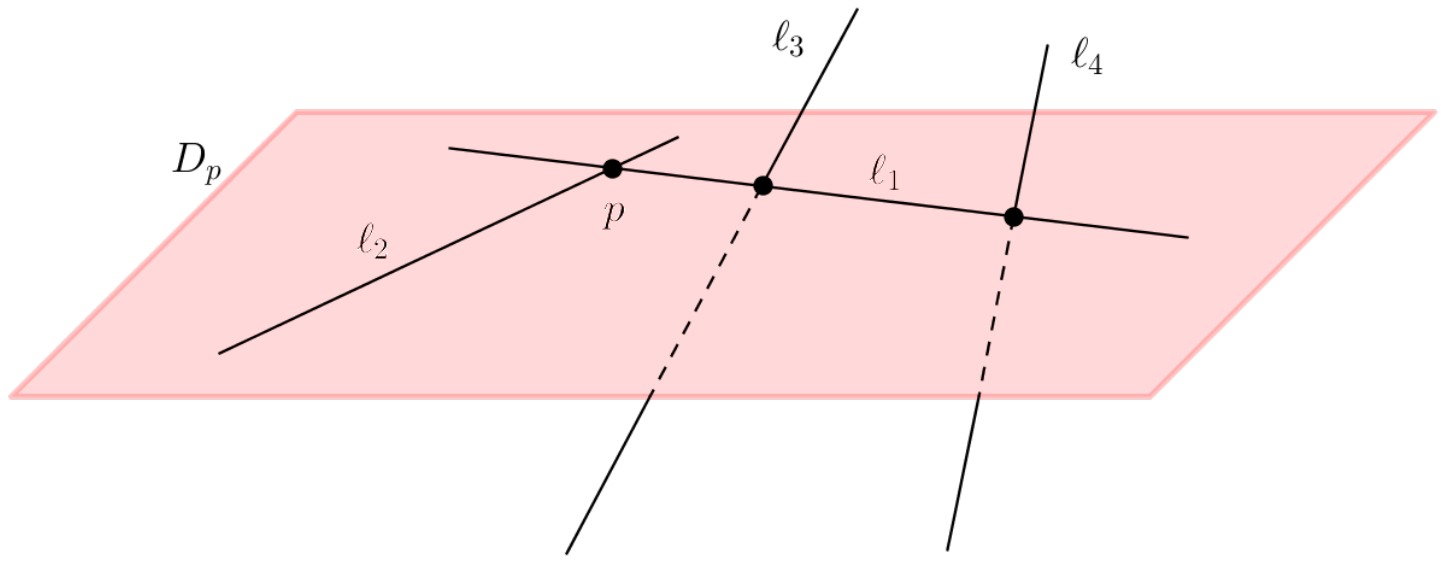

(b) Subconfiguração "w"

Figura 3.3: Subconfigurações de $(1+1+1+1)$

a observação sobre hiperosculação feita na seção (2.3) poderá ser útil agora.

Sabemos que condição de contato implica que a curva oscula a distribuição de contato. Para uma quártica de contato, isso define uma involução sobre seus pontos. De fato, considere a quártica de contato modelo $f(t:s) = (s^4 : \frac{1}{2}t^4 : s^3t : st^3)$, o ponto $q = f(a:b)$ e o plano da distribuição $D_q$. Vimos que o polinômio que determina a interseção da curva com o plano $D_q$ tem fatoração $(bt-as)^3(bt+as)$. Temos então definida uma involução $(a:b) \mapsto (-a:b)$ sobre os pontos de $\mathbb{P}^1$, que associa um ponto $(a:b)$ ao segundo ponto de interseção da curva com o plano. Seus pontos fixos são $(0:1)$ e $(1:0) = \infty$, ou seja, a quártica tem hiperosculação em $D_0 : z_1 = 0$ e $D_\infty : z_0 = 0$.

Temos assim uma alternativa para observar as componente das quárticas de contato $\mathbb{P}^1 \to \mathbb{P}^3$: podemos fixar o plano $D_\infty$ para hiperosculação no ponto $f(1:0)$. Isso nos permite considerar parametrizações afins polinomiais de quárticas, tomando $f_0 = s^4$. Portanto, no aberto afim em que $s \neq 0$ em $\mathbb{P}^1$ e $D_\infty$ considerado como o plano no infinito em $\mathbb{P}^3$, isso nos permite tomar $f(t) = (1 : f_1(t) : f_2(t) : f_3(t))$, onde $f_i(t)$ são polinômios de grau até 4. Por essa redução, a equação diferencial (2.3) toma a forma $f_1' = f_2' f_3 - f_2 f_3'$.

Restrito a esse fechado de $W_4$, realizamos a decomposição primária pelo SINGULAR e obtemos novamente duas componentes, uma delas reduzida referente aos recobrimentos de reta de grau 4. A outra componente, corresponde a quárticas espaciais, não é reduzida. Conjecturamos que isso possa ser por não transversalidade dessa componente quando restrita ao fechado em questão. Porém, ainda é necessário verificar a multiplicidade, além de estender à $W_4$ - aparentemente através de uma involutividade em $\mathbb{P}^3$ para mover o plano hiperosculado e um automorfismo de $\mathbb{P}^1$ que escolha o ponto de pré imagem da hiperosculação.

Retomando as configurações de contato, podemos então decompor

$$N_4 = N_4^{(4)} + N_4^{(3,1)} + N_4^{(1,1,1,1)},$$

onde são possíveis multiplicidades nas componentes de cada parcela.

As parcelas em $N_4$ referentes às configurações $(3+1)$ e $(1+1+1+1)$ podem ser calculadas usando os casos já conhecidos e as configurações de incidência.

Sejam $h_i$, $1 \leq i \leq 9$, nove retas em posição geral no espaço. Podemos associar as quárticas de contato com a condição de incidência às retas com curvas estáveis em $\overline{M}_{0,9}$, análogo ao caso de cúbicas. As possibilidades de interseção de retas $h_i$ em cada componente de contato ficam então restritas à condição de estabilidade das componentes de curvas estáveis em $\overline{M}_{0,9}$: máximo de 3 incidências às retas de contato e 7 às cúbicas de contato.

Para a configuração $(3+1)$, com curva $c \cup \ell$, consideramos a partição das 9 retas $h_i$ na configuração $((6+3))$. Sendo assim, uma construção é: tome 3 dentre as $h_i$ e produza $N_1 = 2$ retas de contato, escolhendo uma delas como $\ell$. Tome as 6 $h_i$ restantes e, juntamente com $\ell$ produza $N_3^{(3)} = 1080$ cúbicas reversas de contato, escolhendo uma delas como $c$. Dessa forma, obtemos

$$\binom{9}{3} \cdot \binom{N_1}{1} \cdot \binom{6}{6} \cdot \binom{N_3^{(3)}}{1} = 181440$$

configurações de cúbica reversa e reta de contato.

Para enumerar as configurações $(1+1+1+1)$ precisamos distribuir 9 marcas dentre as 4 retas de contato, cada uma tendo exatamente 3 pontos de interseção com as retas $h_i$ ou com as outras retas de contato. Isso gera as possibilidades indicadas nas figuras (3.4 e 3.4)

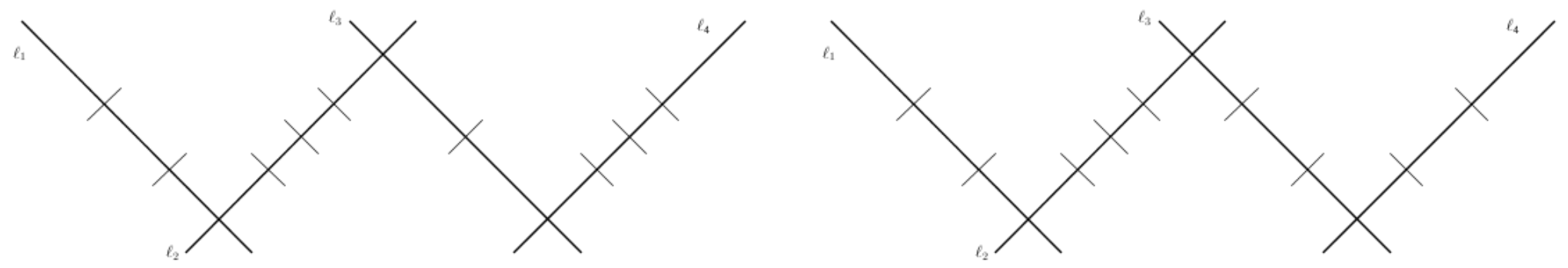

(a) $((2+3+1+3))$, Aut = 1 (b) $((2+3+2+2))$, Aut = 1

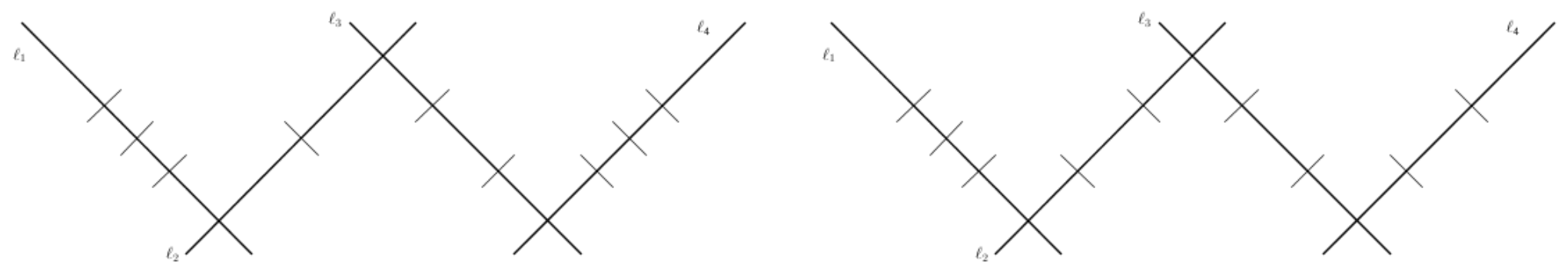

(c) $((3+1+2+3))$, Aut = 1 (d) $((3+2+2+2))$, Aut = 1

Figura 3.4: Partições da subconfiguração "W" de $(1+1+1)$

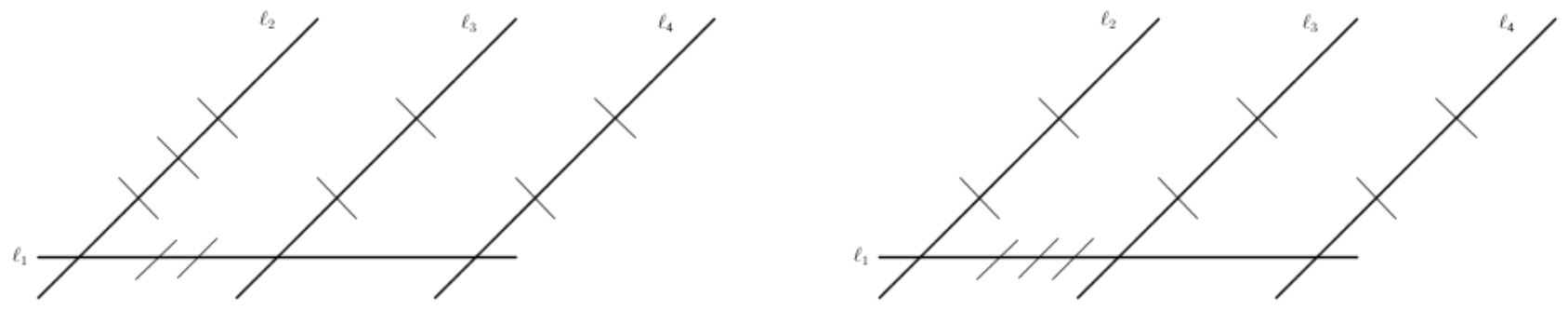

(a) $((2+3+2+2))$, Aut = 2 (b) $((3+2+2+2))$, Aut = 6

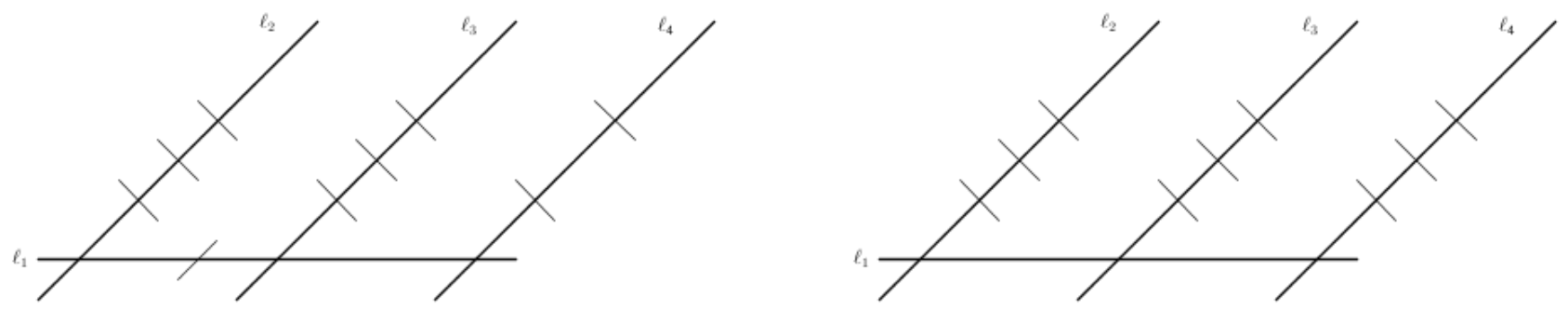

(c) $((1+3+3+2))$, Aut = 2 (d) $((0+3+3+3))$, Aut = 6

Figura 3.5: Partições da subconfiguração "w" de $(1+1+1)$

Observe que temos incluído um caso (o último da lista) onde há uma componente sem interseções com as retas $h_i$, dado que já apresenta interseção com as outras retas de contato. Contabilizando a contribuição de cada partição associada a subconfiguração "W":

- $((3+2+2+2))$: repetindo a construção sequencial das componentes de contato, produzindo $\ell_1$, $\ell_2$, $\ell_3$ e $\ell_4$ nessa ordem, podemos contabilizar

$$\binom{9}{3}\cdot 2\cdot\binom{6}{2}\cdot 2\binom{4}{2}\cdot 2\cdot\binom{2}{2}\cdot 2 = 120960$$

- $((3+1+2+3))$: produzimos inicialmente as retas $\ell_1$ e $\ell_4$, para depois produzir $\ell_3$ e $\ell_2$ e obtemos

$$\binom{9}{3}\cdot 2\cdot\binom{6}{3}\cdot 2\cdot\binom{3}{2}\cdot 2\cdot\binom{1}{1}\cdot 2 = 80640$$

- $((2+3+2+2))$: produzimos inicialmente $\ell_2$, em seguida $\ell_1$, $\ell_3$ e $\ell_4$, obtendo

$$\binom{9}{3}\cdot 2\cdot\binom{6}{2}\cdot 2\binom{4}{2}\cdot 2\cdot\binom{2}{2}\cdot 2 = 120960$$

- $((2+3+1+3))$: produzimos $\ell_2$ e $\ell_4$ para em seguida construir $\ell_1$ e $\ell_3$, contabilizando

$$\binom{9}{3}\cdot 2\cdot\binom{6}{3}\cdot 2\binom{3}{2}\cdot 2\cdot\binom{1}{1}\cdot 2 = 80640.$$

Análogo para a subconfiguração "w":

- $((3+2+2+2))$: produzimos a base $\ell_1$ e em seguida as outras três incidentes à ela, considerando as 3! permutações:

$$\frac{1}{3!}\cdot\binom{9}{3}\cdot 2\cdot\binom{6}{2}\cdot 2\binom{4}{2}\cdot 2\cdot\binom{2}{2}\cdot 2 = 20160$$

- $((2+3+2+2))$: produzimos $\ell_2$, então $\ell_1$. Em seguida, considerando as 2! permutações, produzimos $\ell_3$ e $\ell_4$, obtemos

$$\frac{1}{2!}\cdot\binom{9}{3}\cdot 2\cdot\binom{6}{2}\cdot 2\binom{4}{2}\cdot 2\cdot\binom{2}{2}\cdot 2 = 60480$$

- $((1+3+3+2))$: produzimos $\ell_2$ e $\ell_3$, considerando as 2! permutações, em seguida $\ell_1$ e $\ell_4$, obtendo

$$\frac{1}{2!}\cdot\binom{9}{3}\cdot 2\cdot\binom{6}{3}\cdot 2\binom{3}{1}\cdot 2\cdot\binom{2}{2}\cdot 2 = 40320$$

- $((0+3+3+3))$: produzimos $\ell_2$, $\ell_3$ e $\ell_4$, considerando as 3! permutações, e usamos essas três para em seguida construir a base $\ell_1$, contabilizando

$$\frac{1}{3!}\cdot\binom{9}{3}\cdot 2\cdot\binom{6}{3}\cdot 2\cdot\binom{3}{3}\cdot 2\cdot\binom{0}{0}\cdot 2 = 4480$$

Os números acima totalizam 710080 quárticas redutíveis. Se não houver multiplicidade nas componentes, conjecturamos então $N_4^{(4)} = 378944$.

## 3.5 Desafios para curvas de grau superior

Para curvas de grau $d > 4$, os números $N_d$ ainda podem ser calculados pelas fórmulas apresentadas na seção 3.1. Porém, os algoritmos apresentados no apêndice E ainda não são gerais o suficiente quanto às enumerações, que consistem na geração de árvores ponderadas e colorações, juntamente com suas classes de isomorfismos e automorfismos. Tal listagem torna-se um trabalho enfadonho, mesmo computacionalmente (veja discussão sobre a complexidade do problema de coloração, no apêndice C). Estratégias e algoritmos mais eficientes são citados no próprio artigo de Kontsevich [24]. Também, fórmulas recursivas passam pelo estudo de funções geradoras e cohomologia quântica, como introduzidas em [28] ou [25].

A dimensão de $L_d$ foi confirmada como sendo a esperada para os casos $d = 1, \cdots, 5$ segundo cálculos sobre $W_d$ realizados em SINGULAR, como apresentados no apêndice E. O cálculo para graus superiores foi limitado por recursos computacionais utilizados.

Outro desafio está na descrição do significado enumerativo desses invariantes: é preciso compreender mais da estrutura esquemática e combinatória das componentes de $\mathcal{L}_d$. Conjectura-se que no aberto das curvas de domínio $\mathbb{P}^1$ teremos sempre uma componente correspondente ao recobrimento de grau d de uma reta e outra(s) componente(s) referentes às curvas racionais de grau d não recobertas. Um ponto ainda delicado é entender como as ramificações afetam a dimensão e multiplicidade dessas componentes. Para grau $d = 4$, primeiro caso que permite mapas estáveis com ramificações, observamos apenas uma componente não reduzida, mas ainda restrita a um fechado. Para as componentes dos mapas com curva domínio nodais, as configurações de contato e de incidência apresentadas nos exemplos de cúbicas e quárticas parecem continuar válidas, usando das propriedades do bordo do *stack* de mapas estáveis $\overline{\mathcal{M}}_d$ e o de curvas estáveis $\overline{\mathcal{M}}_{0,2d+1}$. Os cálculos desses exemplos também já sugerem um padrão para fórmulas de contagem em cada configuração, baseada nas partições de grau e marcas.

Contudo, espera-se que atacando ao menos os problemas de implementação, já seremos capazes de expandir a lista de invariantes virtuais $N_d$ conhecidos para mais alguns casos. Nos trabalhos futuros, incluímos ainda a procura por fórmulas recursivas ou funções geradoras para os números $N_d$.

# Apêndice A

# *Stacks* algébricos

Esse apêndice sintetiza, para a conveniência do leitor, construções e propriedades de um *stack* algébrico, baseado principalmente nos textos introdutórios de [29], [1], [19] e também [18]. Outras referências usuais são [7] e [12].

## A.1 O modelo

Dado um esquema $M$ sobre $\mathbb{C}$, temos definido o funtor de pontos de $M$

$$\begin{array}{rcl} h_M := \mathrm{Hom}(\cdot, M): \ (Sch/\mathbb{C}) & \rightarrow & (Set) \\ B & \mapsto & \mathrm{Hom}(B, M) \end{array}$$

Para entender o conceito de *stack*, é importante pensarmos nesse funtor como um 'feixe de conjuntos'. Isso pode ser feito formalmente inserindo uma topologia na categoria $(Sch/\mathbb{C})$, chamada topologia de Grothendieck, através de morfismos de cobertura.

Uma **cobertura** para uma categoria $\mathcal{C}$ é uma família de morfismos $\{f_i : U_i \rightarrow U\}_{i\in I}$ e dizemos que $\mathcal{C}$ tem uma **topologia de Grothendieck** quando

- Se $f$ é um isomorfismo, então $\{f\}$ é uma cobertura;

- Se $\{U_i \rightarrow U\}_i$ é cobertura e $\{U_{ij} \rightarrow U_i\}_j$ é uma cobertura para cada $i$, então $\{U_{ij} \rightarrow U\}_{i,j}$ também é cobertura, por composição.

- Se $\{f_i : U_i \rightarrow U\}$ é uma cobertura e $f_0 : V \rightarrow U$ é um morfismo, então $\{f_{i0,0} : U_i \times_U V \rightarrow V\}$ é uma cobertura.

A categoria $(Sch/\mathbb{C})$ apresenta naturalmente uma topologia de Grothendieck herdada pelos subesquemas abertos e seus morfismos de inclusão.

Quanto esses morfismos são lisos ou *étale*, dizemos que temos uma topologia lisa ou *étale*, respectivamente.

Dizemos que um funtor contravariante do tipo $\mathcal{F} : (Sch/\mathbb{C}) \to (Set)$ é um **pré-feixe de conjuntos**. Para denominarmos tal funtor como um **feixe de conjuntos**, tomamos a partir de uma topologia de Grothendieck em $(Sch/\mathbb{C})$ as seguintes condições válidas para qualquer cobertura $\{f_i : U_i \to U\}$:

- Se $X, Y \in \mathcal{F}(U)$ tais que $\mathcal{F}(f_i)(X) = \mathcal{F}(f_i)(Y)$ para todo $i$, então $X = Y$;
- Se $X_i \in \mathcal{F}(U_i)$, para cada $i$, tal que $\mathcal{F}(f_{ij,i})(X_i) = \mathcal{F}(f_{ij,j})(X_j)$, então existe $X \in \mathcal{F}(U)$ com $\mathcal{F}(f_i)(X) = X_i$ para cada $i$.

O morfismo apresentado na segunda condição acima é $f_{ij,i} : U_i \times_U U_j \to U_i$. Os objetos na imagem $\mathcal{F}(f_i)(X)$ devem ser pensados como uma restrição $X|_i$.

Dadas essas definições, vemos que de fato o funtor $h_M$ é um feixe de conjuntos. Para extrairmos mais informações sobre os objetos parametrizados por esse funtor, como automorfismos, é necessário um pouco mais de estrutura, como descrito a seguir.

## A.2 Problema de *moduli*

Suponha que estamos procurando um espaço de parâmetros para uma certa classe de objetos geométricos baseados em esquemas (curvas sobre um esquema, por exemplo). Isso pode ser decodificado por um pré-feixe de conjuntos

$$\mathcal{F} : (Sch/\mathbb{C}) \to (Set),$$

que envia cada esquema $B$ na classe dos objetos que queremos parametrizar sobre $B$. O estudo da existência e propriedades de tal espaço de parâmetros é o que chamados de **problema de *moduli***.

Dizemos que esse funtor é **co-representado** por um esquema $M$ se existe uma transformação natural $\phi : \mathcal{F} \to h_M$ com propriedade universal. Se além disso temos que para todo corpo algebricamente fechado $k$ vale $\phi(k) : \mathcal{F}(\operatorname{Spec} k) \to \operatorname{Hom}(\operatorname{Spec} k, M)$ bijetor, dizemos que $M$ é um ***moduli* grosseiro** para $\mathcal{F}$. Nesse caso, apresentar uma família de objetos parametrizada por B determina um morfismo $B \to M$ chamado **morfismo classificante**.

Quando um funtor $\mathcal{F}$ é isomorfo a um funtor de pontos $h_M$, dizemos que esse funtor é **representável**, ou representado por M. Isso dá uma propriedade adicional ao espaço de moduli M: existe uma família universal $U \to M$, ou seja, dada uma família de objetos $T \to B$, o morfismo $\phi_T : B \to M$ realiza T como imagem inversa de U. Nesse caso, dizemos que M é um espaço de ***moduli* fino**. Observe que se M representa $\mathcal{F}$, ele também o co-representa, mas a recíproca não é verdadeira.

Porém, nem sempre um esquema traz a informação mais completa para os objetos que tentamos parametrizar. Por exemplo, quanto essas classes de objetos guardam automorfismos, ou seja, 'multiplicidades' ao nível de seus morfismos e não só nos objetos. Para esses casos, surge a necessidade de estudar os funtores de um problema de *moduli* não só como um feixe a partir dos objetos, mas também sobre morfismos. O conceito de *stack* nos leva a essa direção.

## A.3 *Stacks* como funtores

Um **grupoide** é uma categoria cujos morfismos são isomorfismos. A categoria dos grupoides, (Gpd), é na verdade uma 2-categoria, uma vez que os morfismos entre dois grupoides, que são funtores, formam por si uma categoria. De outra forma, a 2-categoria (Gpd) tem como objetos os grupoides, seus 1-morfismos são funtores entre grupoides e seus 2-morfismos são transformações naturais entre esses funtores.

Considere então um 2-funtor contravariante do tipo

$$\mathcal{F} : (Sch/\mathbb{C}) \to (Gpd),$$

onde tomamos $(Sch/\mathbb{C})$ como a 2-categoria associada à categoria dos esquemas. Dizemos que tal funtor é um **pré-feixe de grupoides**. Observe que para um esquema B, $\mathcal{F}(B)$ é um grupoide (e portanto uma categoria) e para um morfismo de esquemas $f : B' \to B$, $\mathcal{F}(f)$ é um funtor entre grupoides.

Tomando uma topologia de Grothendieck para $(Sch/\mathbb{C})$, dizemos que $\mathcal{F}$ é um **feixe de grupoides** se, para toda cobertura $\{f_i : U_i \to U\}$ valem os axiomas de feixe tanto para objetos quanto para morfismos:

- Se X e Y são objetos de $\mathcal{F}(U)$ e $\varphi_i : \mathcal{F}(f_i)(X) \to \mathcal{F}(f_i)(Y)$ são

morfismos tais que $\mathcal{F}(f_{ij,i})(\varphi_i) = \mathcal{F}(f_{ij,j})(\varphi_j)$, então existe morfismo $\eta : X \to Y$ tal que $\mathcal{F}(f_i)(\eta) = \varphi_i$.

- Se $X$ e $Y$ são objetos de $\mathcal{F}(U)$ e $\varphi : X \to Y$ e $\psi : X \to Y$ são morfismos quais que $\mathcal{F}(f_i)(\varphi) = \mathcal{F}(f_i)(\psi)$, então $\varphi = \psi$.
- Se $X_i$ é objeto de $\mathcal{F}(U_i)$ para cada $i$ e $\varphi_{ij} : \mathcal{F}(f_{ij,j})(X_j) \to \mathcal{F}(f_{ij,i})(X_i)$ são morfismos satisfazendo a condição de cociclo

$$\mathcal{F}(f_{ijk,ij})(\varphi_{ij}) \circ \mathcal{F}(f_{ijk,jk})(\varphi_{jk}) = \mathcal{F}(f_{ijk,ik})(\varphi_{ik})$$

  então existe um objeto $X$ em $\mathcal{F}(U)$ e isomorfismos $\varphi_i : \mathcal{F}(f_i)(X) \to X_i$ tais que $\varphi_{ji} \circ \mathcal{F}(f_{ij,i})(\varphi_i) = \mathcal{F}(f_{ij,j})(\varphi_j)$ onde os morfismos acima são como os da seção anterior e também $f_{ijk,ij} : (U_i \times_U U_j) \times_U U_k \to U_i \times_U U_j$. Mais uma vez, as imagens $\mathcal{F}(f_i)(\varphi)$ devem ser pensados como restrições $\varphi|_i$.

Um ***stack*** nada mais é que um feixe de grupoides.

Observe que para um esquema $M$, pode-se verificar que o funtor $h_M$ é um *stack*. De fato, basta considerar que para o esquema $B$, $h_M(B) = \mathrm{Hom}(B, M)$ é um grupoide cujos isomorfismos são apenas as identidades. Nesse caso, dizemos que $h_M$ é *stack* associado a $M$ e, por um abuso bem justificável de notação, o denotamos pelo mesmo nome, $M$. Vemos então que os *stacks* incluem os esquemas, mas podem incluir novos objetos:

**Lema A.3.1.** *Se um* stack *possui objetos com automorfismos além da identidade, ele não pode ser representado por um esquema.*

Um *stack* $\mathcal{F}$ é dito **representado** por um esquema $M$ se é isomorfo ao *stack* associado a $M$. Através desse isomorfismo, $\mathcal{F}$ herda as propriedades de $M$, o que também justifica usarmos a mesma terminologia usada em esquemas para *stacks*.

Quando $\mathcal{F}$ não é representável, ainda assim podemos definir uma estrutura algébrica sobre ele que permite considerar propriedades e construções análogas às de esquemas.

## A.4 *Stacks* como categorias

Antes de definir um *stack* algébrico, vamos descrever os *stacks* a partir de uma **categoria fibrada sobre grupoides** (ou, em sigla, CFG).

Trata-se de uma categoria $\mathcal{C}$ e um funtor contravariante de projeção $p : \mathcal{C} \to (\mathrm{Sch}/\mathbb{C})$ tal que, para um dado esquema $B$, a categoria $\mathcal{C}_B$ formada pelos objetos de $\mathcal{C}$ que se projetam por $p$ sobre $B$ (a fibra de $p$ sobre $B$) é um grupoide. Mais precisamente

- Dado morfismo $f : B' \to B$ e objeto $X$ de $\mathcal{C}_B$, existe ao menos um objeto $X'$ de $\mathcal{C}_{B'}$ tal que o diagrama abaixo comute

$$\begin{array}{ccc} X' & \overset{\varphi}{\dashrightarrow} & X \\ \downarrow & & \downarrow \\ B' & \xrightarrow{f} & B \end{array}$$

- Para todo diagrama comutativo

$$\begin{array}{ccc} B'' & \longrightarrow & B' \\ & \searrow \quad \swarrow & \\ & B & \end{array}$$

  e morfismos $X' \to X$ e $X'' \to X$ em diagramas comutativos sobre $B' \to B$ e $B'' \to B$ respectivamente, existe um único morfismo $X'' \to X'$ sobre $B'' \to B'$ que completa os diagramas comutativos.

Basicamente, a segunda condição implica que o objeto garantido na primeira é único a menos de isomorfismo, o qual denotamos por $f^*X$. Também, $\varphi$ é isomorfismo se e somente se $f$ é isomorfismo.

Se $\mathcal{C}$ é uma CFG, podemos construir o 2-funtor, e portanto um prefeixe de grupoides:

$$\begin{array}{rccc} \mathcal{F}_{\mathcal{C}} : & (\mathrm{Sch}/\mathbb{C}) & \to & (\mathrm{Gpd}) \\ & B & \mapsto & \mathcal{C}_B. \end{array}$$

Por outro lado, a partir de um prefeixe de grupoides $\mathcal{F}$, podemos construir a CFG $\mathcal{C}_{\mathcal{F}}$ com objetos $(B, X)$, onde $B$ é um esquema e $X$ um objeto de $\mathcal{F}(B)$, e morfismos $(f, \varphi) : (B', X') \to (B, X)$, onde $f : B' \to B$ é morfismo de esquemas e $\varphi : f^*X \to X'$ é um isomorfismo de grupoides.

Dessa forma, todas as definições anteriores para construir um *stack* a partir de um 2-funtor têm análogas em termos de uma CFG. Em particular, podemos denotar o objeto $\mathcal{F}(f)(X)$ por $f^*X$.

A partir dessa equivalência de definições, usaremos a forma mais conveniente nas próximas definições.

## A.5 *Stacks* algébricos

Algumas construções e propriedades conhecidas para esquemas podem ser estendidas para *stacks* com as devidas adequações da linguagem de categorias e funtores. Por exemplo, um **morfismo** de *stacks* $f : \mathcal{F} \to \mathcal{G}$ é um funtor que comuta com os funtores de projeção (em *stacks* como 2-categorias). Quando um morfismo de *stacks* é uma equivalência de categorias, é dito um **isomorfismo** de *stacks*. Observe também que dizer que temos um diagrama comutativo de morfismos de *stacks* como abaixo

$$\begin{array}{ccc} \mathcal{F} & \xrightarrow{h} & \mathcal{G} \\ & {\scriptstyle g}\searrow \quad \swarrow{\scriptstyle f} & \\ & \mathcal{H} & \end{array}$$

significa que $h \circ g$ é isomorfo a $f$ como funtores, mas não necessariamente iguais.

Os morfismos de $\mathcal{F}$ em $\mathcal{G}$ formam um grupoide $\mathrm{Hom}(\mathcal{F}, \mathcal{G})$ (cujos morfismos são transformações naturais). O lema de Yoneda tem uma consequência fundamental para a interpretação de um stack:

**Lema A.5.1.** *Dado um* stack $\mathcal{F}$ *e um esquema* $U$*, o funtor*

$$\begin{array}{rcl} \mathrm{Hom}(U, \mathcal{F}) & \to & \mathcal{F}(U) \\ f : (Sch/U) \to \mathcal{F} & \mapsto & f(id_U) \end{array}$$

*é uma equivalência de categorias.*

Ou seja, se $\mathcal{F}$ representa um problema de *moduli* podemos de fato considerá-lo como um *moduli* fino para o problema.

Por um **produto fibrado** entre morfismos de *stacks* $f : \mathcal{F} \to \mathcal{H}$ e $g : \mathcal{G} \to \mathcal{H}$, entendemos um outro *stack* denotado por $\mathcal{F} \times_{\mathcal{H}} \mathcal{G}$ com propriedade universal tal que o diagrama

$$\begin{array}{ccc} \mathcal{F} \times_{\mathcal{H}} \mathcal{G} & \longrightarrow & \mathcal{F} \\ \downarrow & & \downarrow \\ \mathcal{G} & \longrightarrow & \mathcal{H} \end{array}$$

seja compatível com os funtores de projeção.

As imagens de um objeto do produto fibrado em $\mathcal{H}$ pelas duas composições de morfismos não precisam ser necessariamente iguais, mas relacionadas por morfismo em $\mathcal{H}$. Desse modo, fixado um esquema $U$ um objeto do produto fibrado é um trio $(X, Y; \alpha)$ no qual $X$ é objeto de $\mathcal{F}_U$, $Y$ de $\mathcal{G}_U$ e $\alpha$ é morfismo em $\mathcal{H}$ sobre $\mathrm{id}_U$, com $\alpha : f(X) \to g(Y)$. Analogamente, as imagem de um morfismo do produto fibrado em $\mathcal{H}$ pelas duas composições não precisam ser necessariamente um mesmo morfismo, mas satisfazer um diagrama comutativo. Explicitamente, dado um morfismo de esquemas $U \to V$, um morfismo do produto fibrado é um par $(\varphi, \psi)$, que relaciona dois objetos $(X, Y; \alpha)$ sobre $U$ e $(X', Y'; \beta)$ sobre $V$ de modo que $\varphi : X \to X'$ e $\psi : Y \to Y'$ sobre $U \to V$ com $\beta \circ f(\varphi) = g(\psi) \circ \alpha$.

A universalidade do produto fibrado significa então que, se há outro *stack* $\mathcal{P}$ satisfazendo essas condições, há um morfismo de stacks $\mathcal{P} \to \mathcal{F} \times_{\mathcal{H}} \mathcal{G}$ único a menos de um isomorfismo (de funtores). O diagrama é cartesiano se $\mathcal{P} \to \mathcal{F} \times_{\mathcal{H}} \mathcal{G}$ já é um isomorfismo de *stacks*.

Agora dizemos que um morfismo de *stacks* $\mathcal{F} \to \mathcal{G}$ é representável se para todo esquema $U$ e morfismo $U \to \mathcal{G}$ temos o produto fibrado $U \times_{\mathcal{G}} \mathcal{F}$ um *stack* representável. Uma propriedade local e estável por mudança de base para morfismos de esquemas (como separado, quasi-compacta, liso, plano, *étale*, tipo finito...) é transmitida para o morfismo de *stacks* se o morfismo de esquemas $U \times_{\mathcal{G}} \mathcal{F} \to U$ também tem essa propriedade para todo $U$.

Podemos definir o morfismo diagonal $\Delta_{\mathcal{F}} : \mathcal{F} \to \mathcal{F} \times \mathcal{F}$. Esse morfismo é a peça chave para o entendimento do *stack*, uma vez que ele guarda informação sobre os automorfismos de um objeto. De fato, dado um esquema $U$, um morfismo $U \to \mathcal{F} \times \mathcal{F}$ determina, pelo lema de Yoneda, dois objetos $X$ e $Y$ de $\mathcal{F}(U)$. O produto fibrado entre $U$ e $\mathcal{F}$ sobre $\mathcal{F} \times \mathcal{F}$ resulta no grupoide de isomorfismos entre $X$ e $Y$ sobre $U$, $\mathrm{Iso}_U(X, Y)$.

$$\begin{array}{ccc} \mathrm{Iso}_U(X,Y) & \longrightarrow & \mathcal{F} \\ \downarrow & & \downarrow{\scriptstyle \Delta_{\mathcal{F}}} \\ U & \longrightarrow & \mathcal{F} \times \mathcal{F} \end{array}$$

Um ***stack* algébrico** $\mathcal{M}$ (no sentido de Deligne-Mumford) é um *stack* cuja diagonal $\Delta_{\mathcal{M}}$ é representável, quasi-compacta e separável e existe um esquema $U$ equipado com um morfismo sobrejetor *étale* $U \to \mathcal{M}$. O esquema $U$ é dito um **atlas** para $\mathcal{M}$ [1].

[1]Denotaremos por $\mathcal{M}$ ao invés de $\mathcal{F}$ quando o *stack* for algébrico

Para um *stack* algébrico $\mathcal{M}$ é possível mostrar que sua diagonal $\Delta_{\mathcal{M}}$ também é não ramificada. Isso implica que, para X um objeto de $\mathcal{M}_U$, o feixe de grupos de automorfismos $\mathrm{Aut}_U(X) = \mathrm{Iso}_U(X, X)$ é finito em $\mathcal{M}_U$.

Os **pontos** de um *stack* algébrico podem ser tomados como aqueles provenientes de seu atlas. Mais precisamente, se k é um corpo e $x : \mathrm{Spec}\, k \to \mathcal{M}$ morfismo, um ponto do *stack* é definido por uma relação de equivalência entre pares $(k, x)$: $(k, x)$ e $(k', x')$ são equivalentes que existe uma extensão de corpos K de k e $k'$ tal que o diagrama

$$\begin{array}{ccc} \mathrm{Spec}\, K & \longrightarrow & \mathrm{Spec}\, k' \\ \downarrow & & \downarrow \\ \mathrm{Spec}\, k & \longrightarrow & \mathcal{M} \end{array}$$

seja comutativo.

Outras características do *stack* algébrico também são definidas por meio de seu atlas U: propriedades locais válidas para o esquema U são transferidas para o *stack* algébrico $\mathcal{M}$. Propriedades de um morfismo entre *stacks* algébricos $f : \mathcal{M} \to \mathcal{N}$ também são provenientes de morfismos de esquemas lidos de diagramas comutativos

$$\begin{array}{ccc} U & \longrightarrow & \mathcal{M} \\ \downarrow{\scriptstyle g} & & \downarrow{\scriptstyle f} \\ V & \longrightarrow & \mathcal{N} \end{array}$$

onde U e V são atlas de $\mathcal{M}$ e $\mathcal{N}$ e g é morfismo de esquemas com a propriedade analisada.

Um exemplo importante de *stack* algébrico: Se X é um esquema sob ação de um grupo G, o *stack* quociente $[X/G]$ ([12], exemplo 4.8). Como categoria, seus objetos são G-fibrados principais sobre X e seus morfismos são morfismos de G-fibrados compatíveis com a base X.

Um *stack* $\mathcal{L}$ é um ***substack*** algébrico de $\mathcal{M}$ se há um morfismo de *stacks* de $\mathcal{L} \to \mathcal{M}$ representado por mergulhos de esquemas fechados. As definições de *substack* aberto, fechado e localmente fechado, seguem da representabilidade do morfismo de inclusão. Seguem também os conceitos de conexidade e irredutibilidade. Para *stacks* separáveis e próprios, critérios de verificação são apresentados em [18].

Um **feixe** quase-coerente $\mathcal{E}$ sobre um *stack* algébrico $\mathcal{M}$ é uma coleção de feixes quase-coerentes $\mathcal{E}_X$ sobre todo esquema X equipado com isomorfismos $\varphi_f : \mathcal{E}_X \to \mathcal{M}(f)(\mathcal{E}_Y)$ satisfazendo as condições de cociclo sobre cada diagrama comutativo

$$\begin{array}{ccc} X & \xrightarrow{\quad f \quad} & Y \\ & \searrow \qquad \swarrow & \\ & \mathcal{M} & \end{array}$$

Podemos denotar $f^*\mathcal{E}_Y := \mathcal{M}(f)(\mathcal{E}_Y)$.

Ainda, dizemos que $\mathcal{E}$ é feixe coerente ou localmente livre sobre o *stack* algébrico se cada $\mathcal{E}_X$ é coerente ou localmente livre, respectivamente, sobre esquemas. Também, um morfismo entre feixes sobre *stacks* $h : \mathcal{E}' \to \mathcal{E}$ é uma coleção de morfismos de feixes sobre esquemas $h_X : \mathcal{E}'_X \to \mathcal{E}_X$ compatíveis com os isomorfismos $\varphi_f$.

De outra forma, um fibrado vetorial de posto r sobre o *stack* algébrico $\mathcal{M}$ é um morfismo de *stacks* $\mathcal{M} \to \mathcal{V}_r$, onde $\mathcal{V}_r$ é o *stack* (de fato um grupoide) dos fibrados vetoriais de posto r. O fato de ser um morfismo de *stacks* implica que se devem satisfazer condições de cociclos e respeitar mudança de bases. [17]

Os conceitos de dimensão e grau também podem ser estendidos. A **dimensão** de um *stack* algébrico ([18] 2.6) também pode ser definida a partir de seu atlas. Se $u : U \to \mathcal{M}$ é o morfismo do atlas e $x$ é um ponto de $\mathcal{M}$,

$$\dim_x \mathcal{M} = \dim_x U - \dim_x u,$$

onde $\dim u$ é a dimensão relativa de $u$ sobre $x$.

O **grau** de um morfismo de *stacks* algébricos ([29] 1.15) $\mathcal{M} \to \mathcal{N}$ é em geral um número racional dado pela fórmula

$$\deg(\mathcal{M}/\mathcal{N}) = \frac{a(\mathcal{M})}{a(\mathcal{N})}[k(\mathcal{M}) : k(\mathcal{N})],$$

onde $a(\mathcal{M})$ é a ordem do grupo de automorfismos sobre o ponto genérico de $\mathcal{M}$ e $k(\mathcal{M})$ seu corpo de funções racionais.

Temos $a(\mathcal{M}) = 1$ se e somente se $\mathcal{M}$ é representado por um esquema (ou tem um *substack* aberto isomorfo a um esquema). Isso justifica o fato de tomarmos o grupo de ciclos de interseção com coeficientes racionais, como será definido à frente.

## A.6 Teoria de interseção em *stacks*

A descrição de um *stack* algébrico por seu atlas nem sempre é conveniente. Uma maneira mais intuitiva de associar um *stack* algébrico a um esquema é através de seu espaço de *moduli*. Dizemos que um esquema $M$ é um **espaço de *moduli*** para o *stack* algébrico $\mathcal{M}$ se existe morfismo próprio $p : \mathcal{M} \to M$ que estabelece bijeção entre componentes conexas dos grupoides $\mathcal{M}(\operatorname{Spec} k)$ e $M(\operatorname{Spec} k)$, para $k$ qualquer corpo algebricamente fechado. O esquema $M$ não é necessariamente um espaço de *moduli* grosseiro como definido anteriormente, mas é possível verificar que, se $\mathcal{M}$ tem um espaço de *moduli*, então existe esquema $M'$ que é um espaço de *moduli* grosseiro para $\mathcal{M}$, com $M$ e $M'$ relacionados por um morfismo que os torna equivalentes perante a teoria de interseção com coeficientes racionais (veja [29], definição 2.1). Essa construção é particularmente importante para a construção de uma teoria de interseção para *stacks* algébricos.

A partir de então, supomos que $\mathcal{M}$ é um *stack* algébrico (de tipo finito sobre $k$) com um espaço de moduli $M$. Nessas hipóteses, seguem as construções da teoria de interseção de modo análogo para esquemas (como em [1], seção 2): os grupos de ciclos $Z_k(\mathcal{M})_{\mathbb{Q}}$ gerado por substacks integralmente fechados de dimensão $k$, a equivalência racional e o **grupo de Chow** $A_*(\mathcal{M})_{\mathbb{Q}}$. Definimos ainda de modo análogo a imagem direta própria e imagem inversa plana no grupo de Chow, que apresentam as mesmas propriedades usuais.

Para o morfismo próprio $p : \mathcal{M} \to M$, se $\mathcal{V}$ é um *substack* integralmente fechado em $\mathcal{M}$ e $V$ seu espaço de *moduli* temos

$$p_*[\mathcal{V}] = \frac{1}{a(\mathcal{V})}[V],$$

onde $a(\mathcal{V})$ é a ordem do grupo de automorfismos do ponto genérico de $\mathcal{V}$.

Essa relação estabelece um isomorfismo do grupo de Chow do *stack* algébrico com o de seu espaço de *moduli*:

**Proposição A.6.1.** *Seja $\mathcal{M}$ um* stack *algébrico com $p : \mathcal{M} \to M$ seu espaço de* moduli*. Então a imagem direta própria $p_* : A_*(\mathcal{M})_{\mathbb{Q}} \to A_*(M)_{\mathbb{Q}}$ é um isomorfismo (de grupos).*

Se $\mathcal{M}$ é próprio, denotamos a imagem direta sobre um ponto como

$$\int_{\mathcal{M}} : A_*(\mathcal{M}) \to \mathbb{Q}.$$

A construção de uma teoria de interseção bivariante segue a mesma estrutura do caso de esquemas - desde que considerada sobre coeficientes racionais. Temos assim um **anel de Chow bivariante** $A^*(\mathcal{M})_{\mathbb{Q}}$, comutativo com o usual produto $\cup$. Vale também:

**Proposição A.6.2.** *Seja $\mathcal{M}$ um* stack *algébrico com $p : \mathcal{M} \to M$ seu espaço de* moduli. *Então a imagem inversa própria $p^* : A^*(M)_{\mathbb{Q}} \to A^*(\mathcal{M})_{\mathbb{Q}}$ é um isomorfismo (de anéis).*

Definimos também de forma usual o produto $\cap$: dados $\alpha \in A^i(\mathcal{M})_{\mathbb{Q}}$ e $\xi \in A_j(\mathcal{M})_{\mathbb{Q}}$, tomamos $\alpha \cap \xi \in A_{j-i}(\mathcal{M})_{\mathbb{Q}}$. Segue a fórmula de projeção. Em particular, $A_*(\mathcal{M})_{\mathbb{Q}}$ é módulo sobre $A^*(\mathcal{M})_{\mathbb{Q}}$

Se o *stack* algébrico $\mathcal{M}$ é liso, ainda vale a dualidade $A^*(\mathcal{M})_{\mathbb{Q}} \to A_*(\mathcal{M})_{\mathbb{Q}}$ determinada pelo isomorfismo $\alpha \mapsto \alpha \cap [\mathcal{M}]$. Isso ainda estabelece um isomorfismo $A^*(M)_{\mathbb{Q}} \to A_*(M)_{\mathbb{Q}}$ mesmo quando $M$ não for liso, pois pela fórmula de projeção

$$p_*(p^*\alpha \cap [\mathcal{M}]) = \alpha \cap p_*[\mathcal{M}] = \frac{1}{a(\mathcal{M})}\alpha \cap [M].$$

Também, para $\mathcal{E}$ um fibrado vetorial sobre $\mathcal{M}$, define-se as classes de Chern $c_k(\mathcal{E}) \in A^k(\mathcal{M})_{\mathbb{Q}}$ seguindo os axiomas/propriedades usuais, como mostrado em [1], seção 2. Em particular, está definida a classe de Euler $\mathrm{Euler}(\mathcal{E})$.

Em suma, se $M$ é um espaço de *moduli* para um *stack* algébrico liso $\mathcal{M}$, então $M$ se comporta como um esquema liso do ponto de vista da teoria da interseção sobre coeficientes racionais, pois todas propriedades de seu anel de Chow podem ser lidas do anel de Chow do *stack* liso. Em particular, os esquemas e *stacks* algébricos em que estamos interessados apresentam essa característica (em particular, *stacks* algébricos com singularidades quocientes, como definido em [29] 2.8).

### A.6.1 Classes virtuais e Localização

Mesmo para esquemas, definir classes de zeros de uma seção de um fibrado vetorial temos a situação onde não podemos garantir que a dimensão seja constante ao longo das fibras, de onde surge a noção de dimensão esperada, um limite inferior para a dimensão efetiva (como discutido em [14]). Trabalhamos então com a classe fundamental virtual desses subesquemas, considerando apenas a dimensão esperada no anel de Chow. Essa noção também é estendida sobre *stacks* algébricos. A construção de tais classes é feita sobre certas hipóteses em [6] ou [27], usando uma *teoria de obstrução perfeita*. Essa construção também é discutida e exemplificada em [11], 7. O caso específico do *stack* de mapas estáveis também é tratado em [24].

Seja $\mathcal{L}$ um *substack* mergulhado em um *stack* liso $\mathcal{M}$ tal que $\mathcal{L}$ apresente uma teoria de obstrução perfeita. Primeiramente, define-se a classe fundamental virtual $[\mathcal{L}]^{virt}$ em $A_d(\mathcal{L})_{\mathbb{Q}}$, onde d é a dimensão esperada. ([19], parágrafo 1).

Esse é o caso do lugar $\mathcal{Z}$ dos zeros de uma seção s de um fibrado vetorial $\mathcal{E}$ sobre $\mathcal{M}$ ([11], capítulo 7). Em particular, dada a inclusão $i : \mathcal{Z} \to \mathcal{M}$, vale $i_*([\mathcal{Z}]^{virt}) = c_{top}(\mathcal{E}) \cap [\mathcal{M}]^{virt} \in A_d(\mathcal{M})_{\mathbb{Q}}$, onde d é a dimensão esperada.

Se sob o *stack* liso $\mathcal{M}$ temos a ação de um toro, podemos também construir seu anel de Chow equivariante e estender a fórmula de localização de Bott. Apresentaremos aqui tais conceitos de forma geral. Uma descrição mais rigorosa está em [19], no contexto de teoria de obstrução perfeita sob ação de um toro. Outra abordagem está em [11], capítulo 9, no contexto de *orbifolds*.

Dada ação de um toro $T = (\mathbb{C}^*)^k$ sobre o *stack* liso $\mathcal{M}$, define-se o **anel de Chow equivariante** em coeficientes racionais para $\mathcal{L}$ por

$$A_*^T(\mathcal{L})_{\mathbb{Q}} = A_*([(\mathcal{L} \times ET)/T])_{\mathbb{Q}},$$

definido a partir *stack* quociente, onde $ET \to BT$ é o T-fibrado principal da ação. Para $\mathcal{E}$ fibrado linear T-equivariante, define-se as classes de Chern T-equivariantes $c_k^T(\mathcal{E})$. Em particular, também está definida a classe de Euler T-equivariante $\mathrm{Euler}^T(\mathcal{E})$.

Se denotamos por $\mathcal{R}_T$ o corpo de frações do anel T-equivariante de um ponto, temos que a localização do anel de Chow T-equivariante é $A_*^T(\mathcal{L})_{\mathbb{Q}} \otimes \mathcal{R}_T$.

Considere agora $\iota_i : \mathcal{L}_i \to \mathcal{L}$ os mapas de inclusão das componentes do lugar de pontos fixos em $\mathcal{L}$ pela ação de T. Considere também $\mathcal{N}_i^{virt}$

o fibrado normal virtual em $\mathcal{M}$ para cada inclusão ([19], parágrafo 1). Verifica-se que $\mathcal{N}_i^{virt}$ é fibrado T-equivariante e sua classe de Euler $\mathrm{Euler}^T(\mathcal{N}_i^{virt})$ é invertível na localização $A_*^T(\mathcal{L}_i)_{\mathbb{Q}} \otimes \mathcal{R}_T$. Com essas definições, para $\alpha \in A_*^T(\mathcal{L})_{\mathbb{Q}} \otimes \mathcal{R}_T$, temos

$$\alpha = \sum_i \iota_{i*} \frac{\iota_i^* \alpha}{\mathrm{Euler}^T(\mathcal{N}_i^{virt})}.$$

A fórmula de localização, onde $\int$ denota a imagem direta sobre a localização do anel de Chow T-equivariante de um ponto, é então

$$\int_{[\mathcal{L}]^{virt}} \alpha = \sum_i \int_{[\mathcal{L}_i]^{virt}} \iota_{i*} \frac{\iota_i^* \alpha}{\mathrm{Euler}^T(\mathcal{N}_i^{virt})}.$$

Em particular, destacamos que no nosso caso de interesse, o *stack* de mapas estáveis $\overline{\mathcal{M}}_{0,n}(\mathbb{P}^r, d)$ é tal que $[\overline{\mathcal{M}}_{0,n}(\mathbb{P}^r, d)]^{virt} = [\overline{\mathcal{M}}_{0,n}(\mathbb{P}^r, d)]$ ([11], 7.1.5 e 9.3).

# Apêndice B

# Sequências exatas sobre curvas nodais

Explicitamos aqui a base para a construção de sequências exatas sobre curvas nodais, que desempenham papel fundamental para as demonstrações e cálculos ao longo do texto principal.

## B.1 Sequência de excisão

De acordo com [21], II.1, exercício 1.19, seja $X$ um espaço topológico, $Z$ um fechado com inclusão $j : Z \to X$ e $U = X \setminus Z$ com sua inclusão $i : U \to X$. Dado um feixe $\mathcal{F}$, construímos os feixes imagem inversa $\mathcal{F}|_Z := j^{-1}\mathcal{F}$ em $Z$ e $\mathcal{F}|_U := i^{-1}\mathcal{F}$ em $U$, cujos *stalks* são $\mathcal{F}_P$, sobre um ponto $P \in X$ ou $P \in U$ respectivamente.

Por sua vez, tomamos o feixe imagem direta $j_*(\mathcal{F}|_Z)$, cujo *stalk* sobre $P \in X$ é o próprio $\mathcal{F}_P$ se $P \in Z$ e $0$ caso contrário. Por simplicidade, denotamos esse feixe por $\mathcal{F}_Z$.

Analogamente, definimos o feixe $\mathcal{F}_U := i_!(\mathcal{F}|_U)$, o feixe de extensão por zero, de modo que que *stalk* em $P \in X$ seja $\mathcal{F}$ se $P \in U$ e $0$ caso contrário. A partir desses feixes e da propriedade de adjunção entre $j^{-1}$ e $j_*$, construímos a sequência exata

$$0 \longrightarrow \mathcal{F}_U \longrightarrow \mathcal{F} \longrightarrow \mathcal{F}_Z \longrightarrow 0.$$

Além disso, temos que para qualquer feixe $\mathcal{G}$ sobre o fechado $Z$ vale $H^i(X, j_*\mathcal{G}) = H^i(Z, \mathcal{G})$ para todo $i \geq 0$. [21], II.2, lema 2.10.

Nosso principal interesse é usar essas sequência para curvas nodais que são árvores de retas projetivas. Para isso, considere $C$ uma tal curva e sua normalização $\nu : X \sqcup Y \to C$, onde $X$ e $Y$ são componentes de $C$, possivelmente nodais, tais que, temos $P \in X$ e $Q \in Y$ correspondendo a um nó $N \in C$. Sendo assim, tome $\nu|_X$ como mapa de inclusão de $X$ e $\nu|_{Y\setminus Q}$ o mapa de inclusão do aberto complementar $C \setminus X$. Dessa forma, temos a identificação

$$\nu|_{Y\setminus Q_*}(\mathcal{F}|_{C\setminus X}) = \mathcal{F}_Y \otimes \mathcal{O}_C(-N).$$

Resumindo, para $C$ uma árvore de retas projetivas, vale

$$0 \longrightarrow \mathcal{F}_Y \otimes \mathcal{O}_C(-N) \longrightarrow \mathcal{F} \longrightarrow \mathcal{F}_X \longrightarrow 0.$$

Tomando cohomologias, temos a sequência exata longa

$$\cdots \longrightarrow H^i(Y, \mathcal{F}|_Y \otimes \mathcal{O}_Y(-Q)) \longrightarrow H^i(C, \mathcal{F}) \longrightarrow H^i(X, \mathcal{F}|_X) \longrightarrow \cdots$$

# Apêndice C

# Grafos

Ao usar a fórmula de localização para mapas estáveis, a descrição combinatória das componentes de pontos fixos se faz através de problemas clássicos de grafos. Nesta seção apresentamos, de forma sucinta, os conceitos e terminologias da teoria de grafos utilizados ao longo do texto principal. Como referência básica indicamos [8].

Um **grafo** $\Gamma$ é um par de conjuntos $\{V, E\}$ com a seguinte descrição: um elemento $v \in V$ é chamado **vértice**. Um elemento $e \in E$ é chamado **aresta**, estabelecendo uma relação entre 2 vértices $v, w \in V$ e é escrito da forma $e = (vw)$. Dois vértices relacionados por uma aresta são ditos **adjacentes**. O número de arestas adjacentes a um vértice $v$ (com laços contados 2 vezes) será chamado **valência** de $v$ e denotado por $val(v)$.

A representação topológica usual é a de vértices como pontos e as arestas segmentos ligando esses pontos (ou um CW-complexo de dimensão 1). Dessa representação topológica seguem os conceitos de conexidade, laços e ciclos. Um grafo sem ciclos é chamado de **árvore**, devido a sua forma 'ramificada'.

Dizemos que um grafo é **planar** se pode ser mergulhado no plano, ou seja, representado num plano sem que suas arestas se cruzem. Em particular, toda árvore é um grafo planar.

Sobre as arestas de um grafo podemos associar uma **função peso** $\delta : E \to S$. Um grafo $\Gamma = \{V, E, \delta\}$ é dito um **grafo ponderado**. A soma dos pesos das arestas adjacentes a um vértice $v$ será chamada **multiplicidade** de $v$ e denotada por $\mu(v)$.

Podemos analogamente associar também funções $\iota : V \to S$ aos vértices de um grafo, que chamaremos aqui de **rótulos**. Chamaremos também de **grafo rotulado** um grafo $\Gamma = \{V, \iota, E\}$. Um problema importante envolvendo esses rótulos é o da coloração de um grafo.

## C.1 Coloração de grafos

O problema da **coloração** de um grafo $\Gamma$ consiste em distribuir os rótulos de um conjunto finito $K$ entre seus vértices de modo que dois vértices adjacentes não tenham o mesmo rótulo. O nome provém do problema em colorir mapas com cores diferentes de forma que duas regiões vizinhas não tenham a mesma cor. Por essa analogia, nos referimos aos rótulos em $K$ satisfazendo essas condições por **cores**.

O problema de coloração de grafos envolve, computacionalmente, o problema de decisão (dado $\Gamma$ e $K$, se $\Gamma$ pode ser colorido pelas cores em $K$, NP-completo), o problema de otimização (a quantidade mínima de cores suficientes para colorir $\Gamma$, NP-difícil) e o de contagem (quantas colorações possíveis, $\sharp$P-completo). Em nosso estudo, estamos mais interessados no problema de contagem.

O número de colorações possíveis de um grafo é descrito por uma expressão polinomial. O **polinômio cromático** $p_\Gamma(t)$ é uma função real que interpola o número de colorações para $k \in \mathbb{N}$ cores: para $t = k$, o valor de $p_\Gamma(k)$ é o número de colorações possíveis de $\Gamma$ com até $k$ cores. Em particular, o número de colorações com exatamente $k$ cores é $p_\Gamma(k) - p_\Gamma(k-1)$.

Na instância em que estamos interessados, é suficiente saber que

1. Teorema das 4 cores: Todo grafo planar sem laços (em particular árvores) pode ser colorido com 4 cores.

2. Toda árvore pode ser colorida com apenas 2 cores.

3. O polinômio cromático de uma árvore com $\nu$ vértices é $p_\Gamma(t) = t(t-1)^{\nu-1}$.

## C.2 Isomorfismos

Um **isomorfismo** entre grafos $\Gamma_1 = \{V_1, E_1\}$ e $\Gamma_2 = \{V_2, E_2\}$ é uma bijeção $f : V_1 \to V_2$ que preserva as incidências: $(\nu_1 w_1) \in E_1$ se e somente se $(f(\nu_1) f(w_1)) \in E_2$. Isso estabelece uma relação de equivalência entre grafos e define classes de isomorfismos de grafos.

Se os grafos são ponderados, ainda podemos querer que os pesos também sejam preservados pelo isomorfismo, ou seja, para cada aresta $e = (v_1w_1) \in E_1$, vale $\delta_2(f(v_1)f(w_1)) = \delta_1(v_1w_1)$. Condição análoga é feita para grafos rotulados ou coloridos: um isomorfismos preserva rótulos se $\iota_1(v_1) = \iota_2(f(v_1))$ para todo $v_1 \in V_1$. Essas condições também definem classes de isomorfismos de pesos e rótulos.

Um **automorfismo** de um grafo $\Gamma$ é uma permutação $\sigma$ em $V$ que preserva as incidências: $(vw) \in E$ se e somente se $(\sigma(v)\sigma(w)) \in E$. Denotamos o conjunto de automorfismos de $\Gamma$ por $Aut(\Gamma)$ e por construção temos que é um subgrupo do grupo de permutações de vértices $Sym(V)$.

Também temos o conceito de automorfismos de grafos ponderados e rotulados. Se o grafo $\Gamma$ é ponderado, ainda pedimos que os pesos também sejam preservados pelas bijeções ou permutações, ou seja, para cada aresta $e = (vw)$, vale $\delta(\sigma(v_1)\sigma(v_2)) = \delta(v_1v_2)$. Para grafos rotulados ou coloridos: um automorfismos preserva rótulos se $\iota(v) = \iota(\sigma(v))$ para todo $v \in V$.

## C.3 Coloração e automorfismos

O polinômio cromático conta as colorações de grafos a despeito das classe de isomorfismos de coloração. Por exemplo, para um grafo simples de 2 vértices e 1 aresta $\Gamma = \{\{v, w\}, \{(vw)\}\}$, o polinômio cromático é $p_\Gamma(t) = t(t-1)$. Para duas cores $0, 1$ temos então 2 grafos coloridos possíveis $\Gamma_{01}$ e $\Gamma_{10}$, a saber para $\iota(v) = 0, \iota(w) = 1$ e para $\iota(v) = 1, \iota(w) = 0$ respectivamente. Porém, esses dois grafos coloridos são isomorfos a partir de um automorfismo do grafo original $\Gamma$.

Em nosso problema, nos deparamos com uma situação semelhante para contagem. Considere o conjunto de todas colorações com $k$ cores para um grafo fixo $\Gamma$, denotado por $Crom_\Gamma(k)$. Sabemos que sua cardinalidade é $|Crom_\Gamma(k)| = p_\Gamma(k) - p_\Gamma(k-1)$. Fixada também um distribuição de pesos $\delta$ sobre $\Gamma$, consideramos então o grupo de automorfismos de $\Gamma$ que preservem os pesos, $Aut_\delta(\Gamma) \leq Aut(\Gamma)$. A ação de $Aut_\delta(\Gamma)$ sobre $Crom_\Gamma(k)$ determina então classes de isomorfismos de grafos (ponderados) coloridos. Em particular, dado um grafo ponderado colorido $\Gamma_\iota \in Crom_\Gamma(k)$, seu estabilizador é o conjunto de seus automorfismos como grafo ponderado colorido.

Na lista de grafos no apêndice D, usamos esses dados para confirmar nossas contagens, como explicado em sua introdução.

# Apêndice D

# Tabelas de grafos

Nas tabelas, estão apresentados os tipos combinatórios de grafos parametrizando as componentes de mapas fixos $\mathcal{M}_\Gamma$, como construídas em (1.4.1), onde $\Gamma$ é um grafo ponderado colorido. Cada bloco corresponde a uma distribuição de pesos fixada e cada célula corresponde a um tipo combinatório de coloração $\iota$. A primeira linha de cada célula é uma representação desse tipo combinatório: as arestas com seus pesos (graus) e os vértices com sua coloração (imagem). Em seguida, apresentamos o número de automorfismos de $\Gamma$ como grafo ponderado colorido associado a um mapa estável, $a(\Gamma)$, que é o produto da cardinalidade do estabilizador de $\Gamma$ sob a ação descrita no apêndice C com o produto dos pesos (graus) de cada aresta. Indicamos também o número $c$ de classes de colorações possíveis, $c\times$.

Abaixo, o modelo das informações em cada célula:

| |
|---|
| $\iota(v_1) \overset{d_{e_1}}{\text{———}} \cdots \Gamma \cdots \overset{d_{e_m}}{\text{———}} \iota(v_n)$ |
| $a(\Gamma) = \lvert \mathrm{stb}\, \Gamma \rvert \cdot \Pi_e d_e$ |
| $c\times$ |

Como 'tipo combinatório' entendemos não somente a estrutura de vértices e arestas ponderadas do grafo, mas também a disposição de sua coloração (número de cores, cores repetidas e posição relativa). Para isso usamos as letras $i, j, k, l$ variando entre as quatro cores possíveis $0, 1, 2, 3$, correspondentes aos pontos fixos em $\mathbb{P}^3$, seguindo a notação em [11].

Caso $d = 2$

| |
|---|
| $i \overset{2}{\text{———}} j$ |
| $a(\Gamma) = 1 \cdot 2$ |
| $6\times$ |

| | |
|---|---|
| $i \overset{1}{\text{——}} j \overset{1}{\text{——}} k$<br>$a(\Gamma) = 1$<br>$12\times$ | $i \overset{1}{\text{——}} j \overset{1}{\text{——}} i$<br>$a(\Gamma) = 2$<br>$12\times$ |

Caso $d = 3$

| |
|---|
| $i \overset{3}{\text{——}} j$<br>$a(\Gamma) = 1 \cdot 3$<br>$6\times$ |

| | |
|---|---|
| $i \overset{2}{\text{——}} j \overset{1}{\text{——}} k$<br>$a(\Gamma) = 1 \cdot 2$<br>$24\times$ | $i \overset{2}{\text{——}} j \overset{1}{\text{——}} i$<br>$a(\Gamma) = 1 \cdot 2$<br>$12\times$ |
| $i \overset{1}{\text{——}} j \overset{1}{\text{——}} i \overset{1}{\text{——}} j$<br>$a(\Gamma) = 1$<br>$6\times$ | $i \overset{1}{\text{——}} j \overset{1}{\text{——}} i \overset{1}{\text{——}} k$<br>$a(\Gamma) = 1$<br>$24\times$ |
| $k \overset{1}{\text{——}} i \overset{1}{\text{——}} j \overset{1}{\text{——}} k$<br>$a(\Gamma) = 1$<br>$12\times$ | $i \overset{1}{\text{——}} j \overset{1}{\text{——}} k \overset{1}{\text{——}} l$<br>$a(\Gamma) = 1$<br>$24\times$ |

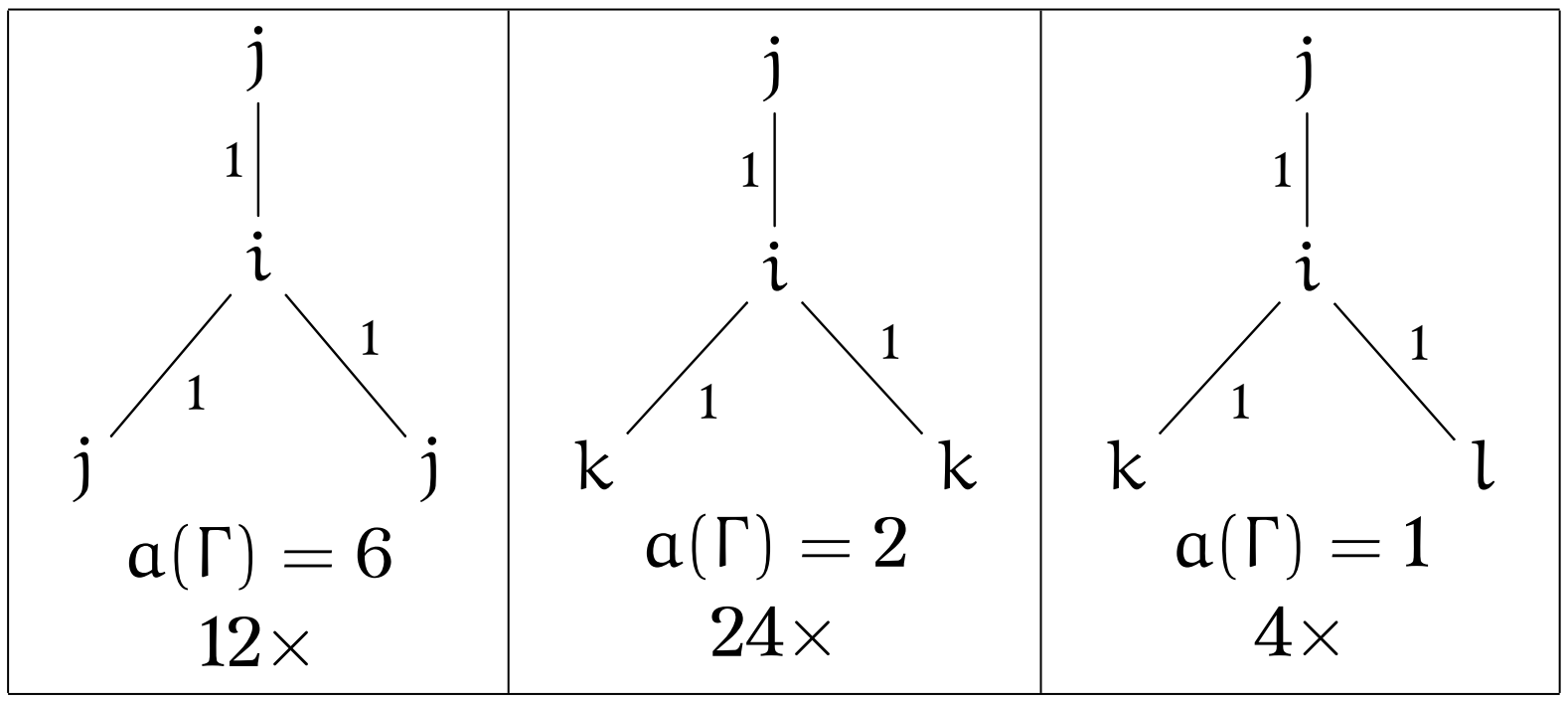

Caso $d = 4$

| |
|---|
| $i \overset{4}{\text{——}} j$<br>$\lvert A \rvert = 1 \cdot 4$<br>$6\times$ |

| | |
|---|---|
| $i \overset{3}{\text{——}} j \overset{1}{\text{——}} k$<br>$\lvert A \rvert = 1 \cdot 3$<br>$24\times$ | $i \overset{3}{\text{——}} j \overset{1}{\text{——}} i$<br>$\lvert A \rvert = 1 \cdot 3$<br>$12\times$ |
| $i \overset{2}{\text{——}} j \overset{2}{\text{——}} k$<br>$\lvert A\rvert = 1 \cdot 4$<br>$12\times$ | $i \overset{2}{\text{——}} j \overset{2}{\text{——}} i$<br>$\lvert A\rvert = 2 \cdot 4$<br>$12\times$ |

| | | |
|---|---|---|
| $i \overset{2}{\text{——}} j \overset{1}{\text{——}} k \overset{1}{\text{——}} l$<br>$\lvert A\rvert = 1 \cdot 2$<br>$24\times$ | $i \overset{2}{\text{——}} j \overset{1}{\text{——}} i \overset{1}{\text{——}} k$<br>$\lvert A\rvert = 1 \cdot 2$<br>$24\times$ | $i \overset{2}{\text{——}} j \overset{1}{\text{——}} k \overset{1}{\text{——}} j$<br>$\lvert A\rvert = 1 \cdot 2$<br>$24\times$ |
| $i \overset{2}{\text{——}} j \overset{1}{\text{——}} k \overset{1}{\text{——}} i$<br>$\lvert A\rvert = 1 \cdot 2$<br>$24\times$ | $i \overset{2}{\text{——}} j \overset{1}{\text{——}} i \overset{1}{\text{——}} j$<br>$\lvert A\rvert = 1 \cdot 2$<br>$12\times$ | |

| | |
|---|---|
| $i \overset{1}{\text{——}} j \overset{2}{\text{——}} k \overset{1}{\text{——}} l$<br>$\lvert A\rvert = 1 \cdot 2$<br>$12\times$ | $i \overset{1}{\text{——}} j \overset{2}{\text{——}} i \overset{1}{\text{——}} k$<br>$\lvert A\rvert = 1 \cdot 2$<br>$24\times$ |
| $i \overset{1}{\text{——}} j \overset{2}{\text{——}} k \overset{1}{\text{——}} i$<br>$\lvert A\rvert = 1 \cdot 2$<br>$12\times$ | $i \overset{1}{\text{——}} j \overset{2}{\text{——}} i \overset{1}{\text{——}} j$<br>$\lvert A\rvert = 1 \cdot 2$<br>$6\times$ |
| $i \overset{1}{\text{——}} j \overset{1}{\text{——}} k \overset{1}{\text{——}} l \overset{1}{\text{——}} i$<br>$\lvert A\rvert = 1$<br>$12\times$ | $i \overset{1}{\text{——}} j \overset{1}{\text{——}} k \overset{1}{\text{——}} l \overset{1}{\text{——}} j$<br>$\lvert A\rvert = 1$<br>$24\times$ |
| $i \overset{1}{\text{——}} j \overset{1}{\text{——}} k \overset{1}{\text{——}} l \overset{1}{\text{——}} k$<br>$\lvert A\rvert = 1$<br>$24\times$ | $i \overset{1}{\text{——}} j \overset{1}{\text{——}} k \overset{1}{\text{——}} i \overset{1}{\text{——}} j$<br>$\lvert A\rvert = 1$<br>$12\times$ |
| $i \overset{1}{\text{——}} j \overset{1}{\text{——}} k \overset{1}{\text{——}} j \overset{1}{\text{——}} i$<br>$\lvert A\rvert = 2$<br>$24\times$ | $i \overset{1}{\text{——}} j \overset{1}{\text{——}} k \overset{1}{\text{——}} j \overset{1}{\text{——}} k$<br>$\lvert A\rvert = 1$<br>$24\times$ |
| $i \overset{1}{\text{——}} j \overset{1}{\text{——}} k \overset{1}{\text{——}} i \overset{1}{\text{——}} k$<br>$\lvert A\rvert = 1$<br>$24\times$ | $i \overset{1}{\text{——}} j \overset{1}{\text{——}} i \overset{1}{\text{——}} j \overset{1}{\text{——}} i$<br>$\lvert A\rvert = 2$<br>$12\times$ |
| $i \overset{1}{\text{——}} j \overset{1}{\text{——}} i \overset{1}{\text{——}} k \overset{1}{\text{——}} i$<br>$\lvert A\rvert = 1$<br>$12\times$ | |

| | | | |
|---|---|---|---|
| 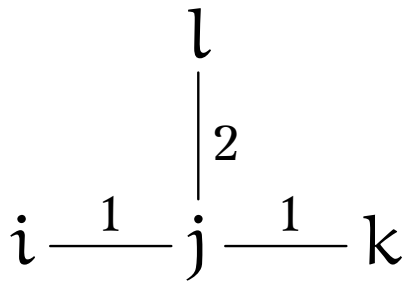<br>$\lvert A\rvert = 2$<br>$12\times$ | 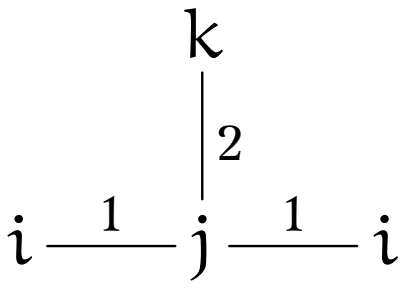<br>$\lvert A\rvert = 2 \cdot 2$<br>$24\times$ | 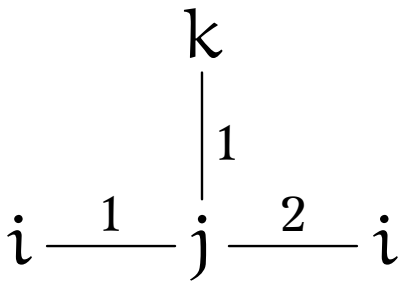<br>$\lvert A\rvert = 2$<br>$24\times$ | 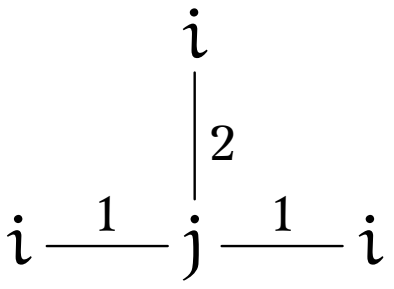<br>$\lvert A\rvert = 2 \cdot 2$<br>$12\times$ |
| 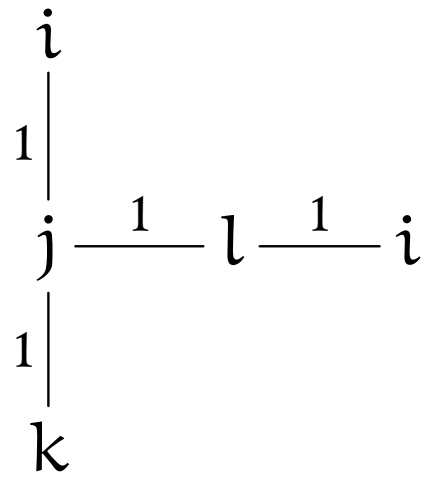<br>$\lvert A\rvert = 1$<br>$24\times$ | 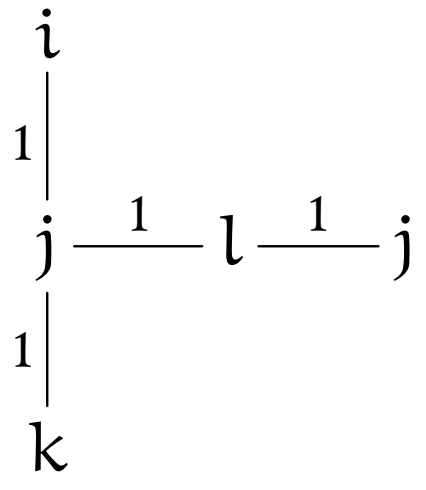<br>$\lvert A\rvert = 1$<br>$12\times$ | 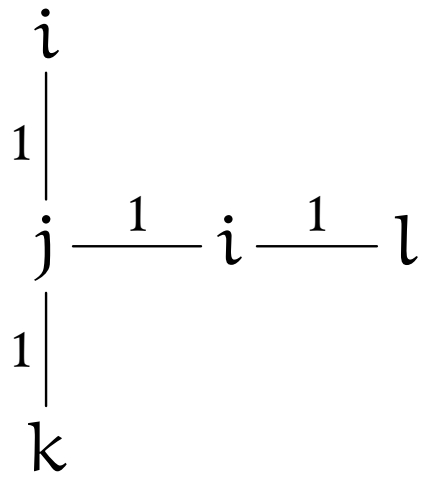<br>$\lvert A\rvert = 1$<br>$24\times$ | 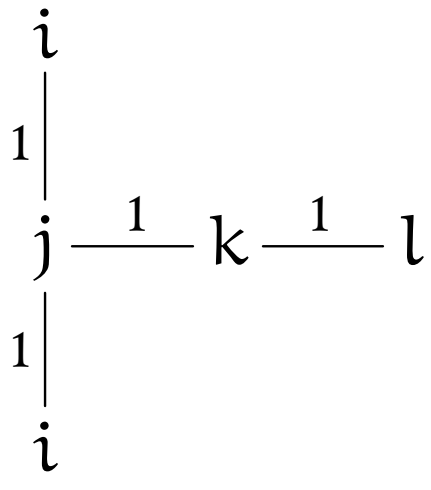<br>$\lvert A\rvert = 2$<br>$24\times$ |
| 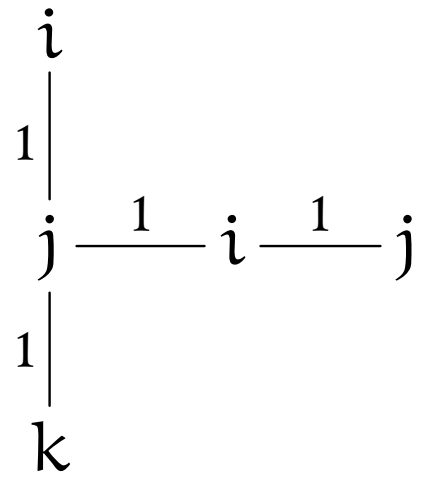<br>$\lvert A\rvert = 1$<br>$24\times$ | 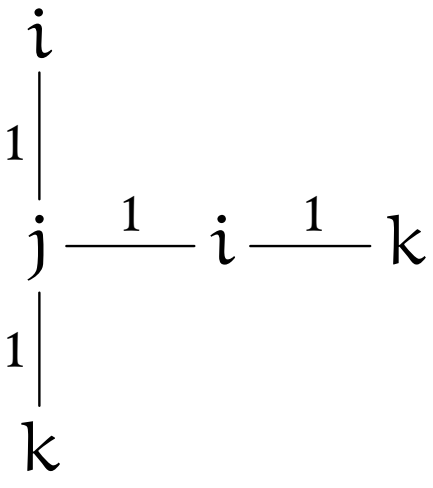<br>$\lvert A\rvert = 1$<br>$24\times$ | 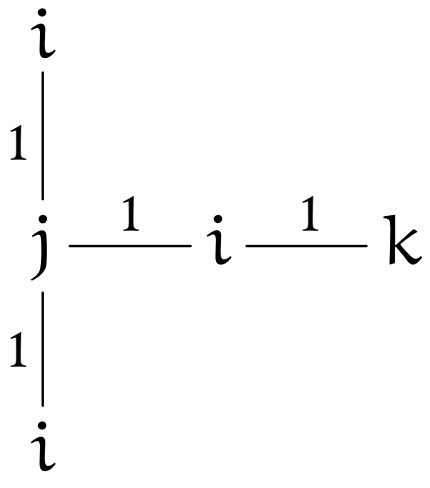<br>$\lvert A\rvert = 2$<br>$24\times$ | 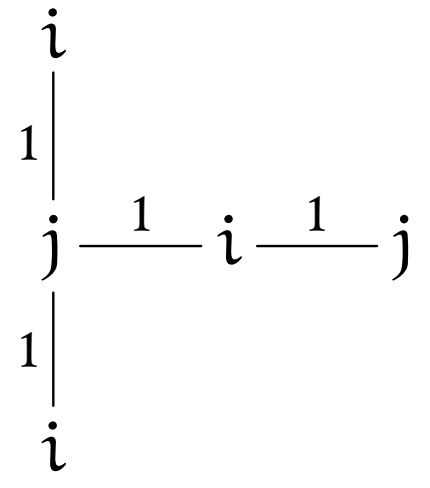<br>$\lvert A\rvert = 2$<br>$12\times$ |
| 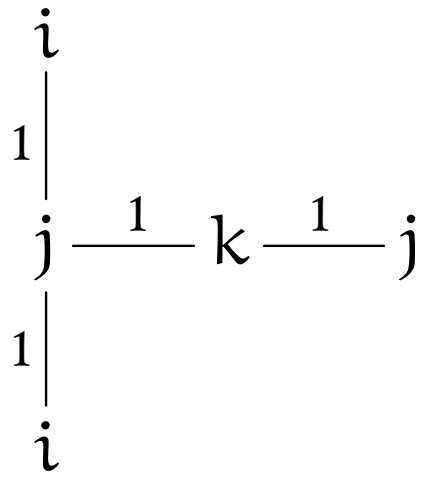<br>$\lvert A\rvert = 2$<br>$24\times$ | 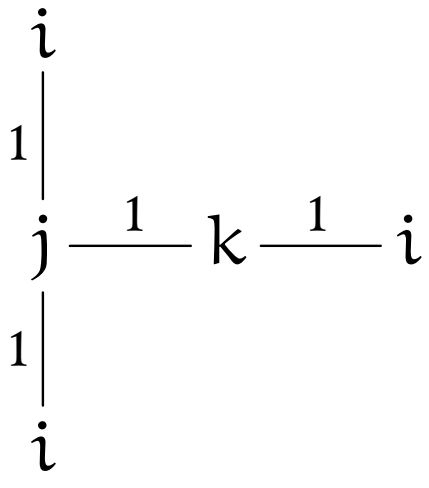<br>$\lvert A\rvert = 2$<br>$24\times$ | | |
| 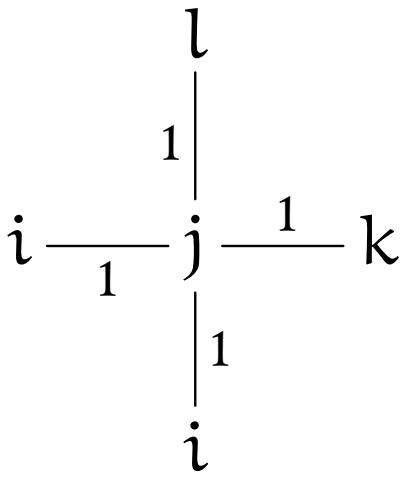<br>$\lvert A\rvert = 2$<br>$12\times$ | 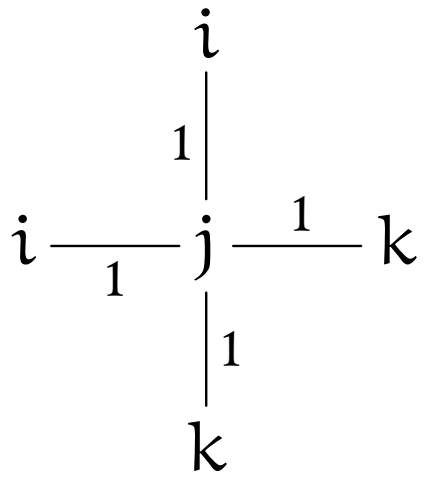<br>$\lvert A\rvert = 4$<br>$12\times$ | 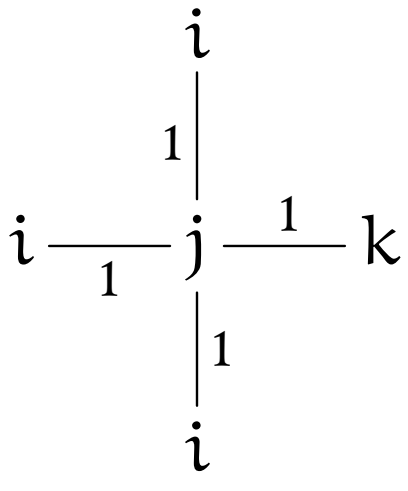<br>$\lvert A\rvert = 6$<br>$24\times$ | 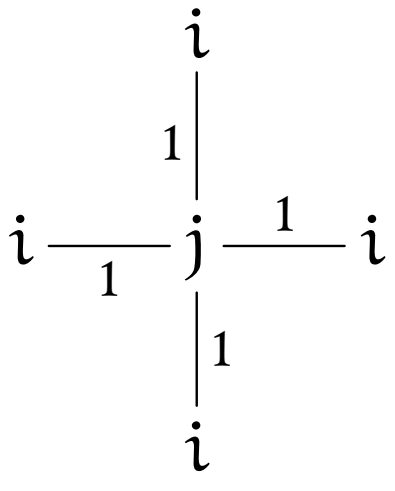<br>$\lvert A\rvert = 24$<br>$12\times$ |

# Apêndice E

# Algoritmos

## E.1 Dimensão e componentes

Os códigos abaixo foram executado em SINGULAR e usando procedimentos próprios compilados no arquivo myprocs[1] (disponível em `http://www.mat.ufmg.br/~israel/Projetos/myprocs.ses`).

O algoritmo abaixo constrói o ideal u em $W_d$ correspondente às parametrizações de curvas racionais de grau d que representam curvas de contato, pela caracterização legendriana (2.3). Também, calcula a dimensão da variedade gerada por esse ideal. Foi possível obter o valor das dimensões esperadas nos casos $d = 1, \cdots, 5$, os quais são iguais a dimensão esperada $2d+1$. Para casos superiores, o tempo de execução e uso da memória da máquina foram excedidos.

```
/* declaração de variáveis
int d=3;
/* a: coeficientes dos polinômios;
/* s,t: parâmetros das curvas;
/* x: variáveis no espaço projetivo;
ring r=0,(a(0..3)(0..d),s,t,x(0..3)),dp;
def aa=ideal(a(0..3)(0..d));
def xx=ideal(x(0..3));
poly f(0..3);
def td = ideal(1,t)^d;
forr(0,3,"f(i)=dotprod([a(i)(0..d)],td)");

/* definição do ideal das curvas de contato pela caracterização legend
def p = diff(f(1),t)*f(0) - f(1)*diff(f(0),t) -
  diff(f(2),t)*f(3) + f(2)*diff(f(3),t);
```

[1] incluem os procedimentos dotprod, um, putsolvform

```
def u = ideal(coeffs(p,t));
/* fixar a(0)(0) = 1
u = um(a(0)(0),u);
/* cálculo da dimensão
dimstd(u);
```

Em continuidade ao algoritmo anterior, o seguinte trecho de código identifica as componentes primárias/primas do ideal u. Foi possível obter a decomposição nos casos $d = 1, \cdots, 3$. Para grau 3 apresenta 2 componentes irredutíveis e reduzidas correspondentes as curvas de contato representadas por mapas de domínio $\mathbb{P}^1$ não ramificados e de recobrimento de retas em grau d, respectivamente.

```
/* identificando as equações lineares e realizando a
/* decomposição primária do ideal
def u0= putsolvform(u);
def l = primdecGTZ(dosubs(u0,u));
```

## E.2 Grafos e invariantes

Baseado na lista de tipos combinatórios de grafos do apêndice D, geramos todos os grafos correspondentes às componentes de pontos fixos de $\overline{\mathcal{M}}_{0,0}(\mathbb{P}^3, d)$, $d = 1, \cdots, 4$, em MAPLE. Para tal, foi usado o pacote GraphTheory.

A geração de grafos e colorações está longe de ser eficiente ou geral, uma vez que cada item é listado um a um, caso a caso. Por esse motivo, os cálculos cobrem apenas os casos $d = 1, \cdots, 4$. Ainda assim, o tempo de execução nesses casos foi irrisório, da ordem de poucos segundos.

Os procedimentos abaixo criam uma estrutura de dois dados: a primeira é o grafo dado por seu conjunto ordenado de vértices e sua matriz de adjacência (padrão do pacote); a segunda são as cores do grafo sob seus vértices ordenados.

Caso $d = 3$

```
> G1 := proc (i, j) options operator, arrow; [Graph([1, 2],
Matrix([[0, 3], [3, 0]]), weighted), [i, j]] end proc;
> G2 := proc (i, j, k) options operator, arrow; [Graph([1, 2, 3],
Matrix([[0, 1, 0], [1, 0, 2], [0, 2, 0]]), weighted),[i, j, k]] end proc;
> G3 := proc (i, j, k, l) options operator, arrow; [Graph([1, 2, 3, 4],
```

```
Matrix([[0, 1, 0, 0], [1, 0, 1, 0], [0, 1, 0, 1], [0, 0, 1, 0]]), weig
[i, j, k, l]] end proc;
> G4 := proc (i, j, k, l) options operator, arrow; [Graph([1, 2, 3, 4]
Matrix([[0, 1, 1, 1], [1, 0, 0, 0], [1, 0, 0, 0], [1, 0, 0, 0]]), weig
[i, j, k, l]] end proc;
```

Caso $d = 4$

```
> Gl1 := proc (i, j) options operator, arrow; [Graph([1, 2],
Matrix([[0, 4], [4, 0]]), weighted), [i, j]] end proc;
> Gl2a := proc (i, j, k) options operator, arrow; [Graph([1, 2, 3],
Matrix([[0, 3, 0], [3, 0, 1], [0, 1, 0]]), weighted), [i, j, k]] end p
> Gl2b := proc (i, j, k) options operator, arrow; [Graph([1, 2, 3],
Matrix([[0, 2, 0], [2, 0, 2], [0, 2, 0]]), weighted), [i, j, k]] end p
> Gl3a := proc (i, j, k, l) options operator, arrow; [Graph([1, 2, 3,
 Matrix([[0, 2, 0, 0], [2, 0, 1, 0], [0, 1, 0, 1], [0, 0, 1, 0]]), wei
[i, j, k, l]] end proc;
> Gl3b := proc (i, j, k, l) options operator, arrow; [Graph([1, 2, 3,
Matrix([[0, 1, 0, 0], [1, 0, 2, 0], [0, 2, 0, 1], [0, 0, 1, 0]]), weig
[i, j, k, l]] end proc;
> Gl4 := proc (i, j, k, l, m) options operator, arrow; [Graph([1, 2, 3
 Matrix([[0, 1, 0, 0, 0], [1, 0, 1, 0, 0], [0, 1, 0, 1, 0], [0, 0, 1,
[0, 0, 0, 1, 0]]), weighted), [i, j, k, l, m]] end proc;
> Ge3 := proc (i, j, k, l) options operator, arrow; [Graph([1, 2, 3, 4
Matrix([[0, 2, 1, 1], [2, 0, 0, 0], [1, 0, 0, 0], [1, 0, 0, 0]]), weig
[i, j, k, l]] end proc;
> Ge4 := proc (i, j, k, l, m) options operator, arrow; [Graph([1, 2, 3
Matrix([[0, 1, 1, 1, 1], [1, 0, 0, 0, 0], [1, 0, 0, 0, 0], [1, 0, 0, 0
[1, 0, 0, 0, 0]]), weighted), [i, j, k, l, m]] end proc;
> Gt4 := proc (i, j, k, l, m) options operator, arrow; [Graph([1, 2, 3
Matrix([[0, 1, 1, 1, 0], [1, 0, 0, 0, 0], [1, 0, 0, 0, 0], [1, 0, 0, 0
[0, 0, 0, 1, 0]]), weighted), [i, j, k, l, m]] end proc;
```

A seguir, os procedimentos para cálculo dos fatores equivariantes da fórmula de localização de $N_d$, $d \geq 1$: T calcula, de forma auxiliar, a classe equivariante do tangente; E o fator correspondente às arestas de $\Gamma$ e V aos vértices; H calcula a classe de incidência a retas; Id a classe de condição de contato; P calcula o produto dos fatores citados.

```
# Cores/rótulos
> S := {0, 1, 2, 3}
```

```
> T := proc (i) #Chern top do tangente ao i-ésimo ponto fixo de P^3
    local p, j;
  p := 1;
  for j in `minus`(S, {i}) do
    p := p*(x[i]-x[j])
  end do
  end proc;

> E := proc (G) #Contribuição das arestas para Bott
    local e, k, a, p, q, d, h, t;
p := 1;
for e in Edges(G[1], weights) do
  d := e[2];
h := G[2][e[1][1]];
t := G[2][e[1][2]];
q := 1;
for k in `minus`(S, {h, t}) do
  for a from 0 to d do
q := q/((a*x[h]+(d-a)*x[t])/d-x[k])
end do
end do;
  p := p*q*(-1)^d*(d/(x[h]-x[t]))^(2*d)/factorial(d)^2
end do
end proc;

> V := proc (G) #Contribuição dos vértices para Bott
    local p, q, s, v, n, h, F, e, t, d;
p := 1;
for v in Vertices(G[1]) do
    n := Degree(G[1], v);
  h := G[2][v];
  F := IncidentEdges(G[1], v);
  q := 1;
  s := 0;
  for e in F do
    t := G[2][op(`minus`(e, {v}))];
  d := GetEdgeWeight(G[1], e);
  q := q*d/(x[h]-x[t]);
```

```
  s := s+d/(x[h]-x[t])
  end do;
  p := p*T(h)^(n-1)*s^(n-3)*q
  end do
end proc;

> H := proc (G) #peso do divisor de incidência à reta
    local s, v, h, e, n;
  s := 0;
  for v in Vertices(G[1]) do
    h := G[2][v];
  n := 0;
  for e in IncidentEdges(G[1], v) do
    n := n+GetEdgeWeight(G[1], e)
  end do;
  s := s+n*x[h]
  end do
  end proc;

> Id := proc (G) #Contribuição da condição de contato
    local p, e, d, h, t, l, v, n, m;
p := 1;
for e in Edges(G[1], weights) do
  d := e[2];
h := G[2][e[1][1]];
t:= G[2][e[1][2]];
for l to 2*d-1 do
  p := p*(l*x[h]+(2*d-l)*x[t])/d
end do
end do;
for v in Vertices(G[1]) do
  n := Degree(G[1], v)-1;
m := G[2][v];
p := p*(2*x[m])^n
end do
end proc;

> P := proc (G) #Produto de todas contribuições
    E(G)*V(G)*H(G)^(2d+1)*Id(G)
```

```
end proc;
```

A seguir, P é calculado sobre cada coloração dos tipos combinatórios gerados, levando em consideração, caso a caso, a ordem do grupo de automorfismo do grafo. Tais colorações são listadas uma a uma.

Caso $d = 3$

```
> R1 := 0; for i from 0 to 3 do for j from i+1 to 3 do
R1 := R1+(1/3)*P(G1(i, j)) end do end do;
> R2a := 0; for i from 0 to 3 do for j in ‘minus‘(S, {i}) do
R2a := R2a+(1/2)*P(G2(i, j, i)) end do end do;
> R2b := 0; for i from 0 to 3 do for j in ‘minus‘(S, {i}) do
for k in ‘minus‘(S, {i, j}) do
R2b := R2b+(1/2)*P(G2(i, j, k)) end do end do end do;
> R3a := 0; for i from 0 to 3 do for j from i+1 to 3 do
R3a := R3a+P(G3(i, j, i, j)) end do end do;
> R3b := 0; for i from 0 to 3 do for j in ‘minus‘(S, {i}) do
for k in ‘minus‘(S, {i, j}) do
R3b := R3b+P(G3(i, j, i, k)) end do end do end do;
> R3c := 0; for i from 0 to 3 do for j from i+1 to 3 do
for k in ‘minus‘(S, {i, j}) do
R3c := R3c+P(G3(k, i, j, k)) end do end do end do;
> R3d := 0; for i from 0 to 3 do for j in ‘minus‘(S, {i}) do
for k in ‘minus‘(S, {i, j}) do l := op(‘minus‘(S, {i, j, k}));
R3d := R3d+(1/2)*P(G3(i, j, k, l)) end do end do end do;
> R4a := 0; for i from 0 to 3 do for j in ‘minus‘(S, {i}) do
R4a := R4a+(1/6)*P(G4(i, j, j, j)) end do end do;
> R4b := 0; for i from 0 to 3 do for j in ‘minus‘(S, {i}) do
for k in ‘minus‘(S, {i, j}) do
R4b := R4b+(1/2)*P(G4(i, j, k, k)) end do end do end do;
> R4c := P(G4(0, 1, 2, 3))+P(G4(1, 0, 2, 3))+
P(G4(2, 0, 1, 3))+P(G4(3, 0, 1, 2));
> R := R1+R2a+R2b+R3a+R3b+R3c+R3d+R4a+R4b+R4c;
```

Caso $d = 4$

```
> R1 := 0; for i from 0 to 3 do for j from i+1 to 3 do
R1 := R1+(1/4)*P(Gl1(i, j)) end do end do;
> R2a := 0; for i from 0 to 3 do for j in ‘minus‘(S, {i}) do
R2a := R2a+(1/3)*P(Gl2a(i, j, i)) end do end do;
```

```
> R2b := 0; for i from 0 to 3 do for j in ‘minus‘(S, {i}) do
for k in ‘minus‘(S, {i, j}) do
R2b := R2b+(1/3)*P(Gl2a(i, j, k)) end do end do end do;
> R2c := 0; for i from 0 to 3 do for j in ‘minus‘(S, {i}) do
R2c := R2c+(1/8)*P(Gl2b(i, j, i)) end do end do;
> R2d := 0; for j from 0 to 3 do for l in ‘minus‘(S, {j}) do
i, k := op(‘minus‘(S, {j, l}));
R2d := R2d+(1/4)*P(Gl2b(i, j, k)) end do end do;
> R3a := 0; for i from 0 to 3 do for j in ‘minus‘(S, {i}) do
R3a := R3a+(1/2)*P(Gl3a(i, j, i, j)) end do end do;
> R3b := 0; for i from 0 to 3 do for j in ‘minus‘(S, {i}) do
for k in ‘minus‘(S, {i, j}) do
R3b := R3b+(1/2)*P(Gl3a(i, j, k, i)) end do end do end do;
> R3c := 0; for i from 0 to 3 do for j in ‘minus‘(S, {i}) do
for k in ‘minus‘(S, {i, j}) do
R3c := R3c+(1/2)*P(Gl3a(i, j, k, j)) end do end do end do;
> R3d := 0; for i from 0 to 3 do for j in ‘minus‘(S, {i}) do
for k in ‘minus‘(S, {i, j}) do
R3d := R3d+(1/2)*P(Gl3a(i, j, i, k)) end do end do end do;
> R3e := 0; for i from 0 to 3 do for j in ‘minus‘(S, {i}) do
for k in ‘minus‘(S, {i, j}) do l := op(‘minus‘(S, {i, j, k}));
R3e := R3e+(1/2)*P(Gl3a(i, j, k, l)) end do end do end do;
> R3f := 0; for i from 0 to 3 do for j from i+1 to 3 do
R3f := R3f+(1/2)*P(Gl3b(i, j, i, j)) end do end do;
> R3g := 0; for i from 0 to 3 do for l in ‘minus‘(S, {i}) do
j, k := op(‘minus‘(S, {i, l}));
R3g := R3g+(1/2)*P(Gl3b(i, j, k, i)) end do end do;
> R3h := 0; for i from 0 to 3 do for j in ‘minus‘(S, {i}) do
for k in ‘minus‘(S, {i, j}) do
R3h := R3h+(1/2)*P(Gl3b(i, j, i, k)) end do end do end do;
> R3m := 0; for i from 0 to 3 do for l from i+1 to 3 do
for j in ‘minus‘(S, {i, l}) do k := op(‘minus‘(S, {i, j, l}));
R3m := R3m+(1/2)*P(Gl3b(i, j, k, l)) end do end do end do;
> R4a := 0; for i from 0 to 3 do for j in ‘minus‘(S, {i}) do
R4a := R4a+(1/2)*P(Gl4(i, j, i, j, i)) end do end do;
> R4b := 0; for i from 0 to 3 do for j in ‘minus‘(S, {i}) do
for k in ‘minus‘(S, {i, j}) do
R4b := R4b+P(Gl4(i, j, k, i, k)) end do end do end do;
> R4c := 0; for i from 0 to 3 do for j in ‘minus‘(S, {i}) do
```

```
for k in `minus`(S, {i, j}) do
R4c := R4c+P(Gl4(i, j, k, j, k)) end do end do end do;
> R4d := 0; for i from 0 to 3 do for j in `minus`(S, {i}) do
for k in `minus`(S, {i, j}) do
R4d := R4d+(1/2)*P(Gl4(i, j, k, j, i)) end do end do end do;
> R4e := 0; for k from 0 to 3 do for l in `minus`(S, {k}) do
i, j := op(`minus`(S, {k, l}));
R4e := R4e+P(Gl4(i, j, k, i, j)) end do end do;
> R4f := 0; for i from 0 to 3 do for j in `minus`(S, {i}) do
for k in `minus`(S, {i, j}) do l := op(`minus`(S, {i, j, k}));
R4f := R4f+P(Gl4(i, j, k, l, k)) end do end do end do;
> R4g := 0; for i from 0 to 3 do for j in `minus`(S, {i}) do
for k in `minus`(S, {i, j}) do l := op(`minus`(S, {i, j, k}));
R4g := R4g+P(Gl4(i, j, k, l, j)) end do end do end do;
> R4h := 0; for i from 0 to 3 do for k in `minus`(S, {i}) do
j, l := op(`minus`(S, {i, k}));
R4h := R4h+P(Gl4(i, j, k, l, i)) end do end do;
> R4n := 0; for i from 0 to 3 do for k in `minus`(S, {i}) do
j, l := op(`minus`(S, {i, k}));
R4n := R4n+P(Gl4(j, k, i, k, l)) end do end do;
> R4m := 0; for i from 0 to 3 do for l in `minus`(S, {i}) do
j, k := op(`minus`(S, {i, l}));
R4m := R4m+P(Gl4(i, j, i, k, i)) end do end do;
> R5a := 0; for i from 0 to 3 do for j in `minus`(S, {i}) do
R5a := R5a+(1/4)*P(Ge3(i, j, j, j)) end do end do;
> R5b := 0; for i from 0 to 3 do for j in `minus`(S, {i}) do
for k in `minus`(S, {i, j}) do
R5b := R5b+(1/4)*P(Ge3(i, j, k, k)) end do end do end do;
> R5c := 0; for i from 0 to 3 do for j in `minus`(S, {i}) do
for k in `minus`(S, {i, j}) do
R5c := R5c+(1/2)*P(Ge3(i, j, j, k)) end do end do end do;
> R5d := 0; for i from 0 to 3 do for j in `minus`(S, {i}) do
k, l := op(`minus`(S, {i, j}));
R5d := R5d+(1/2)*P(Ge3(i, j, k, l)) end do end do;
> R6a := 0; for i from 0 to 3 do for j in `minus`(S, {i}) do
R6a := R6a+(1/24)*P(Ge4(i, j, j, j, j)) end do end do;
> R6b := 0; for i from 0 to 3 do for j in `minus`(S, {i}) do
for k in `minus`(S, {i, j}) do
R6b := R6b+(1/6)*P(Ge4(i, j, j, j, k)) end do end do end do;
```

```
> R6c := 0; for i from 0 to 3 do for l in ‘minus‘(S, {i}) do
j, k := op(‘minus‘(S, {i, l}));
R6c := R6c+(1/4)*P(Ge4(i, j, j, k, k)) end do end do;
> R6d := 0; for i from 0 to 3 do for j in ‘minus‘(S, {i}) do
k, l := op(‘minus‘(S, {i, j}));
R6d := R6d+(1/2)*P(Ge4(i, j, j, k, l)) end do end do;
> R7a := 0; for i from 0 to 3 do for j in ‘minus‘(S, {i}) do
R7a := R7a+(1/2)*P(Gt4(i, j, j, j, i)) end do end do;
> R7b := 0; for i from 0 to 3 do for j in ‘minus‘(S, {i}) do
for k in ‘minus‘(S, {i, j}) do
R7b := R7b+(1/2)*P(Gt4(i, j, j, j, k)) end do end do end do;
> R7c := 0; for i from 0 to 3 do for j in ‘minus‘(S, {i}) do
for k in ‘minus‘(S, {i, j}) do
R7c := R7c+P(Gt4(i, j, k, j, k)) end do end do end do;
> R7d := 0; for i from 0 to 3 do for j in ‘minus‘(S, {i}) do
for k in ‘minus‘(S, {i, j}) do
R7d := R7d+P(Gt4(i, j, k, j, i)) end do end do end do;
> R7e := 0; for i from 0 to 3 do for j in ‘minus‘(S, {i}) do
for k in ‘minus‘(S, {i, j}) do l := op(‘minus‘(S, {i, j, k}));
R7e := R7e+(1/2)*P(Gt4(i, j, j, k, l)) end do end do end do;
> R7f := 0; for i from 0 to 3 do for l in ‘minus‘(S, {i}) do
j, k := op(‘minus‘(S, {i, l}));
R7f := R7f+P(Gt4(i, j, k, l, i)) end do end do;
> R7g := 0; for i from 0 to 3 do for k in ‘minus‘(S, {i}) do
for l in ‘minus‘(S, {i, k}) do j := op(‘minus‘(S, {i, k, l}));
R7g := R7g+P(Gt4(i, j, k, l, j)) end do end do end do;
> R7h := 0; for i from 0 to 3 do for k in ‘minus‘(S, {i}) do
for l in ‘minus‘(S, {i, k}) do j := op(‘minus‘(S, {i, k, l}));
R7h := R7h+P(Gt4(i, j, k, j, l)) end do end do end do;
> R7m := 0; for i from 0 to 3 do for j in ‘minus‘(S, {i}) do
for k in ‘minus‘(S, {i, j}) do
R7m := R7m+(1/2)*P(Gt4(i, j, j, k, i)) end do end do end do;
> R7n := 0; for i from 0 to 3 do for j in ‘minus‘(S, {i}) do
for k in ‘minus‘(S, {i, j}) do
R7n := R7n+(1/2)*P(Gt4(i, j, j, k, j)) end do end do end do;
> R := R6d+R7a+R6c+R7d+R2a+R4e+R4f+R1+R6b+R4c+R7b+R7c+R6a+R5d+R5c+R4b+
R5a+R5b+R3f+R7m+R3g+R3h+R3m+R4a+R2b+R7h+R2c+R2d+R3a+R3b+R3c+R3d+R3e+
R4n+R4m+R7f+R7g+R4d+R7n+R4h+R4g+R7e;
```

# Referências Bibliográficas